\documentclass[draft]{amsart}

\usepackage[utf8]{inputenc}

\usepackage{fixme} % For annotations (contains verbatim)
\fxusetheme{color}

\usepackage{amssymb}
\usepackage{stackengine}[2013-09-11]
\usepackage{graphicx}
\usepackage{xcolor}
\usepackage{subcaption}
\usepackage{cancel}
\usepackage{multirow}
\usepackage{mathrsfs}
\usepackage{tikz-cd}
\usepackage{enumitem}
\usepackage{float}
\usepackage{setspace}
\usepackage{appendix}
\usepackage[colorlinks=true, final=true]{hyperref}

\newcommand{\Ss}{\mathcal{S}}

\newcommand{\Nn}{\mathcal{N}}

\newcommand{\jj}{\mathcal{J}}
\newcommand{\J}{\mathcal{J}}

\newcommand{\T}{\mathbb{T}^2}

\newcommand{\N}{\mathbb{N}}
\newcommand{\Z}{\mathbb{Z}}

\newcommand{\R}{\mathbb{R}}
\newcommand{\C}{\mathbb{C}}

\newcommand{\om}{\omega}

\newcommand{\joml}{\J_{\om_\lambda}}
\newcommand{\oml}{\omega_\lambda}

\newcommand{\jjoml}{\jj_{\om_\lambda}}
\newcommand{\jzoml}{\jj^{\Z_n}_{\om_\lambda}}

\newcommand{\jzomll}{\jj^{\Z_n}_{\om_\lambda,l}}

\newcommand{\del}{\partial}
\newcommand{\delbar}{\bar\partial}

\newcommand{\from}{\leftarrow}
\newcommand{\into}{\hookrightarrow}

\DeclareRobustCommand\longtwoheadrightarrow
     {\relbar\joinrel\twoheadrightarrow}

\newcommand{\SSS}{S^2 \times S^2}

\newcommand{\CP}{\mathbb{C}P}
\newcommand{\CCC}{\CP^2\# \overline{\CP^2}}
\newcommand{\jsoml}{\mathcal{J}^{S^1}_{\om_\lambda, l}}
\newcommand{\jsom}{\mathcal{J}^{S^1}_{\om_\lambda}}

\DeclareMathOperator{\Diff}{Diff}

\DeclareMathOperator{\Symp}{Symp}
\DeclareMathOperator{\Sym}{Sym}

\DeclareMathOperator{\Ham}{Ham}
\DeclareMathOperator{\Stab}{Stab}
\DeclareMathOperator{\Fix}{Fix}

\DeclareMathOperator{\Aut}{Aut}
\DeclareMathOperator{\Iso}{Iso}

\DeclareMathOperator{\AGL}{AGL}
\DeclareMathOperator{\Sp}{Sp}

\DeclareMathOperator{\U}{U}
\DeclareMathOperator{\SO}{SO}
\DeclareMathOperator{\PU}{PU}
\DeclareMathOperator{\PSL}{PSL}
\DeclareMathOperator{\SU}{SU}

\DeclareMathOperator{\pushout}{pushout}

\DeclareMathOperator{\id}{id}

\DeclareMathOperator{\Det}{Det}

\newtheorem{thm}{Theorem}[section]
\newtheorem{prop}[thm]{Proposition}
\newtheorem{cor}[thm]{Corollary}
\newtheorem{lemma}[thm]{Lemma}

\theoremstyle{definition}
\newtheorem{defn}[thm]{Definition}
\newtheorem{exmp}[thm]{Example}

\theoremstyle{remark}
\newtheorem{remark}[thm]{Remark}

\newcommand{\eoesymbol}{$\between$}

\DeclareRobustCommand{\eoe}{%
  \ifmmode \mathqed
  \else
    \leavevmode\unskip\penalty9999 \hbox{}\nobreak\hfill
    \quad\hbox{\eoesymbol}%
  \fi
}

\begin{document}

%%%%%%%%%%%%%%%%%%%%%%%%%%%%%%%%%%%%%%%%%%%%%%%%%%%%%%%%%%%%%%%%%%%%%%%%%%%%%%%%
% Frontmatter
%%%%%%%%%%%%%%%%%%%%%%%%%%%%%%%%%%%%%%%%%%%%%%%%%%%%%%%%%%%%%%%%%%%%%%%%%%%%%%%%

% Title
\title[Centralizers of Hamiltonian finite cyclic groups actions]{Centralizers of Hamiltonian finite cyclic group actions\\ on rational ruled surfaces}

% author one information
\author{Pranav V. Chakravarthy}
\address{Department of Mathematics\\ The Hebrew University of Jerusalem\\ Israel}
\email{pranav.chakravarthy@mail.huji.ac.il}
\thanks{}

% author two information
\author{Martin Pinsonnault}
\address{Department of Mathematics\\  University of Western Ontario\\ Canada}
\email{mpinson@uwo.ca}

\date{\today}

\thanks{This paper contains results from the first author's PhD thesis at the University of Western Ontario under the supervision of the second author. The first author would also like to thank the Hebrew University of Jerusalem where this project was completed. The second author is supported by NSERC grant RGPIN-2020-06428.}

\subjclass[2020]{Primary 53D35; Secondary 57R17,57S05,57T20}
\keywords{symplectic geometry, symplectomorphism group, fundamental group, Hamiltonian circle actions}

\begin{abstract}
Let $M=(M,\om)$ be either $S^2 \times S^2$ or $\CCC$ endowed with any symplectic form $\om$. Suppose a finite cyclic group $\Z_n$ is acting effectively on $(M,\om)$ through Hamiltonian diffeomorphisms, that is, there is an injective homomorphism $\Z_n\into \Ham(M,\om)$. In this paper, we investigate the homotopy type of the group $\Symp^{\Z_n}(M,\om)$ of equivariant symplectomorphisms. We prove that for some infinite families of $\Z_n$ actions satisfying certain inequalities involving the order $n$ and the symplectic cohomology class $[\om]$, the actions extend to either one or two toric actions, and accordingly, that the centralizers are homotopically equivalent to either a finite dimensional Lie group, or to the homotopy pushout of two tori along a circle. 
Our results rely on $J$-holomorphic techniques, on Delzant's classification of toric actions, on Karshon's classification of Hamiltonian circle actions on $4$-manifolds, and on the Chen-Wilczy\'nski classification of smooth $\Z_n$-actions on Hirzebruch surfaces.
\end{abstract}

\maketitle

%\textbf{\Large Version 2023-07-4}

\tableofcontents

%%%%%%%%%%%%%%%%%%%%%%%%%%%%%%%%%%%%%%%%%%%%%%%%%%%%%%%%%%%%%%%%%%%%%%%%%%%%%%%%
% Main content
%%%%%%%%%%%%%%%%%%%%%%%%%%%%%%%%%%%%%%%%%%%%%%%%%%%%%%%%%%%%%%%%%%%%%%%%%%%%%%%%

%%%%%%%%%%%%%%%%%%%%%%%%%%%%%%%%%%%%%%%%%%%%%%%%%%%%%%%%%%%%%%%%%%%%%%%%%%%%%%%%
\section{Introduction}
%%%%%%%%%%%%%%%%%%%%%%%%%%%%%%%%%%%%%%%%%%%%%%%%%%%%%%%%%%%%%%%%%%%%%%%%%%%%%%%%

Let $M=(M,\om)$ be either $\SSS$ or $\CCC$ endowed with any symplectic form $\om$, and let $\Symp(M,\om)$ be the group of symplectomorphisms of $M$ endowed with the $C^{\infty}$ topology. Suppose a finite cyclic group $\Z_n$ is acting effectively and symplectically on $(M,\om)$. We say that such an action is Hamiltonian if the associated injective homomorphism $\rho:\Z_n\into\Symp(M,\om)$ takes values in the subgroup $\Ham(M,\om)\subset\Symp(M,\om)$ of Hamiltonian diffeomorphisms. The main goal of this paper is to determine the homotopy type of the centralizer subgroup $C(\Z_{n})=\Symp^{\Z_n}(M,\om)\subset\Symp(M,\om)$ consisting of symplectomorphisms that commute with such an action. 

One of the main difficulties in determining the homotopy type of centralizers subgroups $\Symp^{\Z_n}(M,\om)$ is that, at the moment of writing, there is no classification of Hamiltonian $\Z_{n}$ actions on $\SSS$ and $\CCC$ up to equivariant symplectomorphisms. This forces us to restrict ourselves to $\Z_n$ actions satisfying some ad hoc inequalities involving the symplectic cohomology class $[\om]$ and the order of the cyclic group $\Z_n$. These inequalities ensure that the action of the centralizer $\Symp^{\Z_n}(M,\om)$ on the contractible space $\jj^{\Z_n}(M,\om)$ of invariant and compatible almost complex structures has a simple orbit structure that can be precisely described despite our lacunary knowledge of Hamiltonian $\Z_n$ actions. As we will see, this is enough to deal with large families of Hamiltonian $\Z_n$ actions, leaving open the cases of actions satisfying some specific numerical conditions.

\section{Overview of the argument} To simplify the exposition, we first explain the structure of our argument in the case of Hamiltonian $\Z_n$ actions on the product $\SSS$. The argument for the non-trivial bundle $\CCC$, which is slightly simpler, follows the same lines and will be addressed in Section~\ref{section:CCC}. We closely follow the $J$-holomorphic framework developed in our previous work~\cite{ChPin-Memoirs} on centralizer subgroups of Hamiltonian circle actions on $\SSS$ and $\CCC$. In order to make the present paper reasonably short, and to avoid unnecessary duplication, we make frequent references to~\cite{ChPin-Memoirs} and we invite the reader to consult this paper for further details.

%%%%%%%%%%%%%%%%%%%%%%%%%%%%%%%%%%%%%%%%%%%%%%%%%%%%%%%%%%%%%%%%%%%%%%%%%%%%%%%%
\subsection{Abelian actions on $(\SSS,\oml)$}
%%%%%%%%%%%%%%%%%%%%%%%%%%%%%%%%%%%%%%%%%%%%%%%%%%%%%%%%%%%%%%%%%%%%%%%%%%%%%%%%
Our first task is to describe all possible $\Z_n$ actions on $\SSS$. To this end, recall that the Lalonde-McDuff classification of symplectically ruled $4$-manifolds~\cite{MR1426534} implies that any symplectic form on $\SSS$ is conformally equivalent to a product form $\oml:=\lambda\sigma\times\sigma$, where $\sigma$ is the standard area form on $S^2$ with total area $1$, and where $\lambda\geq1$. Let's write $(M,\oml)$ for $\SSS$ endowed with such a split form. Suppose a finite cyclic group $\Z_n$ is acting effectively on $(M,\oml)$ through Hamiltonian diffeomorphisms, that is, suppose there is an injective homomorphism $\Z_n\into \Ham(M,\oml)$. By~\cite{Liat}~Theorem~1.1, any such action extends to a Hamiltonian circle action, and we know from~\cite{finTori} that any circle action extends to a toric action. Note that these extensions are generally not unique: up to equivalence, a cyclic action may extend to more than one circle action and, similarly, a circle action may extend to different toric actions. Note also that, a priori, a cyclic action $\Z_n$ may extend simultaneously to some circle action $S^1$ and also to some toric action $\T$, but that the circle action itself may not extend to the same toric action (see Remark~\ref{Rmk:ATFcyclic action}).

We say that two symplectic actions $\rho_{1},\rho_{2}:G\to\Symp(M,\om)$ are \emph{equivariantly symplectomorphic} if there exists a symplectomorphism $\phi$ such that $\rho_2=\phi\circ\rho_1\circ\phi^{-1}$. We say that the two actions are \emph{equivalent} if the subgroups $\rho_{1}(G)$ and $\rho_{2}(G)$ belong to the same $\Symp(M,\om)$-conjugacy class or, equivalently, if there exists an automorphism $\Gamma\in\Aut(G)$ and a symplectomorphism $\phi$ such that $\rho_{2}\circ \Gamma=\phi\circ\rho_{1}\circ\phi^{-1}$. By definition, the subgroup $\Symp^{\Z_n}(M,\om)$ is the centralizer of $\rho(\Z_{n})$ in $\Symp(M,\om)$ and, as a topological group, its homeomorphism type only depends on the equivalence class of the action $\rho:\Z_n\to\Symp(M,\om)$.

Let's write $\lambda=\ell+\delta$ with $\ell$ an integer and $0<\delta\leq 1$. As explained in~\cite{ChPin-Memoirs}~Section~2.2.2, it follows from Delzant's classification that any toric action on $(M,\oml)$ is equivalent to exactly one of the even Hirzebruch actions $\T_{2k}$ with $0\leq k\leq \ell$. Consequently, writing $r=2k$ and fixing an identification $\T_r\simeq S^1\times S^1$, any circle action is equivalent to the restriction of a toric action to a subcircle
\[S^1(a,b;r)=\{(t^a,t^b)~|~t\in S^1\}\subset \T_r\]
where $a,b\in\Z$. Since we only consider effective actions, we can further assume that $\gcd(a,b)=1$. Let $\Z_n\subset S^1$ be the standard inclusion and let $\zeta= e^{2\pi i/n}$ be the first $n^{\text{th}}$ root of unity. By Theorem~1.1 in~\cite{Liat}, any effective Hamiltonian $\Z_n$ action is the restriction to $\Z_n$ of some $S^1(a,b;r)$ action for which $\gcd(a,b)=1$. More concretely, the action $\Z_n(a,b;r)$ is defined via the inclusion \begin{align*}
     \Z_n &\hookrightarrow  \Z_n(a,b;r)\subset\T_r \\
    \zeta &\mapsto (\zeta^a, \zeta^b)
\end{align*}

Note that changing the preferred generator of $\Z_n$ or, equivalently, reparametrizing the $\Z_n$ subgroup, changes the action but determines the same finite subgroup inside $\Symp(M,\om)$, hence the same centralizer subgroup $\Symp^{\Z_n}(M,\om)$.  

\begin{defn}
We say the cyclic action $\Z_n(a,b;r)$ extends to the circle action $S^1(a',b';r')$ iff there exists a cyclic action $\Z_n(a',b';r')$ and a equivariant symplectomorphism $\phi$ that intertwines the $\Z_n(a,b;r)$ and $\Z_n(a',b';r')$ actions.  
We shall denote this by $\Z_n(a,b;r) \subset S^1(a',b';r')$. 
\end{defn}

Given two circles $S^1(a,b;r)$ and $S^1(a',b';r)$ in the same torus $\T_r$, their intersection \\
$S^1(a,b;r) \cap S^1(a',b';r)$ contains a finite cyclic group of order $n$ iff there exists $k\in \mathbb{N}$ such that
\begin{enumerate}
\item $\gcd(k,n)=1$ and 
\item $a \equiv a'k\pmod{n}$ $b \equiv b'k\pmod{n}$
\end{enumerate}
The corresponding actions $\Z_n(a,b;r)$ and $\Z_n(a',b';r)$ differ by the reparametrization $\zeta\mapsto \zeta^k$ and we have equality $\Z_n(a,b;r) = S^1(a,b;r) \cap S^1(a',b';r) =\Z(a',b';r)$ iff $a \equiv a' \pmod{n}$ and $b \equiv b' \pmod{n}$. 
\\

For convenience, we summarize the previous discussion in a single proposition.

\begin{prop}\label{prop:SummaryActionsS2xS2}
Let $M_{\lambda}$ be the product $\SSS$ endowed with the split symplectic form $\oml=\lambda\sigma\times\sigma$, where $1\leq\lambda =\ell+\delta$ with $\ell\geq0$ an integer and $0<\delta\leq 1$.
\begin{enumerate}
\item Up to rescaling, every symplectic form on $\SSS$ is diffeomorphic to a standard product form $\oml$ with $\lambda\geq 1$.
\item Every toric action on $M_{\lambda}$ is equivalent to one of the $(\ell+1)$ standard Hirzebruch actions $\T_{2k}$, where $0\leq k\leq \ell$.
\item Every Hamiltonian circle action on $M_{\lambda}$ is equivalent to the restriction of some toric action $\T_{2k}$ to a subcircle $S^{1}(a,b;2k)$.
\item Similarly, every Hamiltonian action of a cyclic group $\Z_{n}$ on $M_{\lambda}$ is equivalent to the restriction of some toric action $\T_{2k}$ to a subgroup $\Z_{n}(a,b;2k)\subset S^{1}(a,b;2k)\subset\T_{2k}$.\qed 
\end{enumerate}
\end{prop}

\begin{remark}
Although we will not need this in the present work, it is worth pointing out that the isomorphism types of finite groups acting effectively and symplectically on $M_{\lambda}=(\SSS,\oml)$ have been recently determined by C. Sáez-Calvo in his thesis. In particular, Theorem 5.17 of~\cite{CSC} implies that the finite abelian groups acting effectively and symplectically on $M_{\lambda}$ must be isomorphic to a cyclic group $\Z_{p}$, a product $\Z_{p}\times \Z_{q}$, or, in the case of $(\SSS,\oml)$ with $\lambda=1$, to a group $G$ lying in an exact sequence $0 \to \Z_{p} \times \Z_{p} \to G \to \Z/2\Z \to 0$ in which $\Z/2\Z$ acts on $\SSS$ by permuting the factors.
\end{remark}

%%%%%%%%%%%%%%%%%%%%%%%%%%%%%%%%%%%%%%%%%%%%%%%%%%%%%%%%%%%%%%%%%%%%%%%%%%%%%%%%
\subsection{Hamiltonian actions and almost-complex structures} 
%%%%%%%%%%%%%%%%%%%%%%%%%%%%%%%%%%%%%%%%%%%%%%%%%%%%%%%%%%%%%%%%%%%%%%%%%%%%%%%%
Let $r=2k$. The standard toric action $\T_{r}$ is characterized by the existence of an invariant embedded symplectic sphere $C_{r}$ of self-intersection $-r$ and of symplectic area $\lambda-k$. This sphere is holomorphic with respect to the invariant Hirzebruch complex structure $J_{r}$ and represents the homology class $D_{r}:=B-kF$, where $B:=[S^2\times\{*\}]$ is the base class, and where $F:=[\{*\}\times S^2]$ is the fiber class. In particular, if a Hamiltonian group action $G$ extends to a toric action $\T_{r}$, then there exists a $G$-invariant almost complex structure $J$ and a $G$-invariant $J$-holomorphic sphere in class $D_{r}$. As explained in~\cite{ChPin-Memoirs} Chapter~2, the converse holds for Hamiltonian circle actions. The following Lemma by Chiang and Kessler shows that the converse still holds for Hamiltonian actions by finite cyclic groups.  

\begin{prop}[Proposition 4.7 in~\cite{Liat}]\label{Prop:LiatRiver-Strata} Consider a $\Z_n(a,b;r)$ action on $(\SSS,\oml)$. Let $J$ be a $\Z_{n}$-invariant $\omega_{\lambda}$-compatible almost complex structure. Assume that the class  
$D_{2s}:=B-sF$ is represented by a simple $J$-holomorphic sphere. Then the action is equivariantly symplectomorphic to a $\Z_n(a^\prime,b^\prime; 2s)$ action that acts as a subgroup of the torus action $\T_{2s}$.\qed
\end{prop}

Let $\jjoml$ be the space of compatible almost complex structures on $M_{\lambda}=(\SSS,\oml)$ and let $\jzoml$ be the subspace of invariant structures. There is a retraction $\jjoml\to\jzoml$ defined by averaging the metrics associated to compatible pairs $(\om,J)$, showing that the subspace $\jzoml$ is contractible. It is shown in~\cite{MR1775741} that the space $\jjoml$ is stratified into $(\ell+1)$ strata
\begin{equation}\label{eq:Stratification}
\jjoml=U_0\sqcup U_2\sqcup\cdots\sqcup U_{2\ell}
\end{equation}
characterized by the existence of embedded $J$-holomorphic spheres representing the class $D_{r}$, that is,
\[U_{r}=\{J~|~\text{there exists a $J$-holomorphic sphere sphere $C$ with~}[C]=D_{r}\}\]
In particular, the strata are in bijection with the equivalence classes of toric actions on $M_{\lambda}$. It can be shown that $U_0$ is open and dense in $\jjoml$, and that each stratum $U_{r}$, where $r=2k$ and $k\geq 1$, is a  submanifold of codimension $4k-2$. The group $\Symp_h(M_{\lambda})$ of symplectomorphisms acting trivially in homology acts on $\jjoml$ preserving the stratification. Furthermore, each stratum $U_{r}$ is homotopically equivalent to the orbit of $J_{r}$, with stabilizer equal to the identity component of the Kahler isometry group of the pair $(\oml,J_{r})$. Consequently, 
\[U_{0}\simeq\Symp_h(M_{\lambda})/\SO(3)\times\SO(3) \text{~and~} U_{2k}\simeq\Symp_h(M_{\lambda})/\SO(3)\times S^{1} \text{~for~} 1\leq k\leq\ell.\] 
As explained in~\cite{MR1775741} and~\cite{AGK}, applying the Borel construction to the $\Symp_h(M_{\lambda})$ -invariant stratification~\eqref{eq:Stratification} shows that the classifying space $B\Symp_h(M_{\lambda})$ is homotopically equivalent to an iterated pushout of finitely many classifying spaces of Kahler isometry subgroups.

In the presence of a Hamiltonian $\Z_{n}$ action on $M_{\lambda}=(\SSS,\oml)$, the stratification of $\jjoml$ induces a partition of the space of invariant, compatible, almost-complex structures
\begin{align}
\jzoml&=\left(\jzoml\cap U_0 \right)\sqcup \left(\jzoml\cap U_2 \right)\sqcup\cdots\sqcup \left(\jzoml\cap U_{2\ell} \right)\notag\\
&=U_0^{\Z_n}\sqcup U_2^{\Z_n}\sqcup\cdots\sqcup U_{2\ell}^{\Z_n}\label{eq:stratification}
\end{align}
In particular, by Proposition~\ref{Prop:LiatRiver-Strata}, the invariant strata are in bijection with the non-equivalent 2-tori the $\Z_n$ action extends to.  
\begin{cor}\label{cor:strata-cyclic-case}
Let $\Z_n(a,b;r)$ be a symplectic action on $M_{\lambda}=(\SSS,\oml)$ with $\lambda >1$. Then the space of $\Z_n$-invariant, compatible, almost complex structures $\jzoml$ intersects the strata $U_{r'}$ if, and only if, $\Z_n(a,b; r)$ is equivariantly symplectomorphic to a $\Z_n(a^\prime,b^\prime;r')$ action that acts as a subgroup of the torus action $\T_{r'}$.
\end{cor}
\begin{proof}
Since $\lambda>1$, every symplectomorphism of $(\SSS,\oml)$ acts trivially on homology and hence acts on $\joml$ preserving the stratification. If $\Z_n(a,b;r)$ is symplectomorphic to a $\Z_n(a^\prime,b^\prime; r^\prime)$ action via a symplectomorphism $\phi$, then the pullback $\phi^*J_{r^\prime} \in U_{r^\prime}$ of the standard almost complex structure $J_{r^\prime} \in U_{r^\prime}$ is invariant under the $\Z_n(a,b;r)$ action. Thus  $\jzoml \cap U_{r^\prime}\neq\emptyset$. Conversely, if the space of $\Z_n$-equivariant complex structures $\jzoml$ intersects the strata $U_{r^\prime}$, then there is an invariant curve in class $D_{r'}$, and the result follows from Proposition~\ref{Prop:LiatRiver-Strata}.
\end{proof}

Let $\J_{\oml,l}$ be the space of compatible almost complex structures of regularity $C^l$ on $M$, endowed with $C^l$ topology, and let $\J_{\oml,l}^{\Z_n}$ be the subspace of invariant structures. Being a space of sections, $\J_l$ is a smooth Banach manifold (an explicit atlas can be constructed using the Cayley transform, see for instance~\cite{Smolentsev}). As in the non-equivariant case, we use $J$-holomorphic curves techniques to show that each invariant stratum is a submanifold of $\J_{\oml,l}^{\Z_n}$.
\begin{prop}\label{prop:InvariantStrataAreSubmanifolds} Each non-empty invariant stratum $U^{\Z_{n}}_{r,l}$ is a co-oriented Banach submanifold of $\J_{\oml,l}^{\Z_n}$ of finite codimension.
\end{prop}
The group $\Symp_h^{\Z_n}(M_{\lambda})$ of equivariant symplectomorphisms acting trivially on homology acts on $\jzoml$ preserving the stratification above. Under some conditions on the symplectic form $\oml$ and on the action $\Z_{n}(a,b;r)$ that we will specify later, we can show that each stratum is homotopically equivalent to a $\Symp_h^{\Z_n}(M_{\lambda})$ orbit. For the moment, we simply make this  into a definition.
\begin{defn}\label{defn:HomogeneityProperty} We say that a Hamiltonian $\Z_{n}$ action satisfies the \emph{homogeneity property} if each non-empty invariant stratum in the decomposition
\[\jzoml= U_0^{\Z_n}\sqcup U_2^{\Z_n}\sqcup\cdots\sqcup U_{2\ell}^{\Z_n}\]
is homotopically equivalent to the $\Symp_h^{\Z_n}(M_{\lambda})$ orbit of a standard structure $J_{r}$ with stabilizer homotopically equivalent to the subgroup $\Iso^{\Z_{n}}_h(\oml,J_{r})$ of equivariant Kahler isometries. 
\end{defn}
Following~\cite{MR1775741} and~\cite{AGK}, the homotopy decomposition of the space $\jzoml$ into $\Symp_h^{\Z_n}(M_{\lambda})$ orbits should yield an homotopy decomposition of the classifying space $B\Symp_h^{\Z_n}(M_{\lambda})$ in terms of the classifying spaces of the possible non-equivalent toric extensions of the $\Z_{n}$ action. However, due to the lack of a complete classification of Hamiltonian $\Z_n$ actions on $\SSS$ up to equivariant symplectomorphisms, we do not know necessary and sufficient conditions for two actions $\Z_{n}(a,b;r)$ and $\Z_{n}(a',b';r')$ to be equivalent. It follows that we cannot determine all the possible toric extensions of an arbitrary Hamiltonian $\Z_n$ action. Nevertheless, by again imposing some conditions on the symplectic form $\oml$ and on the action $\Z_{n}(a,b;r)$, it is possible to ensure that the action admits either one or two toric extensions. Assuming this for the moment, we now explain how the homotopy type of $B\Symp_h^{\Z_n}(M_{\lambda})$ is determined from the orbit decomposition of $\jzoml$ in these two special cases. 

%%%%%%%%%%%%%%%%%%%%%%%%%%%%%%%%%%%%%%%%%%%%%%%%%%%%%%%%%%%%%%%%%%%%%%%%%%%%%%%%
\subsubsection{Case of a unique toric extension}
%%%%%%%%%%%%%%%%%%%%%%%%%%%%%%%%%%%%%%%%%%%%%%%%%%%%%%%%%%%%%%%%%%%%%%%%%%%%%%%%
If we assume the $\Z_n$ action only extends to a single toric action $\T_{r}$ with $r\neq0$, and if the action has the homogeneity property, then it immediately follows that there are homotopy equivalences
\[\{*\}\simeq\jzoml = U_{r}^{\Z_n}\simeq \Symp_h^{\Z_n}(M,\oml)/\Iso^{\Z_{n}}_h(\oml, J_{r})\]
\begin{thm}\label{thm:ConditionalHomotopyTypeOnlyOneStratum} Suppose a Hamiltonian $\Z_n$ action only extends to a single toric action $\T_{r}$ with  $r\neq 0$. Suppose also that the action has the homogeneity property. Then
\[\Symp_h^{\Z_n}(M,\oml)\simeq\Iso^{\Z_{n}}_h(\oml, J_{r}).\]
\end{thm}

The case when $r=0$ is more intricate due to the existence of exceptional symmetries. The homotopy type of the entire group of equivariant symplectomorphisms is listed in Table~\ref{table:Actions on SSS r=0}.
%%%%%%%%%%%%%%%%%%%%%%%%%%%%%%%%%%%%%%%%%%%%%%%%%%%%%%%%%%%%%%%%%%%%%%%%%%%%%%%%
\subsubsection{Case of exactly two toric extensions along a common circle}
%%%%%%%%%%%%%%%%%%%%%%%%%%%%%%%%%%%%%%%%%%%%%%%%%%%%%%%%%%%%%%%%%%%%%%%%%%%%%%%%
Let us assume that the $\Z_n(a,b;r)$ action extends to exactly two non-equivalent toric actions $\T_{r}$ and $\T_{r'}$. 
\begin{prop}\label{prop:CodimensionTwo}
Suppose the $\Z_n(a,b;r)$ action extends to exactly two non-equivalent tori $T_r$ and $T_{r'}$ that intersect along the circle $K=S^1(a,b;r)$. 
Then
\begin{enumerate}
\item $a=\pm 1$. 
\item There is a decomposition $\jzoml=U_{r}^{\Z_n}\sqcup U_{r'}^{\Z_n}$ into an open stratum $U_{r}^{\Z_n}$ and  stratum $U_{r'}^{\Z_n}$ of positive codimension.
\item Moreover, the isotropy representation of $\T_{r'}$ on the normal fiber $V$ of $U_{r'}^{\Z_n}$ at $J_{r'}$ has global stabilizer $K$.
\end{enumerate}
\end{prop}
\begin{proof}
This is proven in~\cite{ChPin-Memoirs}. The first statement is Corollary~2.22 and Corollary~2.23. The second and third statements follow from Corollary 5.13, Theorem~5.17, and Corollary~5.18.
\end{proof}

The decomposition of $\jzoml$ into two submanifolds yields a pushout diagram of topological manifolds
\[
\begin{tikzcd}
\Nn U_{r'}^{\Z_n}\setminus\, U_{r'}^{\Z_n} \arrow{r}{} \arrow{d}{} & \Nn U_{r'}^{\Z_n} \arrow{d}{} \\
U_{r}^{\Z_n} \arrow{r}{} & \jzoml=U_{r}^{\Z_n}\cup \Nn U_{r'}^{\Z_n}
\end{tikzcd}
\]
where $\Nn U_{r'}^{\Z_n}$ is a tubular neighborhood of $U_{r'}^{\Z_n}$ in $\jzoml$. To avoid technicalities, let us further assume that we can choose $\Nn U_{r'}^{\Z_n}$ to be $\Symp_h^{\Z_n}(M_\lambda)$-invariant and equivariantly homeomorphic to a tube $\Symp_h^{\Z_n}\times_{\T_{r'}} V$. Then, the previous diagram is an equivariant diagram of $\Symp_h^{\Z_n}(M_\lambda)$-spaces, and applying the homotopy orbit functor (i.e. the Borel construction), we  obtain a homotopy decomposition of $B\Symp_h^{\Z_n}(M_\lambda)$ as a pushout
\[
\begin{tikzcd}
\left(S^{1}\right)_{h\T_{r'}}\simeq BK \arrow{r}{} \arrow{d}{} & B\T_{r'} \arrow{d}{} \\
B\T_{r} \arrow{r}{} & X\simeq B\Symp_h^{\Z_n}(M_\lambda)
\end{tikzcd}
\]
where $\left(S^{1}\right)_{h\T_{r'}}$ is the homotopy orbit of the isotropy representation of $\T_{r'}$ on the unit circle $S^{1}\subset V$. After looping, we finally obtain a pushout of H-spaces
\[
\begin{tikzcd}
K=S^1 \arrow{r}{} \arrow{d}{} & \T_{m'} \arrow{d}{} \\
\T_{r} \arrow{r}{} & P\simeq\Omega X\simeq\Symp_h^{\Z_n}(M_\lambda)
\end{tikzcd}
\]
that describes the homotopy type of the topological group $\Symp_h^{\Z_n}(M_\lambda)$ as a pushout of Lie subgroups $P:=\pushout(\T_r\from S^1\to\T_{r'})$. Although this discussion can be made rigorous using the framework developed in~\cite{AGK}, we prefer to establish this homotopy equivalence using a slightly different approach. Following Anjos and Granja~\cite{AG}, we show in Section~\ref{section:Centralizers-Zn(1,b;r)} that there is an isomorphism of Pontryagin rings between $H_*(\Symp_h^{\Z_n}(M_{\lambda});\Z)$ and $H_*(P;\Z)$. Then, standard arguments imply that there is a homotopy equivalence of H-spaces between $\Symp_h^{\Z_n}(M_{\lambda})$ and $P$. In particular, we obtain the following description of the homotopy type of $\Symp_h^{\Z_n}(M_{\lambda})$ as a topological space.
\begin{thm}\label{thm:ConditionalHomotopyTypeTwoStrata}
Suppose a Hamiltonian $\Z_n$ action extends to exactly two toric actions $\T_{r}$ and $\T_{r'}$, with $r'>r$,  intersecting along a common circle. Suppose also that the action has the homogeneity property. Then, as a topological space,
\[\Symp_h^{\Z_n}(M,\oml)\simeq\Omega S^{3} \times S^{1}\times S^{1}\times S^{1}.\]
\end{thm}
%

%%%%%%%%%%%%%%%%%%%%%%%%%%%%%%%%%%%%%%%%%%%%%%%%%%%%%%%%%%%%%%%%%%%%%%%%%%%%%%%%
\subsection{Tractable Hamiltonian $\Z_{n}$ actions}
%%%%%%%%%%%%%%%%%%%%%%%%%%%%%%%%%%%%%%%%%%%%%%%%%%%%%%%%%%%%%%%%%%%%%%%%%%%%%%%%
To complete our analysis of $\Symp_{h}^{\Z_{n}}(M_{\lambda})$, we need to find conditions on the symplectic form $\oml$ and on the action $\Z_{n}(a,b;r)$, $r=2k$, ensuring that the following statements hold:
\begin{enumerate}
\item The $\Z_{n}$ action admits either exactly one or exactly two toric extensions.
\item The $\Z_{n}$ action has the homogeneity property.
\item In the case of exactly two toric extensions, the two tori intersect along a common circle $K$.
\end{enumerate}
A Hamiltonian $\Z_{n}$ action satisfying these properties will be called \emph{tractable}. From the discussion above, the homotopy type of the symplectic centraliser $\Symp_{h}^{\Z_{n}}(M_{\lambda})$ of a tractable $\Z_{n}$ action is given by Theorem~\ref{thm:ConditionalHomotopyTypeOnlyOneStratum} or by Theorem~\ref{thm:ConditionalHomotopyTypeTwoStrata} depending on whether $\Z_{n}$ admits one or two toric extensions.
 
We will now describe two distinct sets of conditions on the parameters $\{\lambda, n, a, b, r\}$, coming from two different kinds of considerations, ensuring tractability of $\Z_{n}$ actions. 

\subsubsection{$\Z_n(a,b;r)$ actions with $\gcd(a,n)\neq 1\pmod{n}$} The simplest condition ensuring tractability is to assume $\gcd(a,n)\neq 1\pmod{n}$.

\begin{prop}\label{prop:gcdConditionSingleToricExtension}
Consider a finite cyclic group action $\Z_n(a,b;r)$  with $\gcd(a,n)\neq 1 \pmod{n}$. Then the invariant curve $C_{r}$ of non-positive self-intersection $-r$ is $J$-holomorphic for all $J\in\jzoml$. Consequently, $\jzoml=U_{r}^{\Z_{n}}$ so that the $\Z_{n}$ action admits a single toric extension $\T_{r}$.
\end{prop}
The proof is given in Section~\ref{ProofofPropgcdcondition}.

\begin{prop}\label{prop:gcdConditionTractability}
Any finite cyclic group action $\Z_n(a,b;r)$ with $\gcd(a,n)\neq 1 \pmod{n}$, has the homogeneity property. Consequently, it is tractable. 
\end{prop}
\begin{proof}
For $r\neq0$, this is Corollary~\ref{Cor:SinglecurveInclassD2s}. For $r=0$, this is  Theorem~\ref{Thm:HomtopytypeZn(a,b;0)}.
\end{proof}

\subsubsection{$\Z_n(a,b;r)$ actions with $\gcd(a,n)= 1\pmod{n}$} When $\gcd(a,n)=1 \pmod{n}$, we first note that because $a$ is invertible in $\Z_n$, we can reparametrize the action $\Z_n(a,b;r)$ to get an action of the form $\Z_n(1,a^{-1}b;r)$.

\begin{exmp} For $n$ prime, any action $\Z_n(a,b;r)$ with $a\neq 0\pmod{n}$ is equivalent to $\Z_n(1,a^{-1}b;r)$.\eoe
 
\end{exmp}

Furthermore as shown in Example~\ref{example:NormalizerChenWilczynski}, the actions $\Z_n(-1,b;r)$ are equivariantly symplectomorphic to the action $\Z_n(1,b+r;r)$. Hence, when $\gcd(a,n)=1$, we can assume $a=1$.\\ 

For actions of the form $\Z_n(1,b;r)$, we use the Chen--Wilczy\'nski classification~\cite{Chen} of smooth $\Z_n$ actions on $\SSS$ to find suitable conditions on the parameters $n$, $b$, and $r$ to ensure tractability. As explained in Section~\ref{Subsection:Consequences of the classification}, this classification gives necessary and sufficient conditions for two actions $\Z_{n}(a,b;r)$ and $\Z_{n}(a',b';r')$ to belong to the same conjugacy class in the group $\Diff^{+}(\SSS)$ of orientation-preserving diffeomorphisms of $\SSS$ and, therefore, restricts the possible toric extensions of symplectic $\Z_{n}(1,b;r)$ actions. The following proposition follows from Lemma~\ref{lemma:Toric-Extensions-(1,b;r)} and Corollary~\ref{cor:NumberOfStrataE1E2}. 

\begin{prop}\label{prop:Chen-WilczynskiInequalities}
The following Hamiltonian $\Z_n(1,b;r)$ actions on $(\SSS,\oml)$ admit at most two toric extensions.
\begin{itemize} 
\item Actions with no fixed surfaces:

\begin{enumerate}
\item If $a=1$, $b \neq  \{0,r/2,r\}$, $n\geq 2\lambda$ and $2\lambda < \min\{|r-2b|,r-2b+2n \}$, then the $\Z_{n}$ action only extends to $\T_{r}$.
\item If $a=1$, $b\neq\{0,r/2,r\}$,  $n>2\lambda$ and  $2\lambda > |2b-r|$ then the action extends to the two tori $\T_r$ and $\T_{|r-2b|}$. 
\item If $a=1$, $b\neq\{0,r/2,r\}$,  $n>2\lambda$ and  $2\lambda > r-2b+2n$ then the action extends to the two tori $\T_r$ and $\T_{r-2b+2n}$. 
\end{enumerate}

\item Actions that contain surfaces in their fixed points set: Suppose the parameters $r\neq0$ and $\lambda\geq 1$ satisfy the inequality $n>2\lambda > r$. If either 

\begin{enumerate}[resume]
\item $a=1$ and $b =0$, or 
\item $a=1$ and $b=r$,
\end{enumerate}
then the $\Z_{n}$ action only extends to $\T_{r}$.
\end{itemize}
\end{prop}

To simplify the exposition, we say that actions of the form $\Z_n(1,b;r)$ that extend to a single toric action are of type E1, while actions that extend to two tori are of type E2. Not all the above actions satisfy the homogeneity property though. We now list three types of actions of the form $\Z_n(1,b;r)$ which are tractable.

\begin{enumerate}[label=(H{\arabic*}), start=1]

\item Actions $\Z_n(1,b;r)$ of type E1 satisfying $2b-r\neq n$ and $b\not\in\{0,r\}$
\item Actions $\Z_n(1,b;r)$ of type E2 satisfying $2b-r\neq n$ with the exclusion of $\Z_{2r}(1,b;r)$, $\Z_{2r+2}(1,r+1;r)$, and $\Z_{2r+2}(1,2r+1;r)$. 
\item Actions of type E1 with $b\in\{0,r\}$, that is, $\Z_n(1,0;r)$ and $\Z_n(1,r;r)$ with $n\geq 2\lambda>r\geq2$.

\end{enumerate}

Using the fact that $\Z_n(-1,b;r) \cong \Z_n(1,b+r,r)$ (See Example~\ref{example:NormalizerChenWilczynski}) we can similarly produce a list of all tractable actions of the form $\Z_n(-1,b;r)$. The full homotopy type of actions on $\SSS$ for which we can prove tractability is listed in Tables~\ref{table:Actions on SSS r not 0} and~\ref{table:Actions on SSS r=0} in Section~\ref{section:Summary}. Similar results holds for the non-trivial bundle $\CCC$ as well and can be found in Table~\ref{table:Actions on CCC}.
\\

\subsection{Remaining cases} 
Because there is no complete classification of $\Z_n(a,b;r)$ actions up to equivariant symplectomorphisms on $(\SSS,\oml)$ and $(\CCC,\oml)$, we are presently unable to determine the homotopy type of the centralizer $\Symp^{\Z_n}(\SSS,\oml)$ in some specific cases. The complete lists of the remaining actions is given in sections~\ref{section:RemainingSSS} and~\ref{section:RemainingCCC}.

%%%%%%%%%%%%%%%%%%%%%%%%%%%%%%%%%%%%%%%%%%%%%%%%%%%%%%%%%%%%%%%%%%%%%%%%%%%%%%%%
\section{Toric extensions of finite cyclic groups actions}\label{section:Toric extensions}
%%%%%%%%%%%%%%%%%%%%%%%%%%%%%%%%%%%%%%%%%%%%%%%%%%%%%%%%%%%%%%%%%%%%%%%%%%%%%%%%

%%%%%%%%%%%%%%%%%%%%%%%%%%%%%%%%%%%%%%%%%%%%%%%%%%%%%%%%%%%%%%%%%%%%%%%%%%%%%%%%
\subsection{Hamiltonian actions on $\SSS$ and $\CCC$}\label{Section:finite abelian groups}
%%%%%%%%%%%%%%%%%%%%%%%%%%%%%%%%%%%%%%%%%%%%%%%%%%%%%%%%%%%%%%%%%%%%%%%%%%%%%%%%

For any given integer $r\in\Z$, recall that the $r^{\text{th}}$ Hirzebruch surface $W_r$ is the complex submanifold of $\mathbb{C}P^1 \times \mathbb{C}P^2$ defined by the equation
\[
W_r:=  \left\{ \left(\left[x_1,x_2\right],\left[y_1,y_2, y_3\right]\right) \in \mathbb{C}P^1 \times \mathbb{C}P^2 ~|~  x^r_1y_2 - x_2^ry_1 = 0 \right\}
\]
The projection map $\mathbb{C}P^1 \times \mathbb{C}P^2 \longrightarrow \mathbb{C}P^1$ gives $W_r$ the structure of a $\mathbb{C}P^1$ bundle over $\mathbb{C}P^1$. It is diffeomorphic to $S^2 \times S^2$ when $r$ is even, while it is diffeomorphic to the non-trivial $S^2$ bundle over $S^2$, that is, to $\CP^2 \# \overline{\CP^2}$, when $r$ is odd. The complex surfaces $W_{r}$ and $W_{r'}$ are biholomorphically equivalent if, and only if, $r'=\pm r$. Consequently, for each non-negative integer $r\geq 0$ we have an integrable complex structure $J_r$ induced on $\SSS$ or $\CCC$ that is well defined up to oriented diffeomorphisms. Moreover, any complex structure on $\SSS$ or $\CCC$ is diffeomorphic to one of these Hirzebruch structures, see~\cite{Qin}.\\

When $r=2k\geq 0$ is even, we can choose any real number $\lambda > k$ and endow $\mathbb{C}P^1 \times \mathbb{C}P^2$ with the symplectic form $(\lambda -k) \sigma_1 \times \sigma_2$, where $\sigma_1$ is the standard symplectic form on $\mathbb{C}P^1$ of total area 1, and where $\sigma_2$ is the standard symplectic form on $\mathbb{C}P^2$ normalized in such a way that the symplectic area of a line is $1$. 

Restricting this symplectic form to $W_r$ makes it a K\"ahler manifold symplectomorphic to  $(\SSS,\oml)$. Similarly, when $r=2k+1\geq 1$ is odd, we can choose any $\lambda > k+1$ and restrict the symplectic form $(\lambda - (k+1)) \sigma_1 \oplus \sigma_2$ to $W_r$ to obtain a K\"ahler form $\oml$ on $\CCC$. The Lalonde-McDuff classification theorem~\cite{MR1426534} implies that, after rescalling, any symplectic form on either $\SSS$ or $\CCC$ is diffeomorphic to one of these normalized forms $\oml$.\\

\begin{figure}[H]
\centering       
\begin{tikzpicture}
\node[left] at (0,2) {$Q=(0,1)$};
\node[left] at (0,0) {$P=(0,0)$};
\node[right] at (4,2) {$R= (\lambda - k ,1)$};
\node[right] at (6,0) {$S=(\lambda + k ,0)$};
\node[above] at (2,2) {$D_r=B-kF$};
\node[right] at (5.15,1) {$F$};
\node[left] at (0,1) {$F$};
\node[below] at (3,0) {$B+ kF$};
\draw (0,2) -- (4,2) ;
\draw (0,0) -- (0,2) ;
\draw (0,0) -- (6,0) ;
\draw (4,2) -- (6,0) ;
\end{tikzpicture}
\caption{Even Hirzebruch polygon}
\label{hirz}
\end{figure}
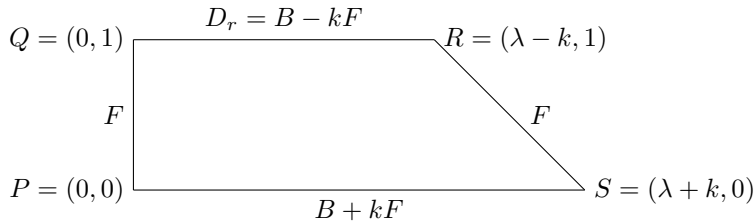

\begin{figure}[H]
\centering   
\label{hirZ1}
\begin{tikzpicture}
\node[left] at (0,2) {$Q=(0,1)$};
\node[left] at (0,0) {$P=(0,0)$};
\node[right] at (4,2) {$R= (\lambda - (k+1) ,1)$};
\node[right] at (6,0) {$S=(\lambda + k ,0)$};
\node[above] at (2,2) {$D_r=B-(k+1)F$};
\node[right] at (5.15,1) {$F$};
\node[left] at (0,1) {$F$};
\node[below] at (3,0) {$B+kF$};
\draw (0,2) -- (4,2) ;
\draw (0,0) -- (0,2) ;
\draw (0,0) -- (6,0) ;
\draw (4,2) -- (6,0) ;
\end{tikzpicture}
 \caption{Odd Hirzebruch polygon}\label{fig:OddHirzebruch}
\end{figure}
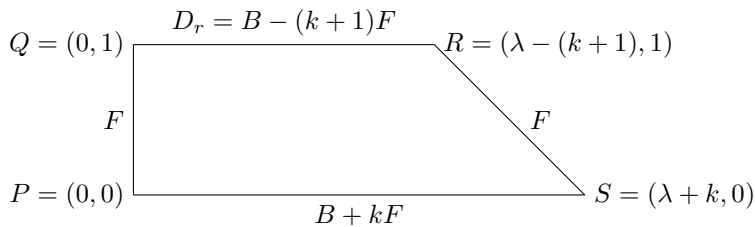

Given an integer $r\geq 0$, the torus $\T$ acts on $\mathbb{C}P^1 \times \mathbb{C}P^2$ via
\[
 \left(u,v\right) \cdot  \left(\left[x_1,x_2\right],\left[y_1,y_2, y_3\right]\right) = \left(\left[ux_1,x_2\right],\left[u^ry_1,y_2,vy_3\right]\right)
\]
with moment map 
\[\mu([x_1:x_2],[y_1:y_2:y_3])=\left(\frac{M|x_1|^2}{|x_1|^2+|x_2|^2}+\frac{r|y_1|^2}{|y_1|^2+|y_2|^2+|y_3|^2},~\frac{|y_3|^2}{|y_1|^2+|y_2|^2+|y_3|^2}\right)\]
where $M=\lambda-k$ when $r=2k$ and $M=\lambda-(k+1)$ when $r=2k+1$. This action leaves $W_r$ invariant and preserves both the complex and the symplectic structures. Its restriction to $W_r$ defines a toric action that we denote $\T_r$. When $r=2k$ is even, the image of the moment map (defined up to translations) is the polytope of Figure~\ref{hirz}

\noindent where the labels along the edges refer to the homology classes of the  $\T_r$ invariant spheres in $\SSS$, and where the vertices $P$,$Q$,$R$,$S$ are the fixed points of the torus action. 
Similarly, when $r=2k+1$ is odd, the moment map image is given in Figure~\ref{fig:OddHirzebruch}
\noindent where $B$ now refers to the homology class of a line $L$ in $\CP^2 \# \overline{\CP^2}$, $E$ is the class of the exceptional divisor, and $F$ is the fiber class $B - E$.\\

Near any fixed point, a Hamiltonian action by a subgroup of $\T_{r}$ is characterized by its weights. The following proposition follows from \cite{Montaldi} Theorem~2.4, see also~\cite{Karshon} Corollary~A.7 or \cite{Lerman-Tolman} Lemma~A.1.
\begin{prop}\label{prop:DefinitionOfWeightsZ_n}
Fix a generator $\zeta$ of $\Z_n$. Let $p$ be a fixed point of a symplectic action of $\Z_n$ on a $4$-manifold $(M,\om)$. Then there exist complex coordinates $(z,w)$ on a neighborhood of $p$ in $M$, and unique integers $\alpha_p$ and $\beta_p$, well defined modulo $n$, called the isotropy weights at $p$, such that
\begin{enumerate}
\item the action is $\zeta\cdot (z, w) = (\zeta^{\alpha_p}z, \zeta^{\beta_p}w)$,
\item the symplectic form is $\om = \frac{i}{2}(dz\wedge d\bar z + dw \wedge d\bar w)$,
\end{enumerate}
\end{prop}

Since for each fixed point $p$ the unordered set $\{\alpha_p,\beta_p\}$ is a symplectic invariant of the action, the collection of weights imposes strong restrictions on equivariant symplectomorphisms.
\begin{cor}\label{cor:ActionPreservesWeightsZ_n}
Consider a Hamiltonian action by a torus or by a finite cyclic group $\Z_n$ of order $n$ on a symplectic $4$-manifold $(M,\om)$. Let $\phi$ be an equivariant symplectomorphism. Then $\phi$ acts on the fixed point set preserving their weights as unordered sets. In particular, equivalent actions have the same set of weights.\qed
\end{cor}

The weights of the $\Z_n(a,b;r)$ action at the four toric fixed points are given in the following table.
\begin{table}[H]
\begin{center}
\begin{tabular}{ |p{2cm}||p{4cm}|p{5cm}|  }
\hline
Vertex & Weights for $\T_r$ action & Weights for the $\Z_n(a,b;r)$ action \\
\hline
\hline
P & $\big\{\{1,0\},\{0,1\}\big\}$   & $\{a,b\} \pmod{n}$ \rule{0pt}{10pt}\\
\hline
Q & $\big\{\{1,0\},\{0,-1\}\big\}$  & $\{a,-b\} \pmod{n}$\rule{0pt}{10pt}\\
\hline
R & $\big\{\{-1,r\},\{0,-1\}\big\}$ & $\{-a,ar-b\} \pmod{n}$ \rule{0pt}{10pt}\\
\hline
S & $\big\{\{-1,-r\},\{0,1\}\big\}$ & $\{-a,-ar+b\} \pmod{n}$ \rule{0pt}{10pt}\\
\hline
\end{tabular}
\end{center}
\caption{Weights of $\T_{r}$ and $\Z_n(a,b;r)$ actions}
\label{table_weights_Z_n}
\end{table}

\begin{remark}
Observe that under the $\Z_n(a,b;r)$ action, the stabilizer of the invariant curve $\overline{D}_r$ has order $\gcd(a,n)$, while the stabilizer of an invariant fiber has order $\gcd(b,n)$. 
\end{remark}

\begin{remark}
In what follows, we will often impose conditions on cyclic actions $\Z_n(a,b;r)$ on $(S^2\times S^2,\oml)$ by imposing numerical inequalities involving the parameters $\lambda\in\R_{\geq 1}$, $n\in\N$, $r\in\N$ and $a,b\in\Z_n$. These inequalities must be understood as inequalities in $\R$, with $a$ and $b$ taken in $\{0,\ldots,n-1\}$.
\end{remark}

The following lemma will be useful in our analysis. 

\begin{lemma}\label{weightZ_n}
Suppose $\Z_n$ acts by rotations on $S^2$ with weight $k$, fixing the north and south poles. Suppose that the action lifts to a complex line bundle $E$ over $S^2$. Then $\Z_n$ acts linearly on the fibers over the north and south poles; let $\alpha$ and $\beta$ be the weights for these actions. Then
\[\alpha - \beta = \pm ek \pmod{n}\]
where $e$ is the self-intersection of the zero section. 
\end{lemma}
\begin{proof}
See Lemma~5.4 in~\cite{Karshon}. 
\end{proof}
\begin{remark}
    Unlike the circle action case, the weights only determine the self-intersection of the curve modulo the order of the cyclic group $n$. Hence for exceptional actions such as $\Z_{2r}(1,b;r)$, it is not possible to distinguish between a curve with self-intersection  $-2r$ and a curve with self-intersection  $2r$  using just the weight calculations. Thus we impose further conditions on $a,b,n$ as in Lemma~\ref{lemma:ConditionsForDrCurvesOnlyThroughQandR} and Proposition~\ref{prop:ExistenceOfGoodFixedPointH1H2} in order to use the self-intersection formula above to uniquely determine the self-intersection of the curves given the weights.
\end{remark}
We close this section by recalling the classification of toric extensions of Hamiltonian circle actions.

\begin{prop} \label{prop:ToricExtensionsOfCircleActions}
Consider a circle action $S^1(a,b;r)$ on $(\SSS,\oml)$  or $(\CCC,\oml)$. Let $\epsilon_r=r\mod 2$. 
\begin{enumerate}
\item The action $S^1(a,b;r)$ only extends to the toric action $\T_r$ unless $a=\pm1$ and $b\not\in\{0,ar\}$.
\item Any action of the form $S^1(-1,b;r)$ is equivariantly symplectomorphic to $S^1(1,-b;r)$ through an action of the symplectic normalizer of $\T_r$.
\item When $a=1$, $b\not\in\{0,b\}$, and $|2b-r|+\epsilon_r\geq 2\lambda$, the action only extends to $\T_r$.
\item When $a=1$, $b\not\in\{0,b\}$, and $2\lambda>|2b-r|+\epsilon_r$, the action extends to both $\T_r$ and $\T_{|2b-r|}$.
\end{enumerate}
Moreover, in the case the $S^1(1,b;r)$ action extends to $\T_r$ and $\T_{|2b-ar|}$, it is equivariantly symplectomorphic to the following subcircle in $\T_{|2b-m|}$\,:
\begin{enumerate}[label=(4\alph*)]
    \item $S^1(1,b; |2b-r|)$ if $b>0$ and $b>r$, 
    \item $S^1(1,-b; |2b-r|)$ if $b>0$, $r>b$, and $2b-r < 0$,
    \item $S^1(1,b;|2b-r|)$ if $b>0$, $r>b$, and $2b-r > 0$, 
    \item $S^1(1,-b;|2b-r|)$ if $b<0$.\qed
\end{enumerate}
\end{prop}
\begin{proof}
All statements are taken from~\cite{ChPin-Memoirs}. The first one is Proposition~2.21, the second one is Proposition 2.18, and the other statements are taken from Corollaries~2.22 and~2.23.
\end{proof}

\subsection{Proof of Proposition~\ref{prop:gcdConditionSingleToricExtension}}\label{ProofofPropgcdcondition}
Recall that Proposition~\ref{prop:gcdConditionSingleToricExtension} states that an action $\Z_{n}(a,b;r)$ admits a single toric extension whenever $\gcd(a,n)\neq 1\pmod{n}$. This follows from the following observation.

\begin{lemma}[Lemma 2.5 in \cite{Liat}]{\label{LR}}
Let $(M,J)$ be an almost complex manifold such that $J$ is invariant under the action of a compact Lie group $G$ on $M$. Let $S$ be an embedded sphere in $M$. Assume that $S$ is a connected component of the fixed point set of a non-trivial subgroup $H\subseteq G$. Then $S$ is a $J$-holomorphic sphere.
\end{lemma}
\begin{proof}
To show $S$ is $J$-holomorphic we need to show that for any vector  $v \in T_xS$ $Jv \in T_xS$. As $S$ is pointwise fixed by a non-trivial subgroup $H$, all tangent vectors $v \in T_xS$ are characterised by the property that $dh \cdot v = v$ for all $h \in H$. Thus in order to show that that  $Jv \in T_xS$, it suffices to prove that $dh \cdot Jv = Jv$ for all $h \in H$. But this immediately follows from the equivariant of $J$ as
\[dh \cdot Jv = J(dh \cdot v) = Jv\]
\end{proof}

\begin{proof}[Proof of Proposition~\ref{prop:gcdConditionSingleToricExtension}]
From Table~\ref{table_weights_Z_n} in Section~\ref{Section:finite abelian groups}, we see that the stabilizer $H$ of the invariant curve $\overline{D}_r$ under the $\Z_n(a,b;r)$ action has order $\gcd(a,n)$, and is therefore non-trivial whenever $\gcd(a,n)\neq 0\pmod{n}$. The statement follows from Lemma~\ref{LR} and from the fact that each stratum $U_k$ in the partition~\eqref{eq:stratification} is characterized by the existence of a $J$-holomorphic curve in class~$D_k$.
\end{proof}

\begin{remark}
    We use the convention that $\gcd(0,n)=n$. In particular, Proposition~\ref{prop:gcdConditionSingleToricExtension} is also true for the action $\Z_n(0,1;r)$. 
\end{remark}

%%%%%%%%%%%%%%%%%%%%%%%%%%%%%%%%%%%%%%%%%%%%%%%%%%%%%%%%%%%%%%%%%%%%%%%%%%%%%%%%
\subsection{The Chen--Wilczy\'nski classification}
%%%%%%%%%%%%%%%%%%%%%%%%%%%%%%%%%%%%%%%%%%%%%%%%%%%%%%%%%%%%%%%%%%%%%%%%%%%%%%%%

As described in the works of W. Chen \cite{Chen} and D. Wilczy\'nski \cite{W}, one can establish the existence of $\Z_n$-equivariant orientation-preserving diffeomorphisms between different actions $\Z_n(a,b;r)$ and $\Z_n(a',b';r')$ defined on Hirzerbruch surfaces. In Chen's terminology, these equivariant diffeomorphisms are of six different types which are characterised by the following six conditions on the triples $(a,b;r)$ and $(a',b';r')$.
\begin{itemize}
    \item Type $c_1$: When $a^\prime = -a$, $b^\prime = -b$ and $r^\prime = r$. 
    \item Type $c_2$: When $a^\prime = -a$, $b^\prime = b - ra$ and $r^\prime = r$. 
    \item Type $c_3$: When $a^\prime = a$, $b^\prime = -b$ and $r^\prime = -r$. 
    \item Type $c_4$: When $a^\prime = -b$, and $b^\prime = -a$, and $r^\prime = r = 0$. 
    \item Type $c_5$: When $a^\prime = a$, $b^\prime = b$, and $r^\prime \equiv r $ (mod 2n), assuming $\gcd(a,n)=\gcd(a',n)=1$.
    \item Type $c_6$: When  $a^\prime = a$, $b^\prime = b$, and $r^\prime a^\prime \equiv 2b - ra$ (mod 2n), assuming $\gcd(a,n)=\gcd(a',n)=1$.
\end{itemize}
We call the above $\Z_{n}$ equivariant diffeomorphism \emph{\textit{standard}} of type $c_1, \cdots, c_6$. One of the main results of W. Chen~\cite{Chen} is the following classification theorem that builds on the previous work of D.~Wilczy\'nski~\cite{W}. 

\begin{thm}[Chen-Wilczy\'nski classification, see~\cite{Chen} Theorem 1.5]\label{Chen} Two cyclic actions $\Z_n(a,b;r)$ and $\Z_n(a',b';r')$ are orientation-preserving equivariantly diffeomorphic if, and only if, there is a composition of standard equivariant diffeomorphisms between them.
\end{thm}

\begin{remark}
The notation used in Chen~\cite{Chen} and Wilczy\'nski~\cite{W} is slightly different from the one used in the present document. In particular, the action we denote by $\Z_n(a,b;r)$ corresponds to the action $F_r(a,-b)$ in \cite{Chen} and \cite{W}. This explains the differences between the conditions $c_1 \cdots c_6$ listed above and the ones given in \cite{Chen} Section~5.
\end{remark}

\begin{remark}\label{rmk:Chen-WilczynskiNegativeIndex}
For convenience, the classification is stated for Hirzebruch surfaces $W_{r}$ with indices $r\in\Z$. The diffeomorphisms of type $c_{3}$ identify $W_{r}$ with $W_{-r}$ and induce a map diffeotopic to $z\mapsto z^{-1}$ between the fibers $F_{r}$ and $F_{-r}$.\
\end{remark}

\begin{remark}\label{rmk:Chen-WilczynskiConditionsMod2n}
When $\gcd(a,n)=1$, the existence of equivariant diffeomorphisms of types $c_5$ shows that the action $\Z_{n}(a,b;r)$ is equivariantly diffeomorphic to $\Z_{n}(a,b;r')$ whenever $r'\equiv r\pmod{2n}$. In particular, when $\gcd(a,n)=1$, we can ignore diffeomorphisms of type $c_{5}$ from the above list and, instead, consider the conditions on $r$ and $r'$ as given modulo $2n$ in the descriptions of the other five types.
\end{remark}

\begin{exmp}\label{example:NormalizerChenWilczynski}
Let $C(\T_{r})$ and $N(\T_{r})$ be the centralizer and normalizer of the torus $\T_{r}$ in $\Symp(W_{r},\oml)$. We know from~\cite{ChPin-Memoirs}~Proposition 2.16 that for $r\geq 1$, the symplectic Weyl group $W(\T_{r})=N(\T_{r})/C(\T_{r})$ is isomorphic to $\Z_2$, and that its action on subcircles of $\T_r$ sends $S^{1}(a,b;r)$ to $S^{1}(-a,b-ar;r)$. This shows that equivariant diffeomorphisms of type $c_2$ with $r=r'$ are always realizable by symplectomorphisms of~$W_r$. In particular, the Hamiltonian action $\Z_n(-1,b;r)=\Z_n(n-1,b;r)$ is symplectically equivalent to $\Z_n(1,b+r;r)$.\eoe %This will prove useful later.
\end{exmp}

\begin{exmp}\label{example:Type c2 is Kahler when n=2a}
In the last example, if we further assume $2a=n$ and $r$ is even, then the $\Z_{2a}(a,b;r)$ equivariant symplectomorphism of type $c_2$ is realized by the $\Z_n$ equivariant Kahler involution of $W_r$ given by
\[([x_1:x_2],[y_1:y_2:y_3])\mapsto ([x_2:x_1],[y_2:y_1:y_3])\]
that swaps the fixed points $Q$ and $R$, and also swaps $P$ and $S$.
\end{exmp}

\begin{exmp}\label{example:c5NotSymplectic}
Consider the Hamiltonian $\Z_n(1,b;2k)$ action on $M_\lambda=(\SSS, \oml)$, where $\lambda=\ell+\delta$, $0<\delta\leq 1$. From Proposition~\ref{prop:SummaryActionsS2xS2}, we know that $M_\lambda$ admits exactly $(\ell+1)$ inequivalent Hamiltonian toric actions given by $\T_0,\ldots,\T_\ell$. Therefore, any Hamiltonian toric extension $\T_r$ of $\Z_n(1,b;2k)$ must satisfy $0\leq r\leq\ell$. On the other hand, according to the Chen--Wilczy\'nski classification, actions of the form $\Z_n(1,b;2k\pm2nq)\subset\T_{2k\pm2nq}$ are all equivariantly diffeomorphic to $\Z_n(1,b;2k)$ via diffeomorphisms of type $c_5$. This shows that, in general, diffeomorphisms of type $c_5$ cannot be chosen to be symplectic.\eoe
\end{exmp}

\begin{exmp}\label{example_obstruction} Consider Hamiltonian actions of the form $\Z_3(2,1;6k)\subset\T_{6k}$ on $M_\lambda=(\SSS, \oml)$, where $\lambda=\ell+\delta$, $0<\delta\leq 1$ such that $0\leq6k\leq\ell$ and $2\lambda >> 3$. These actions have four isolated fixed points $P$, $Q$, $R$, and $S$ (as in fig~\ref{hirz}) whose weights are all the same for every action in the family $\{\Z_3(2,1;6k)\}_{k\geq0}$. Moreover,  the Chen--Wilczy\'nski classification implies that these actions are all equivariantly diffeomorphic to $\Z_3(2,1;0)$, via diffeomorphisms of type $c_5$. However, we are unable to tell whether they are equivariantly symplectomorphic. Consequently, we cannot determine the Hamiltonian toric extensions of $\Z_3(2,1;0)$. 

This is in sharp contrast with the case of Hamiltonian circle actions for which the set of weights together with the cohomology class of the symplectic form determine the Hamiltonian toric extensions. This later fact directly follows from Karshon's equivariant classification of Hamiltonian $S^1$ actions on $4$-manifolds, see~\cite{ChPin-Memoirs} Pro\-po\-si\-tion~3.22.\eoe
\end{exmp}

\begin{exmp} Contrary to $S^1$ actions, it is possible to produce $\Z_n$ actions that extend to more than two toric actions. Consider the action $\Z_n(1,1;2)$. Note that this $\Z_n$ action is contained in infinitely many circles $S^1(a,b;r) \subset \T_r$ satisfying $a \equiv 1$ mod($n$) and $b \equiv 1$ mod($n$). In particular $\Z_n(1,1;2)$ is contained in the family of circles $S^1(1,1+kn;2)$, $k\in\N$. By Proposition~\ref{prop:ToricExtensionsOfCircleActions}, we know that $S^1(1,1+kn;2)$ is symplectomorphic to $S^1(1,1+kn;2kn)$ as long as $2\lambda > |2(1+kn)-2|$. It follows that $\Z_n(1,1;2)$ is equivariantly symplectomorphic to $\Z_n(1,1;2kn)$ for each $0\leq k\leq M$ where $M$ is the largest integer such that $2\lambda > |2(1+nM)-2|$. Consequently, $\Z_n(1,1;2)$ extends to the toric actions $\T_0, ~ \ldots~,T_{2kn},~\ldots \T_{2Mn}$ among others.\\

When we further specialize to the case $n=2$, we can show that the action $\Z_2(1,1;r)$ extends to all tori $T_k$ such that $k < 2\lambda$ and $k \equiv r \mod 2$. Let $m$ denote the maximum integer such that $m < 2\lambda$. We have the following sequence of isomorphism: $\Z_2(1,-1;m) \cong \cdots \cong \Z_2(1,-1;r+4)\cong \Z_2(1,-1;r+2) \cong \Z_2(1,-1;r)\cong \Z_2(1,1;r) \cong \Z_2(1,-1;|r-2|) \cong \Z_2(1,1;|r-2|) \cong \Z_2(1,-1;|r-4|) \cdots \Z_2(1,1,0)$ where we use the isomorphism $S^1(1,-1,r) \cong S^1(1,-1;r+2)$ to show that $\Z_2(1,1;r)$ intersects tori $\T_k$ for all $k \geq r$ and we use the isomorphism $S^1(1,1;r) \cong S^1(1,-1;|r-2|)$ to show that $\Z_2(1,1;r)$ intersects tori $\T_k$ for all $k \leq r$ (See Corollary 2.22 in \cite{ChPin-Memoirs} for more details about the $S^1$-equivariant isomorphisms). The isomorphisms $\Z_2(1,-1;r) \cong \Z_2(1,1;r)$ hold because we consider $\Z_2$ actions. The homotopy type of $\Symp^{\Z_n}(\SSS,\oml)$ for $\Z_2(1,1;r)$ actions with $\lambda \geq 1$ will be explored in a future work. \eoe
\end{exmp}

%%%%%%%%%%%%%%%%%%%%%%%%%%%%%%%%%%%%%%%%%%%%%%%%%%%%%%%%%%%%%%%%%%%%%%%%%%%%%%%%
\subsection{Consequences of the smooth classification for toric extensions}\label{Subsection:Consequences of the classification}
%%%%%%%%%%%%%%%%%%%%%%%%%%%%%%%%%%%%%%%%%%%%%%%%%%%%%%%%%%%%%%%%%%%%%%%%%%%%%%%%
As Example~\ref{example_obstruction} illustrates, the Chen-Wilczy\'nski classification only gives obstructions to symplectic toric extensions. In order to determine exactly which tori an action $\Z_n(a,b;r)$ extends to, one has to impose further constraints on $a$, $b$, $r$ and $n$. By analogy with Hamiltonian circle actions, we consider actions of the form 
\begin{itemize}
\item $\Z_n(a,b;r)$ with $r\geq 2$ and $a=\pm1\pmod{n}$.
\end{itemize}
For such an action to be Hamiltonian, we must have $2\lambda>r$. By Example~\ref{example:NormalizerChenWilczynski}, any Hamiltonian action of the form $\Z_n(-1,b;r)$ is symplectomorphic to a Hamiltonian action $\Z_n(1,b+r;r)$. Consequently, we can restrict ourselves to $a=1$. 
We further impose the inequality 
\begin{itemize}
\item $n\geq2\lambda$. 
\end{itemize}

Note that the above inequalities imply that the order of the cyclic group acting on $\SSS$ is at least $3$. Our goal is to prove an analogue of Proposition~\ref{prop:ToricExtensionsOfCircleActions} for this type of $\Z_n$ actions.\\

By the Chen-Wilczynski classification, the set of $\Z_{n}(a',b';r')$ actions diffeomorphic to $\Z_{n}(a,b;r)$ can be identified with all triples $(a',b';r')\in\Z_{n}\times\Z_{n}\times\Z$ that can be obtained from $(a,b;r)$ by applying finitely many transformations of types $c_{1},\ldots,c_{6}$. For the two types of actions that interest us, $a\equiv\pm1\pmod{n}$, and $n > 2\lambda > r>0$, so that $\gcd(a,n)=1$ and $r\not\equiv 0\pmod{2n}$. Consequently, we can apply the transformation $c_5$ to the triple $(a,b;r)$, which amounts to taking $r$ modulo $2n$. On the other hand, the transformation $c_{4}$ is excluded as long as the triple $(a,b;r)$ is not equivalent to any triple of the form $(a',b';0)$ under the action of $c_1$, $c_2$, $c_3$, and $c_6$. We are lead to investigate the action of the group generated by transformations of types $c_1$, $c_2$, $c_3$, and $c_6$ onto the set of triples 
\[\mathcal{T}_{n}=\{(a,b;r)\in\Z_{n}\times\Z_{n}\times\Z_{2n}~|~r\not\equiv 0\mod2n, \text{~and~}\gcd(a,n)=1\}\]
whose orbits describes actions diffeomorphic to $\Z_n(a,b;r)$.

\begin{lemma}\label{lemma:OrbitChen-Wilczynski} Let $B$ be the group generated by transformations of types $c_1$, $c_2$, $c_3$, and $c_6$.
\begin{enumerate}
\item $B$ is isomorphic to
\[\Z_{2}\times D_{4}=\langle c_{1}, c_{2}, c_{6}~|~c_{1}^{2}, c_{2}^2, c_{6}^2, (c_{1}c_{2})^2, (c_{1}c_{6})^2, (c_{2}c_{6})^4\rangle\]
\item In particular, when $a=\pm1$, the orbit of $(a,b;r)$ in $\mathcal{T}_{n}$ consists of triples $(a',b';r')$ for which $r'\equiv \pm r\mod2n$ or $r'\equiv \pm (2b-ar)\pmod{2n}$.
\item Assuming $a=\pm1$, $r\not\equiv 0\pmod{2n}$, and $2b-ar\not\equiv 0\pmod{2n}$, the $\Z_n$ actions diffeomorphic to $\Z_n(a,b;r)$ are of the form $\Z_n(a',b';r'+2kn)$, where $k\in\Z$, and where $(a',b';r')$ is in the orbit of $(a,b;r)$ under the action of $B$.
\end{enumerate}
\end{lemma}
\begin{proof}
It is easy to check that $c_1$, $c_2$, $c_3$, and $c_6$ are involutions, that $c_1$ commutes with the other generators, and that $c_3=(c_2c_6)^2$. This proves the first statement. The second statement is a straightforward computation. The third statement follows from the second. Indeed, under the conditions $a=\pm1$, $r\neq 0\pmod{n}$, and $2b-ar\neq 0\pmod{2n}$, all triples $(a',b';r')$ in the orbit of $(a,b;r)$ have $r'\neq0\pmod{2n}$, so that the transformation of type $c_4$ is not applicable to any $(a',b';r')$.
\end{proof}

\begin{lemma}\label{lemma:Toric-Extensions-(1,b;r)}
Consider a Hamiltonian $\Z_n(1,b;r)$ action on $(\SSS,\oml)$ with $n\geq 2\lambda>r\geq2$, and $2b-r\neq 0$. %, and $b \neq \{0,r\}$. 
Then, either
\begin{enumerate}
\item $(2b-r)$, $(r-2b)$, and $(r-2b+2n)$ are outside the symplectic range $(0,2\lambda)$ and $\Z_{n}(1,b;r)$ only extends to the toric action $\T_r$, or
\item there is exactly one index $r'$ among $(2b-r)$, $(r-2b)$, or $(r-2b+2n)$ in the range $(0,2\lambda)$. In this case, $\Z_{n}(1,b;r)$ extends to $\T_r$ and $\T_{r'}$, and to no other toric actions.
\end{enumerate}
\end{lemma}
\begin{proof}
 
Under the conditions $n\geq 2\lambda>r>0$, the modular equalities $r=0\pmod{2n}$ and $2b-r=0\pmod{2n}$ are only possible if $r=0$ and $2b-r=0$, which are excluded by hypothesis. Consequently, Lemma~\ref{lemma:OrbitChen-Wilczynski}~(3) applies, showing that the only Hamiltonian $\Z_n(a',b';r')$ actions symplectomorphic to $\Z_n(1,b;r)$ must have $r'=\pm r \pmod{2n}$  or $r'=\pm(2b-r)\pmod{2n}$, with $0<r'<2\lambda$. It is easy to see that the only possible representatives of $r'$ in the symplectic range $(0,2\lambda)$ are $(2b-r)$, $(r-2b)$, and $(r-2b+2n)$. It is also easy to see that these three possibilities are mutually exclusive. Note also that if either $b=0$ or $b=r$, then $\pm(2b-r)=\pm r\pmod{2n}$, so that $r'=\pm r$ is the only possibility for $r'$ in these two cases.

If no representative of $r'$ in the range $(0,2\lambda)$ exists, then the only toric extension of $\Z_{n}(1,b;r)$ is $\T_{r}$. If exactly one representative $r'$ exists, then we claim that the corresponding action is Hamiltonian. Indeed, from Proposition~\ref{prop:ToricExtensionsOfCircleActions}, we know that the circle action $S^1(1,b;r)$ with $r\geq 2$ and $b\not\in\{0,r\}$ extends to a second tori $\T_{r'}$ if, and only if, $r'=|2b-ar|$, and $2\lambda>|2b-r|$. This shows that in case $r'=\pm(2b-r)$,  the action $\Z_n(1,b;r)$ does extends to $\T_{r'}$. In the case the only representative is $r'=r-2b+2n$, we have the following equivalences between Hamiltonian $\Z_n$ actions:
\[\Z_n(1,b;r)\simeq \Z_n(1,b-n;r)\simeq \Z_n(1,n-b;|2b-2n-r|)\simeq\Z_n(1,n-b;r-2b+2n)\]
where the first equivalence holds since $b$ is only defined modulo $n$, and where the second equivalence is given by Proposition~\ref{prop:ToricExtensionsOfCircleActions}~4d).
\end{proof}

\begin{defn}
Given a Hamiltonian action $\Z_{n}(1,b;r)$ with $n\geq 2\lambda>r\geq2$, and $2b-r\neq 0$, %and $b \neq \{0,r\}\pmod{n}$, 
we say that it is of type E1 if it only extends to the toric action $\T_{r}$. We say that it is of type E2 if it extends to exactly two toric actions $\T_{r}$ and $\T_{r'}$.
\end{defn}

\begin{cor}\label{cor:NumberOfStrataE1E2} For a Hamiltonian $\Z_n(1,b;r)$ action of type E1 on $(\SSS,\oml)$, the space of invariant almost complex structures $\jzoml$ only intersects the stratum $U_{r}$. For an action of type E2, the space $\jzoml$ intersect the strata $U_r$ and $U_{r'}$, where $r'$ is either equal to $|r-2b|$ or $r-2b +2n$.
\end{cor}
\begin{proof}
By Corollary~\ref{cor:strata-cyclic-case}, the space $\jzoml$ intersects the stratum $U_{r'}$ iff the action $\Z(a,b;r)$ is equivariantly symplectomorphic to some $\Z(a',b';r')$ action. Hence the Corollary follows from Lemma~\ref{lemma:Toric-Extensions-(1,b;r)}.
\end{proof}

\begin{remark}\label{Rmk:ATFcyclic action} 
Using the theory of almost toric fibrations, it is possible to show that the cyclic action $\Z_{r-2}(1,1;r)$ on $(\SSS,\oml)$ with $r>4$ is isomorphic to some action $\Z_{r-2}(a,b;r-4)\subset\T_{r-4}$. However, this isomorphism is only $\Z_{r-2}(1,1;r)$ equivariant and does not extend to any $S^1$-equivariant isomorphism for any circle $S^1\subset\T_{r}$ that contains $\Z_{r-2}(1,1;r)$. In other words, there is no circle action inside the intersection of $\T_{r}$ and $\T_{r-4}$ that contains the $\Z_{r-2}(1,1;r)$ action.  This is in sharp contrast to the cyclic actions discussed above where all the toric extensions are in fact $S^1$-equivariant for some $S^1$ extensions. This action also provides us with a concrete example of a $\Z_{r-2}(1,1;r)$-invariant sphere which is not $S^1$ invariant for any $S^1\subset\T_r$ extending the $\Z_{r-2}(1,1;r)$ action.  
Note that this fact doesn't contradict Lemma~\ref{lemma:Toric-Extensions-(1,b;r)} as the order of the group $\Z_{r-2}(1,1;r)$ is less than $r$.  
\end{remark}

%%%%%%%%%%%%%%%%%%%%%%%%%%%%%%%%%%%%%%%%%%%%%%%%%%%%%%%%%%%%%%%%%%%%%%%%%%%%%%%%
\section{Homogeneity property of \texorpdfstring{$\Z_n(a,b;r)$}{Z\_n(a,b;r)} actions on \texorpdfstring{$(\SSS,\oml$)}{S2xS2} when $r\neq 0$.}\label{Section:HomogeneityProperty}
%%%%%%%%%%%%%%%%%%%%%%%%%%%%%%%%%%%%%%%%%%%%%%%%%%%%%%%%%%%%%%%%%%%%%%%%%%%%%%%%
Recall that a Hamiltonian $\Z_{n}(a,b;r)$ action on $(\SSS,\oml)$ satisfies the \emph{homogeneity property} if each non-empty invariant stratum in the decomposition
\[\jzoml= U_0^{\Z_n}\sqcup U_2^{\Z_n}\sqcup\cdots\sqcup U_{2\ell}^{\Z_n}\]
is homotopically equivalent to the $\Symp_h^{\Z_n}(M_{\lambda})$ orbit of a standard structure $J_{r}$ with stabilizer homotopically equivalent to the subgroup $\Iso^{\Z_{n}}_h(\oml,J_{r})$ of equivariant Kahler isometries. In this section we investigate the homogeneity property of the cyclic actions having at most two toric extensions described in Proposition~\ref{prop:gcdConditionSingleToricExtension} and Lemma~\ref{lemma:Toric-Extensions-(1,b;r)}. We find further conditions on the parameters $\lambda$, $a$, $b$, $r$, and $n$ implying the homogeneity property. This is done by first replacing the action of $\Symp^{\Z_n}_{h}(\SSS,\oml)$ on the stratum $U_{2s}^{\Z_n}$ by the homotopy equivalent action of $\Symp^{\Z_n}_{h}(\SSS,\oml)$ on the space $\mathcal{S}^{\Z_n}_{D_{2s}}$ of invariant curves in class $D_{2s}$, and then by investigating the latter action through a sequence of fibrations (See Lemma~\ref{Lemma:Homotopyequivalencestratawithcurves}). Our analysis closely follow the exposition in Section 3.5.2 of~\cite{ChPin-Memoirs} where the homogeneity property is proven for all Hamiltonian circle actions.

%%%%%%%%%%%%%%%%%%%%%%%%%%%%%%%%%%%%%%%%%%%%%%%%%%%%%%%%%%%%%%%%%%%%%%%%%%%%%%%%
\subsection{The action of $\Symp^{\Z_n}_{h}(\SSS,\oml)$ on the space $\mathcal{S}^{\Z_n}_{D_{2s}}$ for $2s>0$}
%%%%%%%%%%%%%%%%%%%%%%%%%%%%%%%%%%%%%%%%%%%%%%%%%%%%%%%%%%%%%%%%%%%%%%%%%%%%%%%%

\begin{lemma}\label{Lemma:Homotopyequivalencestratawithcurves}
Given $r=2s>0$, the map $U_r^{\Z_n}\to \mathcal{S}^{\Z_n}_{D_{r}}$ that takes an almost complex structure $J$ to the unique $J$-holomorphic curve in class $D_r$ is a homotopy equivalence.
\end{lemma}
\begin{proof}
The proof is exactly as in the proof of Lemma 3.41 in \cite{ChPin-Memoirs}.
\end{proof}
The rest of this subsection is devoted to the proof that the evaluation map $\Symp^{\Z_n}_{h}(\SSS,\oml) \to {\mathcal{S}^{\Z_n}_{D_{2s}}}$ is a fibration for all $\Z_n(a,b;r)$ actions of the following three types: 
\begin{enumerate}[label=(H{\arabic*}), start=0]
\item Actions $\Z_n(a,b;r)$ such that $\gcd(a,n) \neq 1$ as considered in Proposition~\ref{prop:gcdConditionSingleToricExtension} and satisfying the extra condition $r \neq 0$.
\item Actions $\Z_n(1,b;r)$ of type E1 satisfying $2b-r\neq n$ and $b\not\in\{0,r\}$.
\item Actions $\Z_n(1,b;r)$ of type E2 satisfying $2b-r\neq n$ and $b\not\in\{0,r\}$, with the exclusion of $\Z_{2r}(1,b;r)$, $\Z_{2r+2}(1,r+1;r)$, and $\Z_{2r+2}(1,2r+1;r)$. 
\item Actions of type E1 with $b\in\{0,r\}$, that is, $\Z_n(1,0;r)$ and $\Z_n(1,r;r)$ with $n\geq 2\lambda>r\geq2$.
\end{enumerate}

%%%%%%%%%%%%%%%%%%%%%%%%%%%%%%%%%%%%%%%%%%%%%%%%%%%%%%%%%%%%%%%%%%%%%%%%%%%%%%%%
\subsubsection{Actions of type H0:}\label{subsection:Actions of type H0} When $\gcd(a,n)\neq 1\pmod{n}$, the $\Z_n(a,b;r)$ action on an invariant curve $\overline{D}_r$ in the homology class $D_{r}$ is not free, and by Lemma~\ref{LR}, the curve $\overline{D}_r$ is $J$-holomorphic for all invariant almost complex structure $J$. Because $r>0$, the homological self-intersection $D_r\cdot D_r=-r$ is strictly negative, and positivity of intersections implies that the space $\mathcal{S}^{\Z_n}_{D_{r}}$ of $\Z_n$ invariant symplectic spheres in class $D_r$ only contains the curve $\overline{D}_r$. 

\begin{cor}\label{Cor:SinglecurveInclassD2s}
Consider a finite cyclic group action $\Z_n(a,b;r)$ with $a\not\equiv \pm1\pmod{n}$ that leaves invariant a curve of negative self-intersection $-2s$. Then $2s=r$, and the space $\mathcal{S}^{\Z_n}_{D_{r}}$ of all $\Z_n$-invariant symplectic embedded spheres representing the class $D_{r}$ contains a single curve $\overline{D}_{r}$. Consequently when $a\not\equiv \pm1\pmod{n}$, the action map 
\[\Symp^{\Z_n}_h(\SSS,\oml) \to {\mathcal{S}^{\Z_n}_{D_{r}}}\]
is a fibration in a trivial way, and there are weak homotopy equivalences 
\[
\{*\}\simeq\jzoml=\jzoml \cap U_{r}
\xrightarrow{~\simeq~}
\mathcal{S}^{\Z_n}_{D_{r}}=\{\overline{D}_{r}\}.
\]
\qed
\end{cor}

%%%%%%%%%%%%%%%%%%%%%%%%%%%%%%%%%%%%%%%%%%%%%%%%%%%%%%%%%%%%%%%%%%%%%%%%%%%%%%%%
\subsubsection{Actions of types H1 and H2:} For actions of type $\Z_{n}(1,b;r)$, the restriction of the action on the invariant curve $\overline{D}_r$ is effective, so that the space $\mathcal{S}^{\Z_n}_{D_{r}}$ no longer reduces to a single curve. Moreover, it is a priori possible for different curves in class $D_r$ to pass through different pairs of fixed points. If this is the case, then the space $\mathcal{S}^{\Z_2}_{D_{r}}$ is not connected. Similarly, if the action of the centralizer $\Symp_h^{\Z_n}$ permutes isolated fixed points, then $\Symp_h^{\Z_n}$ is not connected.

\begin{exmp}
Consider the action $\Z_2(1,1,2)$. Recall from Example~\ref{example_obstruction} that this action has isolated fixed points $P$, $Q$, $R$, $S$ (as in fig~\ref{hirz}) with weights all equal to $\{1,1\}\mod 2$. We claim that there exists a $\Z_2(1,1;2)$ invariant curve between $Q$ and $S$. Indeed, by Corollary~2.22 in~\cite{ChPin-Memoirs}, there exists an $S^1(1,1;2)$-equivariant symplectomorphism $\phi$ intertwining the action $\Z_2(1,1;2)$ with the product action $\Z_2(1,1,0)$ on $S^2\times S^2$. If $P'$, $Q'$, $R'$, $S'$ denote the isolated fixed points for the $\Z_2(1,1;0)$ action, this symplectomorphism can be chosen such that $\phi(P)= P'$, $\phi(Q)= Q'$, $\phi(R)=S'$, and $\phi(S)=R'$. Finally, we note that there exists a $\Z_2(1,1;0)$ equivariant symplectomorphism $\psi$ that swaps the pairs $P'$ and $Q'$ and $R'$ and $S'$ fixed. Let $\overline{D}$ denote the standard $\T_r$ invariant curve in class $B-F$. Then $\phi^{-1}\psi \phi(\overline{D})$ is the required $\Z_2(1,1;2)$ invariant curve passing through $Q$ and $S$. Observe that the existence of such a curve implies that that the centralizer $\Symp_h^{\Z_{2}}$ and the space of invariant curves $\mathcal{S}^{\Z_2}_{B-F}$ are not connected. Note, however, that the action $\Z_2(1,1;2)$ is not of type H1 nor H2 since $n=2\leq r=2$.\eoe
\end{exmp}

\begin{prop}\label{prop:ExistenceOfGoodFixedPointH1H2}
For a $\Z_n(1,b;r)$ action of type H1 or H2, and for each invariant stratum $U_{2s}^{\Z_n}$, there is a fixed point $p$ such that
\begin{enumerate}
\item all invariant curves in class $D_{2s}$ pass through $p$,
\item $\phi(p)=p$ for all $\phi\in\Symp^{\Z_n}_h(\SSS,\oml)$.
\end{enumerate}
\end{prop}

The proof is broken in a sequence of lemmas.

\begin{lemma}\label{Lemma:SwapFixedPointsSphere}
Let $\Z_n$ act through Hamiltonian diffeomorphisms on the sphere $(S^2,\om)$. Let $p_1$ and $p_2$ denote the two fixed points for the action. Suppose there exists an equivariant symplectomorphism $\phi \in \Symp^{\Z_n}(S^2, \om)$ such that $\phi(p_1) = p_2$ and $\phi(p_2) =p_1$, then n=2. 
\end{lemma}

\begin{proof}
Assume that $n\geq3$. Then by Lemma~\ref{EquivariantSO(3)}, the centralizer subgroup $\Symp^{\Z_n}(S^2, \om)$ is connected. Hence such a $\phi$ cannot exist. For $n=2$, the action of $\Z_2$ is equivalent to the half-turn rotation fixing the north and south poles. Then, any rotation by $\pi$ radians about a perpendicular axis commutes with the $\Z_2$ action and switches the poles. 
\end{proof}

\begin{lemma}\label{lemma:TransitivityDrCurves}
Given a Hamiltonian $\Z_{n}(a,b;r)$ action, let $\Ss(r;p,q)$ be the space of invariant embedded spheres representing the class $D_{r}$ and having a fixed point where the weight of the tangent representation of $\Z(a,b;r)$ is $p$, and the weight of the normal representation is $q$. Then the group $\Symp_{h}^{\Z_{n}}(\SSS,\oml)$ of equivariant symplectomorphisms acts transitively on $\Ss(r;p,q)$.
\end{lemma}
\begin{proof}
This follows from Proposition~4.7 in~\cite{Liat}. See also the proof of Lemma~\ref{lemma:EvaluationFibrationConfigurations} below.
\end{proof}

\begin{lemma}\label{lemma:WeightsAlongDr}
Let $r>0$, and let $C$ be an $\Z_{n}(a,b;r)$ invariant sphere representing the class $D_{r}$. Suppose the ordered weights at one fixed point are $(p,q)$ (modulo $n$), where $p$ is the weight of the tangent representation of $\Z(a,b;r)$ along $C$, and $q$ is the weight of the normal representation. Then the ordered weights at the other fixed point are $(-p, pr+q)$, with $-p$ corresponding to the tangent representation, and $pr+q$ corresponding to the normal representation.
\end{lemma}
\begin{proof}
This follows directly from Table~\ref{table_weights_Z_n} and Lemma~\ref{weightZ_n}.
\end{proof}

\begin{lemma}\label{lemma:gcdImplies0}
Suppose $n\geq 3$, $n>r$, and $\gcd(b,n)=1$. Then $rb=0\pmod{n}$ implies $b=0$.\qed
\end{lemma}

In the next lemma, we use the extra condition $2b-r\neq n$ imposed in the definition of cyclic actions of types H1 and H2.

\begin{lemma}\label{lemma:ConditionsForQandRInvariant}
Let $\Z_{n}(1,b;r)$ be an action of type E1 or E2. Assume also that $2b-r\neq n$ and $b \neq \{0,r\}$. Then any symplectomorphism $\phi$ commuting with the action acts trivially on the fixed points $Q$ and $R$, that is, $\phi(Q)=Q$ and $\phi(R)=R$.
\end{lemma}
\begin{proof}
First recall that an action of type E1 or E2 satisfies $n\geq 2\lambda>r\geq2$ and $2b-r\neq 0$. Together with the conditions $2b-r\neq n$ and $b \neq \{0,r\}$, this implies the modular inequalities $r\neq0\pmod{n}$, $2b-r\neq0\pmod{n}$, and $b \neq \{0,r\}\pmod{n}$. 

For an equivariant symplectomorphism to permute $Q$ and $R$, the sets of weights at $Q$ and $R$ must be the same, that is, $\{1,-b\}=\{-1,r-b\}$. Then either $1=-1\pmod{n}$ and $b=b-r\pmod{n}$, which is excluded, or $1=r-b\pmod{n}$ and $-b=-1\pmod{n}$, which implies $2b-r=0\pmod{n}$, which is also excluded. 

We now compare the weights at $Q$ and $R$ with the weights at $P$ and $S$.

\begin{itemize}
\item $w(Q)=w(P)\iff \{1,-b\}=\{1,b\}$. If $b=-b\pmod{n}$, then $2b=0\pmod{n}$. If $b=1\pmod{n}$, then we must also have $b=-1\pmod{n}$, which is only possible when $n=2$.

\item $w(Q)=w(S)\iff \{1,-b\}=\{-1,b-r\}$. If $1=-1\pmod{n}$, then $n=2$. If $b=1\pmod{n}$ and $b-r=1\pmod{n}$, then $r=0\pmod{n}$.

\item $w(R)=w(P)\iff \{-1,r-b\}=\{1,b\}$. If $-1=1\pmod{n}$, then $n=2$. If $r-b=1\pmod{n}$ and $b=-1\pmod{n}$, then $r=0\pmod{n}$.

\item $w(R)=w(S)\iff \{-1,r-b\}=\{-1,b-r\}$. If $b-r=r-b\pmod{n}$, then $2(b-r)=0\pmod{n}$.
\end{itemize}
We conclude that under the assumptions of the lemma, $w(Q)\neq w(S)$, and $w(Q)=w(P)$ only if $2b=0\pmod{n}$. Similarly, $w(R)\neq w(P)$ and $w(R)=w(S)$ only if $2(b-r)=0\pmod{n}$.

Suppose $2b=0\pmod{n}$ and $w(Q)=w(P)$. If $\gcd(b,n)=1$, then Lemma~\ref{lemma:gcdImplies0} implies that $b=0$. If $\gcd(b,n)\neq 1$, then the invariant fiber $F$ connecting $Q$ and $P$ is $J$-holomorphic for all invariant $J$. In particular, it is invariant under any equivariant symplectomorphism $\phi$ that swaps $Q$ and $P$. Restricting $\phi$ to $F$, we obtain an equivariant symplectomorphism of $S^{2}$ that swaps the north and south poles which, by Lemma~\ref{Lemma:SwapFixedPointsSphere}, is only possible if $n=2$.

Similarly, suppose $2(b-r)=0\pmod{n}$ and $w(R)=w(S)$, If $\gcd(b-r,n)=1$, then Lemma~\ref{lemma:gcdImplies0} implies that $b-r=0\pmod{n}$. If $\gcd(b-r,n)\neq 1$, then Lemma~\ref{lemma:WeightsAlongDr} implies that $n=2$.
\end{proof}

\begin{lemma}\label{lemma:ConditionsForDrCurvesOnlyThroughQandRH1}
Let $\Z_{n}(1,b;r)$ be an action of type H1. Then any invariant curve $C$ in class $D_{r}$ passes through the fixed points $Q$ and~$R$.
\end{lemma}
\begin{proof}
By definition, actions of type H1 only extend to the toric action $\T_r$, so that $\jzoml=U_r^{\Z_n}$. In particular, the unique strata is path connected. Given any curve $C$ in class $D_r$, choose a $J\in U_r^{\Z_n}$ such that $C$ is $J$-holomorphic. Let $J_t$, $0\leq t\leq1$, be a path connecting the standard Hirzebruch structure $J_r$ to $J$. For each $J_t$, there is a unique $J_t$  holomorphic invariant curve $C_t$ in class $D_r$ that varies continuously with $t$. Since $C_0$ passes through $Q$ and $R$, and since the fixed points are isolated, it follows that all the invariant curves $C_t$ must pass through $Q$ and~$R$.
\end{proof}

A similar result holds for actions of type H2. However, the previous argument does not apply since strata of strictly positive codimensions could be disconnected. Instead, we again consider the possible weights at the 4 fixed points. This explains the exclusion of some actions in the definition of cyclic actions of type H2.

\begin{lemma}\label{lemma:ConditionsForDrCurvesOnlyThroughQandR}
Let $\Z_{n}(1,b;r)$ be an action of type H2. Then any $\Z_{n}(1,b;r)$ invariant curve $C$ in class $D_{r}$ passes through the fixed points $Q$ and~$R$.
\end{lemma}
\begin{proof}
Recall that an action of type H2 is of type E2, satisfies $2b-r\neq n$, $b \neq \{0,r\}$, and is not of the form $\Z_{2r}(1,b;r)$, $\Z_{2r+2}(1,2r+1;r)$, or $\Z_{2r+2}(1,r+1;r)$. Suppose a curve $C$ in class $D_r$ passes through $Q$ and $P$. If $\gcd(b,n)\neq 1\pmod{n}$, then the invariant fiber $F$ through $P$ and $Q$ has a non-trivial stabilizer and, by Lemma~\ref{LR}, is holomorphic for every invariant $J$, in particular for a $J$ for which $C$ is also holomorphic. Since the homological intersection is $F\cdot C=1$, positivity of intersections yields a contradiction. Let's assume $\gcd(b,n)=1\pmod{n}$. If the ordered weights along $C$ are the same as the standard curve $\overline{D}_{r}$, then Lemma~\ref{lemma:TransitivityDrCurves} implies that there is an equivariant symplectomorphism sending $C$ to $\overline{D}_{r}$, hence, sending $P$ to $R$. By Lemma~\ref{lemma:ConditionsForQandRInvariant} this is impossible. Consequently, the ordered weights of $C$ at $Q$ must be inverted, that is, the tangent representation at $Q$ along $C$ must have weight $-b$, and the normal representation must have weight $1$. By Lemma~\ref{lemma:WeightsAlongDr}, the ordered weights at $P$ must be congruent to $(b,1-br)$. Since the weights at $P$ are $\{1,b\}$, either $br=0\pmod{n}$, or $b=1\pmod{n}$ and $1-br=b\pmod{n}$. In both cases, Lemma~\ref{lemma:gcdImplies0} implies $b=0$, which is excluded.

Now suppose a curve $C$ in class $D_{r}$ passes through $Q$ and $S$. As before, the weight at $Q$ of the tangent representation must be $-b$, and the weight of the normal representation at $Q$ must be $1$. By Lemma~\ref{lemma:WeightsAlongDr}, it follows that the tangent weight at $S$ is $b$, while the normal weight at $S$ is $1-br$. Since the weights at $S$ are $\{-1,b-r\}$ (modulo $n$), we have either $b=b-r$ and $1-br=-1$, or $b=-1$ and $1-br=b-r$. In the former case, we obtain $r=0\pmod{n}$, which is excluded. In the latter case, we must have $2(r+1)=0\pmod{n}$. Because $2\leq r<n$, this implies $n=2r+2$ and $b=-1\equiv 2r+1\pmod{n}$, so that $2b-r=3r+2> 2r+2=n$ which is excluded.

Similarly, suppose a curve $C$ in class $D_{r}$ passes through $R$ and $S$. If $\gcd(b-r,n)\neq 1$, then the invariant fiber through $R$ and $S$ is $J$-holomorphic for every invariant $J$ and we reach a contradiction with positivity of intersections. If $\gcd(b-r,n)=1$, then the tangent weight along $C$ at $R$ must be $r-b$, while the normal weight must be $-1$. By Lemma~\ref{lemma:WeightsAlongDr}, the tangent weight along $C$ at $S$ must be $b-r$ and the normal weight must be $(r-b)r-1$. Since $w(S)=\{-1,b-r\}$, it follows that either $r=0\pmod{n}$, which is excluded, or that $r-b=0\pmod{n}$. In the latter case, Lemma~\ref{lemma:gcdImplies0} implies that either $r=0$ or $r-b=0$, which are also excluded. 

Now suppose a curve $C$ in class $D_{r}$ passes through $R$ and $P$. By Lemma~\ref{lemma:TransitivityDrCurves} and Lemma~\ref{lemma:ConditionsForQandRInvariant}, the tangent weight along $C$ at $R$ must be $r-b$, while the normal weight must be $-1$. By Lemma~\ref{lemma:WeightsAlongDr}, the tangent weight along $C$ at $P$ must be $b-r$ and the normal weight must be $(r-b)r-1$. Since $w(P)=\{1,b\}$, it follows that $b-r=b$ or $b-r=1$ and $(r-b)r-1=b$. In the former case we get $r=0\pmod{n}$. In the latter case we must have $b=r+1$ and $2(r+1)=0\pmod{n}$. Since $2\leq r<n$, this implies $n=2r+2$. As before, these two possibilities are excluded.

Finally, suppose a $D_{r}$ curve $C$ passes through $P$ and $S$. Assume that its tangent and normal weights at $P$ and $S$ coincide with the tangent and normal weights of the invariant $D_{-r}$ curve passing through these points. Then Lemma~\ref{lemma:WeightsAlongDr} implies that $2r=0\pmod{n}$, which is excluded.
Now suppose instead that the tangent weight at $P$ is $b$ and the normal weight is $1$. Then, by Lemma~\ref{lemma:WeightsAlongDr}, the tangent and normal weights at $S$ are respectively $-b$ and $br+1$. Since the weights at $S$ are $\{-1,b-r\}$, we either have $b=1$ and $br+1=b-r$ or $br+1=-1$ and $-b=b-r$. In the former case, this implies $2r=0\pmod{n}$ which is excluded. In the latter case, $r=2\pmod{n}$ and $b=1\pmod{n}$. Then $2b-r=0$, which is also excluded.
\end{proof}

\begin{proof}[Proof of Proposition~\ref{prop:ExistenceOfGoodFixedPointH1H2}] For the invariant stratum $U_{2s}$, take $p$ to be one of the fixed points $Q$ or $R$ of the toric extension~$\T_{2s}$ and apply lemmas~\ref{lemma:ConditionsForQandRInvariant},~\ref{lemma:ConditionsForDrCurvesOnlyThroughQandRH1}, ~\ref{lemma:ConditionsForDrCurvesOnlyThroughQandR}.
\end{proof}

\begin{lemma}\label{lemma:DrCurvesForH3Action} Let $\Z_{n}(1,b;r)$ be an action of type E1 with $b=0$ or $b=r$.
\begin{enumerate}
\item When $b=0$ then any $\Z_{n}(1,b;r)$ invariant curve $C$ in class $D_{r}$ passes through the fixed point $R$. Consequently, any symplectomorphism $\phi$ commuting with the action acts trivially on the fixed $R$, that is, $\phi(R)=R$.
\item When $b=r$ then any $\Z_{n}(1,b;r)$ invariant curve $C$ in class $D_{r}$ passes through the fixed point $Q$. Consequently, any symplectomorphism $\phi$ commuting with the action acts trivially on the fixed $Q$, that is, $\phi(Q)=Q$.
\end{enumerate}
\end{lemma}
\begin{proof}
We prove the statement for $b=0$. The case $b=r$ is analogous. We mimick the proof of Lemma~\ref{lemma:ConditionsForDrCurvesOnlyThroughQandRH1}. By Lemma~\ref{lemma:Toric-Extensions-(1,b;r)} we know that these actions extend to only one torus $\T_r$. Consequently, $U^{\Z_n}_r = \jzoml$ is path-connected. Given any curve $C$ in class $D_r$, choose a $J\in U_r^{\Z_n}$ such that $C$ is $J$-holomorphic. Let $J_t$, $0\leq t\leq1$, be a path connecting the standard Hirzebruch structure $J_r$ to $J$. For each $J_t$, there is a unique $J_t$  holomorphic invariant curve $C_t$ in class $D_r$ that varies continuously with $t$. Since $C_0$ passes through $Q$ and $R$, and as $R$ is an isolated fixed point, any $J_t$ holomorphic curve in class $D_r$ has to pass through $R$ as well. Furthermore, the other fixed point of $C$ must lie in the pointwise invariant fiber passing through $Q$. Consequently, in order to show that $\phi(R)=R$ for all equivariant symplectomorphisms $\phi$, it is enough to rule out the possibility that $\phi(R) = S$. But if $\phi$ is an equivariant symplectomorphism such that $\phi(R) = S$, then $\phi(C_0)$ is a $D_r$ curve that doesn't pass through $R$, contradicting the first part of the claim. 
\end{proof}

Let $J$ be an invariant almost complex structure in the stratum $U_{2s}^{\Z_{n}}$. Let $\overline{D}_{2s}$ be the unique invariant $J$-curve in class $D_{2s}$. Choose a fixed point $p\in\overline{D}_{2s}$ and let $\overline{F}$ be the unique invariant $J$-holomorphic fiber through $p$. If the two weights at $p$ are distinct, these two symplectic curves intersect orthogonally. Conversely, each orthogonal configuration of such invariant curves through $p$ is $J$-holomorphic for some invariant $J\in U_{2s}$. Let $\mathcal{C}(D_{2s}\vee F,p)^{\Z_{n}}$ be the space of all such orthogonal configurations.

\begin{lemma}\label{lemma:EvaluationFibrationConfigurations} 
Let $\Z_{n}(1,b;r)$ be an action of type H1, H2, or H3, which admits a toric extension $\T_{2s}$. There is a fixed point $p$ satisfying the following conditions
\begin{itemize}
\item $\phi(p)=p$ for all symplectomorphism commuting with $\Z_{n}$,
\item all invariant curves in class $D_{2s}$ pass through $p$, and
\end{itemize}
Moreover, the evaluation map at a standard configuration $\overline{D}_{2s}\vee \overline{F}$ through $p$
\[\Symp^{\Z_{n}}_h(\SSS,\oml) \to \mathcal{C}(D_{2s}\vee F,p)^{\Z_{n}}\]
is a Serre fibration.
\end{lemma}
\begin{proof}
Consider an action of type H1 or H2. By Proposition~\ref{prop:ExistenceOfGoodFixedPointH1H2}, the toric fixed points $Q$ and $R$ satisfy the first two conditions. For an action of type H3, Lemma~\ref{lemma:DrCurvesForH3Action} implies the first two conditions. 

We now show that the action is transitive. Given an invariant configuration $C\vee A \in \mathcal{C}(D_{2s}\vee F,p)^{\Z_{n}}$, using Theorem C.2 in \cite{ChPin-Memoirs}, there exists an equivariant Hamiltonian isotopy $\eta$ of $(\SSS,\oml)$ supported in a small neighbourhood of the fixed point $p$ such that $\eta(C)$ and $\eta(A)$ intersect $\oml$-orthogonally. By abuse of notation let us denote the deformed curves by $C$ and $A$ as well. Using the equivariant symplectic neighbourhood theorem implies we can find an invariant neighbourhood $V$ of $C \vee A$, an invariant neighbourhood $V'$ of $\overline{D}_{2s} \cup \overline{F}$, and an equivariant symplectomorphism $\alpha:V \to V'$. We claim that $\alpha$ can be extended to an ambient diffeomorphism $\beta$ of $\SSS$. Assume this for the moment. By construction, the pullback form $\om_\beta:=\beta^*\om_\lambda$ is invariant under the conjugate action $\beta^{-1}\rho\beta$.

We observe that the complement of the standard configuration $\overline{D}_{2v} \vee\overline{F}$ in $W_{2s}$ is symplectomorphic to $\C^2$ with the symplectic form $\om_f:=\frac{1}{2\pi}\del\delbar f$ where $f=\log\left(\left(1+||w||^2\right)^{\lambda}\left(1+||w||^{4k}+||z||^2\right)\right)$, see~\cite[Lemma~3.5]{abreu}. Under this identification, the standard $\T_{2s}$ action on $W_{2s}$ becomes linear on $\C^2$. It follows that near infinity, the form $\om_\beta$ is equal to $\om_f$ and that the action $\beta^{-1}\rho\beta$ is linear. By Proposition 3.27 in \cite{ChPin-Memoirs}, we get that there is an equivariant symplectomorphism $\gamma$ that is equal to the identity near infinity, and that identifies $(\C^2,\om_f, \rho)$ with $(\C^2,\beta^*\om_\lambda,\beta^{-1}\rho\beta)$. By construction, the equivariant symplectomorphism $\phi := \gamma \circ \beta$ takes the configuration $C\vee A$ to $\overline{D}_{2s}\vee\overline{F}$. 

It remains to see that the local diffeomorphism $\alpha$ can be extended to an ambient diffeomorphism $\beta$ of $\SSS$. By the Isotopy Extension Theorem (see~\cite[Theorem~1.4, p.180]{Hi}), it suffices to show that any two configurations of embedded spheres in classes $\overline{D}$ and $F$ intersecting transversely and positively are isotopic. In turn, this follows from the fact that any two $F$-foliations corresponding to two almost-complex structures $J$ and $J'$ are diffeotopic, and that any two sections of the product $\SSS$ are diffeotopic through sections iff they are homotopic. This shows that $\Symp^{\Z_{n}}_h(\SSS,\oml)$ acts transitively on $\mathcal{C}(D_{2s}\vee F,p)^{\Z_{n}}$.

To prove the homotopy lifting property, consider any family of maps $\gamma: D^n \times [0,1] \rightarrow \mathcal{C}(D_{2s}\vee F,p)^{\Z_{n}}$ from a $n$ dimensional disk $D^n$ to $\mathcal{C}(D_{2s}\vee F,p)^{\Z_{n}}$, and choose a lift $\overline{\gamma_0}:D^n \rightarrow \Symp^{\Z_{n}}_h(\SSS,\oml)$ of $\gamma_0$. Since the complement of a configuration is contractible, the equivariant version of Banyaga's Extension Theorem for families (\cite[Theorem A.9]{ChPin-Memoirs}) implies that there exists a lift $\overline{\gamma}: D^n \times [0,1] \rightarrow \Symp^{\Z_{n}}_h(\SSS,\oml) $ extending $\overline{\gamma_0}$. 
\[
\begin{tikzcd}
D^n \times \{0\} \arrow[d,hookrightarrow] \arrow[r,"\overline{\gamma_0}"]    &\Symp^{\Z_{n}}_h(\SSS,\oml) \arrow[d,"\theta"] \\
    D^n \times [0,1]\arrow[r,"\gamma"] \arrow[ur,dashrightarrow,"\exists ~ \overline{\gamma}"] &\mathcal{C}(D_{2s}\vee F,p)^{\Z_{n}}
\end{tikzcd}
\]
Alternatively, one can apply the equivariant Gromov-Auroux Lemma \cite[Lemma A.10]{ChPin-Memoirs} to show the existence of the lift  $\overline{\gamma}$. In both cases, this concludes the proof.
\end{proof}

\begin{cor}\label{corFibrationSpaceOfCurves} 
Let $\Z_{n}(1,b;r)$ be an action of type H1, H2, or H3 which admits a toric extension $\T_{2s}$. Then $\Symp^{\Z_{n}}_h(\SSS,\oml)$ acts transitively on ${\mathcal{S}^{\Z_{n}}_{D_{2s}}}$. Moreover, let $\overline{D}\in {\mathcal{S}^{\Z_{n}}_{D_{2s}}}$ be an invariant symplectic sphere in the class $D_{2s}$, then the evaluation map
\begin{align*}
\theta: \Symp^{\Z_{n}}_h(\SSS,\oml) &\to {\mathcal{S}^{\Z_{n}}_{D_{2s}}} \\
\phi &\mapsto \phi(\overline{D})
\end{align*}
is a Serre fibration with fibre over $\overline{D}$ given by
\[\Stab(\overline{D}):= \left\{ \phi \in \Symp^{\Z_{n}}(\SSS,\oml)~|~ \phi(\overline{D}) = \overline{D}\right\}\]\qed
\end{cor}

\begin{proof}
The evaluation map is transitive by Corollary~\ref{corFibrationSpaceOfCurves}. The homotopy lifting property follows from the equivariant Gromov-Auroux Lemma as in the proof of Lemma~\ref{lemma:EvaluationFibrationConfigurations}. Alternatively, one can also note that the action map factors through the restriction map 
\[\mathcal{C}(D_{2s}\vee F,p)^{\Z_{n}}\to\mathcal{S}^{\Z_{n}}_{D_{2s}}\]
which is itself a fibration. To see this, note that the restriction map fits into a commuting diagram
\[
\begin{tikzcd}
U_{2s}^{\Z_{n}}\arrow[r,"f_1"]\arrow[rd,"f_2"] & \mathcal{C}(D_{2s}\vee F,p)^{\Z_{n}} \arrow[d]\\
& \mathcal{S}^{\Z_{n}}_{D_{2s}}
\end{tikzcd}
\]
where the maps $f_1$ and $f_2$ are fibrations. Observe that the map $f_1$ is well defined because the weights at the chosen fixed point $p$ are not equal. Hence, for any choice of invariant almost complex structure $J \in U_{2s}^{\Z_{n}}$, the unique invariant $J$-holomorphic curve $C$ in class $D_{2s}$ intersects the unique invariant $J$-holomorphic fiber through $p$ $\oml$-orthogonally. 
\end{proof}

\subsection{The homogeneity property for actions of types H0, H1, H2, and H3}

Let $\Z_{n}(1,b;r)$ be an action of type H0, H1, H2, or H3 which admits a toric extension $\T_{2s}$.

We now proceed with the proof that the space of $\Z_n(a;b,r)$ invariant almost complex structures $\U_{2s}^{\Z_{n}}$ is homotopically equivalent to the orbit of the almost complex structure $J_{2s}\in U_{2s}^{\Z_{n}}$ under the action of the equivariant symplectomorphism group $\Symp^{\Z_n}(\SSS,\oml)$. The arguments used in~\cite[Section~3.5.2]{ChPin-Memoirs} in the case of Hamiltonian circle actions work mutatis mutandis for $\Z_n$ actions. Hence we refer the reader to this paper for most of the proofs and only give details for the homotopy equivalences that differ from the $S^1$ case.\\

Starting with the fibration
\begin{equation}\label{eq:Fibration1}
\Stab^{\Z_n}(\overline{D}) \to \Symp^{\Z_n}_{h}(\SSS,\oml) \to \mathcal{S}^{\Z_n}_{D_{2s}} \simeq U_{2s}^{\Z_{n}},
\end{equation}
the restriction of an element of $\phi\in\Stab^{\Z_n}(\overline{D})$ to the curve $\overline{D}$ defines a map $\Stab^{\Z_n}(\overline{D}) \to \Symp^{\Z_n}(\overline{D})$ to the group of symplectomorphisms of $\overline{D}$ commuting with the restricted action.
\begin{lemma}\label{EquivariantSO(3)}
Given a $\Z_{n}(a,b;r)$ action on $\SSS$, the group $\Symp^{\Z_n}(\overline{D})$ of symplectomorphisms of the invariant sphere $\overline{D}$ that are equivariant with respect to the restricted action is homotopy equivalent~to
\begin{itemize}
\item $\SO(3)$ for the action $\Z_n(0,\pm 1;r)$,
\item $S^{1}\times \Z_{2}$ for actions of the type $\Z_{2a}(a,b;r)$,
\item $S^1$ for all other $\Z_n(a,b;r)$ actions.
\end{itemize}
\end{lemma}
\begin{proof}  
The restriction of the $\Z_{n}(a,b;r)$ action to the curve $\overline{D}$ is trivial if, and only if, $a=0$, in which case the only effective $\Z_{n}$ action on $\SSS$ is $\Z_n(0,\pm 1;r)$. Then we clearly have $\Symp^{\Z_n}(\overline{D})=\Symp(\overline{D})\simeq\SO(3)$.

Now suppose $\Z_{n}(a,b;r)$ acts non-trivially on $\overline{D}$. After identifying $\overline{D}$ with $S^2$, we can assume that $\Z_{n}$ acts by isometries fixing the two poles. Since the stabilizer of a generic point in $\overline{D}$ is a cyclic group of order $\gcd(a,n)$, the image of $\Z_{n}(a,b;r)$ in $\SO(3)$ is a cyclic group of order $k=n/\gcd(a,n)$. Let $\SO(3)^{\Z_k}$ denote the centraliser of $\Z_k$ inside $\SO(3)$. Given a $\Z_k$ equivariant symplectomorphisms $\psi$ of the sphere, consider the graph $\tilde \psi$ of $\psi$ i.e
\begin{align*}
    \tilde \psi: S^2 &\rightarrow S^2 \times S^2 \\
    z &\mapsto (z,\psi(z))
\end{align*}
Choose a $\Z_k$ equivariant metric for the product $\Z_n$ action on $\SSS$ coming from the $\Z_n$ action on $S^2$. Then by Theorem C, Corollary C and Corollary 4.1 in \cite{wang}, the mean curvature flow with respect to this equivariant metric gives us a canonical homotopy of $\psi$ to an element inside $\SO(3)^{\Z_k}$. Since this homotopy is the identity on all the elements of $\SO(3)^{\Z_k}$,  it is a deformation retraction of $\Symp^{\Z_k}(S^2)$ onto $\SO(3)^{\Z_k}$.

Any element $\phi$ of $\SO(3)^{\Z_k}$ either fixes the poles or permute them. Conversely, any element of $\Z_{k}$ must fix or permute the two fixed points of $\phi$. This implies that $\SO(3)^{\Z_k}$ is isomorphic to $S^1$ when $k\neq 2$, while it is isomorphic to $S^{1}\times \Z_{2}$ for $k=2$. The rest of the claim follows easily.
\end{proof}

The extension of Hamiltonian isotopies implies that the restriction map $\Stab^{\Z_n}(\overline{D}) \to \Symp^{\Z_n}(\overline{D})$ is surjective when $\Stab^{\Z_n}(\overline{D})$ is connected. For the actions $\Z_{2a}(a,b;r)$, the $\Z_{2}$ component of $\Stab^{\Z_n}(\overline{D})$ permutes the toric fixed points $Q$ and $R$ and is induced by the non-trivial element of the symplectic normalizer $N(\T_{r})$ described in Example~\ref{example:Type c2 is Kahler when n=2a}. It follows that the restriction map is a fibration   
\begin{equation}\label{eq:Fibration2}
\Fix^{\Z_n}(\overline{D}) \to \Stab^{\Z_n}(\overline{D}) \to  \Symp^{\Z_n}(\overline{D})
\end{equation}
whose fiber $\Fix^{\Z_n}(\overline{D})$ consists of all $\phi\in\Stab(\overline{D})$ that are the identity on $\overline{D}$. 

For each $\phi\in\Fix^{\Z_n}(\overline{D})$, the restriction of the derivative $d\phi$ to $\overline{D}$ defines an equivariant symplectic automorphism of the symplectic normal bundle $N(\overline{D})$. This defines another fibration
\begin{equation}\label{eq:Fibration3}
\Fix^{\Z_n} (N(\overline{D})) \to \Fix^{\Z_n}(\overline{D}) \to  \Aut^{\Z_n}(N(\overline{D}))
\end{equation}
whose fiber is made of symplectomorphisms whose derivatives restrict to the identity on the tangent bundle $T_{\overline{D}}(\SSS)$ along $\overline{D}$. 
\begin{lemma}
The group $\Aut^{\Z_n}(N(\overline{D}))$ of $\Z_n(a,b;r)$ equivariant automorphisms of the symplectic normal bundle $N(\overline{D})$ deformation retracts onto $\Sp(2)\simeq U(1)=S^{1}$.
\end{lemma}
\begin{proof}
The curve $\overline{D}$ is invariant under the extension of the $\Z_n$-action to the circle action $S^1(a,b;r)$. The argument is then identical to the proof of  Lemma~B.1 in~\cite{ChPin-Memoirs}.
\end{proof}
Choose a fixed point $p_{0}$ on ${\overline{D}}$ and let ${\overline{F}}$ be an invariant fiber intersecting ${\overline{D}}$ $\oml$-orthogonally at $p_{0}$. Let $\overline{\mathcal{S}^{\Z_n}_{F,p_0}}$ be the space of unparametrized $\Z_n$-invariant symplectic spheres in the homology class $F$ that coincide with ${\overline{F}}$ in a neighbourhood of $p_0$. Let $\J_{\oml}^{\Z_n}(\overline{D})$ be the contractible space of $\Z_n$-equivariant $\oml$ compatible almost complex structures for which the curve $\overline{D}$ is holomorphic. As explained in~\cite[Lemma~3.45]{ChPin-Memoirs}, there is a natural homotopy equivalence $\overline{\mathcal{S}^{\Z_n}_{F,p_0}}\simeq \J_{\oml}^{\Z_n}(\overline{D})$ and the group $\Fix^{\Z_n}(N(\overline{D}))$ acts transitively on 
$\overline{\mathcal{S}^{\Z_n}_{F,p_0}}$ inducing a fibration
\begin{equation}\label{eq:Fibration4}
\Stab^{\Z_n}(\overline{F}) \cap \Fix^{\Z_n}(N(\overline{D})) \to \Fix^{\Z_n}(N(\overline{D})) \to  \overline{\mathcal{S}^{\Z_n}_{F,p_0}}\simeq \mathcal{J}_{\oml}^{\Z_n}(\overline{D})\simeq \{*\}
\end{equation}
As before, restricting symplectomorphisms and their derivatives to $\overline{F}$ and $N(\overline{F})$ yields two more fibrations
\begin{equation}\label{eq:Fibration5}
\Fix^{\Z_n}(\overline{F}) \to \Stab^{\Z_n}(\overline{F}) \cap \Fix^{\Z_n}(N(\overline{D})) \to  \Symp^{\Z_n}(\overline{F}, N(p_0)) 
\end{equation}
where $\Symp^{\Z_n}({\overline{F}}, N(p_0))$ is the group of equivariant symplectomorphisms of the sphere $\overline{F}$ that are the identity in an open neighborhood of $p_0$ in $\overline{F}$, and
\begin{equation}\label{eq:Fibration6}
\Fix^{\Z_n}(N(\overline{D} \vee \overline{F})) \to \Fix^{\Z_n}(\overline{F}) \to  \Aut^{\Z_n}(N(\overline{D} \vee \overline{F}))
\end{equation}
where $\Fix^{\Z_n}(N(\overline{D} \vee {\overline{F}}))$ is the group of all $\Z_n$-equivariant symplectomorphisms that are the identity in the neighbourhood of the union $\overline{D} \vee {\overline{F}}$, and where 
$\Aut^{\Z_n}(N(\overline{D} \vee \overline{F}))$ is the group of $\Z_n(a,b;r)$ equivariant automorphisms of the symplectic normal bundle $N(\overline{D} \vee \overline{F})$ defined as the symplectic plumbing of the normal bundles $N(\overline{D})$ and $N(\overline{F})$.
\begin{lemma}\label{lemma:ContractibilitySympAwayFromPole}
The group $\Symp^{\Z_n}(\overline{F}, N(p_0))$ of $\Z_n(a,b;r)$ is contractible.
\end{lemma}
\begin{proof}
The group $\Symp^{\Z_n}(\overline{F}, N(p_0))$ is homeomorphic to the group $\Symp^{\Z_{n}}(D^{2},N(\partial))$ of area preserving diffeomorphisms of the disk $D^{2}$ which are the identity on a fixed neighborhood $N(\partial)$ of its boundary, and which commute with a linear cyclic action. Extend each diffeomorphism by the identity outside $D^{2}$ to get an area preserving diffeomorphism of $\R^{2}$. For any such diffeomorphism $\phi$, the Alexander isotopy defined by setting $\phi_{0}=\id$ and $\phi_{t}(x)=t\phi(x/t)$ for $0<t\leq 1$ yields a deformation retraction of $\Symp^{\Z_{n}}(D^{2},N(\partial))$ onto $\{\id\}$.
\end{proof}
\begin{lemma}\label{lemma:ContractibilityAutConfiguration}
The group $\Aut^{\Z_n}(N(\overline{D} \vee \overline{F}))$ of $\Z_n(a,b;r)$ equivariant automorphisms of the symplectic normal bundle $N(\overline{D} \vee \overline{F})$ of a standard orthogonal configuration is weakly contractible.
\end{lemma}
\begin{proof}
Because the configuration $\overline{D} \vee \overline{F}$ is invariant under the circle action $S^1(a,b;r)$, this follows from  \cite[Lemma~B.2]{ChPin-Memoirs}.
\end{proof}
\begin{lemma}\label{lemma:ContractibilityFixConfiguration}
$\Fix^{\Z_n}(N(\overline{D} \vee {\overline{F}}))$ is contractible.
\end{lemma}
\begin{proof}
This follows from an equivariant version of Gromov’s theorem on the contractibility of the group of compactly supported symplectomorphisms of an open polydisc. See~\cite[Theorem 3.29]{ChPin-Memoirs}.
\end{proof}
\begin{cor}
Given any $\Z_{n}(a,b;r)$ action, the group $\Fix^{\Z_n}(\overline{D})$ in fibration~\eqref{eq:Fibration2} is homotopy equivalent to $S^{1}(0,1;r)$, that is, to the $S^{1}$ factor in $\T_{r}$ that fixes the curve $\overline{D}$ pointwise.
\end{cor}
\begin{proof}
This follows from the above fibrations together with Lemmas~\ref{lemma:ContractibilitySympAwayFromPole}, \ref{lemma:ContractibilityAutConfiguration}, and~\ref{lemma:ContractibilityFixConfiguration}. 
\end{proof}

Let $\Iso^{\Z_{n}}(\oml,J_{r})$ be the group of Kahler isometries of the compatible pair $(\oml,J_{r})$ that commute with the $\Z_{n}(a,b;r)$ action. 
\begin{lemma}\label{lemma:EquivariantIsometryGroups}
The group of Kahler isometries commuting with the $\Z_{n}(a,b;r)$ action is isomorphic to
\[\Iso^{\Z_{n}}(\oml,J_{r}) \simeq
\begin{cases}
\SO(3)\times S^{1}(0,1;r)& \text{for the action~} \Z_n(0,\pm 1;r),\\
\T_{r}\times \Z_{2} & \text{for the actions~} \Z_{2a}(a,b;r),\\
\T_{r} & \text{for all other~} \Z_n(a,b;r) \text{~actions}.
\end{cases}
\]
\end{lemma}
\begin{proof}
For $r=2k\geq2$, the Kahler isometry group of $\SSS$ with the metric defined from the pair $(\oml,J_{r})$ is isomorphic to $S^{1}\times\SO(3)$, where the action of $\SO(3)$ is a lift of the standard action on $S^{2}$ to the Hirzebruch surface $W_{r}$ along the projection $W_{r}\to S^{2}$. In particular, $\SO(3)$ acts effectively on the two $J_{r}$ holomorphic sections $s_{0}$ and $s_{\infty}$ in classes $B-kF$ and $B+kF$. The $S^{1}$ factor acts by rotating the fibers fixing $s_{0}$ and $s_{\infty}$. The statement follows from this description together with Lemma~\ref{EquivariantSO(3)}. For more details, see~\cite{AGK} Theorem 3.1 and Section 4.
\end{proof}

\begin{cor}\label{cor:HomotopyTypeStab(D)}
The group $\Stab^{\Z_n}(\overline{D})$ is weakly homotopy equivalent to the equivariant Kahler isometry group $\Iso^{\Z_{n}}(\oml,J_{r})$.
\end{cor}
\begin{proof}
The equivariant Kahler isometry group fits into a commutative ladder of fibrations 
\begin{equation}\label{eq:Fibration2Ladder}
\begin{tikzcd}
\Fix^{\Z_n}(\overline{D}) \arrow[r]{}  & \Stab^{\Z_n}(\overline{D}) \arrow{r} & \Symp^{\Z_n}(\overline{D})\\
S^{1}(0,1;r)\arrow[u,hook,"\simeq"] \arrow{r}{} & \Iso^{\Z_{n}}(\oml,J_{r})\arrow[u,hook]\arrow{r} & \SO(3)^{\Z_n}\arrow[u,hook,"\simeq"]
\end{tikzcd}
\end{equation}
in which the righmost and leftmost inclusions are homotopy equivalences. It follows that the middle inclusion is a weak homotopy equivalence.
\end{proof}

\begin{prop}
Hamiltonian cyclic actions of types H0, H1, and H2 satisfy the homogeneity property. In particular, for every invariant stratum $U_{2s}^{\Z_n}$, there is a homotopy equivalence
\[U_{2s}^{\Z_n}\simeq \Symp^{\Z_n}_{h}(\SSS,\oml)/ \Iso^{\Z_{n}}(\oml,J_{r}). \]
Moreover, these actions are tractable.
\end{prop}
\begin{proof}
The first statement follows from Corollary~\ref{cor:HomotopyTypeStab(D)} and from the fact that for actions of types H0, H1, and H2, and for any invariant curve $\overline{D}$ in class $D_{2s}$, the evaluation map
\begin{equation}\label{eq:Fibration1Again}
\Iso^{\Z_{n}}(\oml,J_{r})\simeq\Stab^{\Z_n}(\overline{D}) \to \Symp^{\Z_n}_{h}(\SSS,\oml) \to \mathcal{S}^{\Z_n}_{D_{2s}} \simeq U_{2s}^{\Z_{n}}
\end{equation}
is a fibration, as shown in Corollary~\ref{corFibrationSpaceOfCurves}. The last statement follows from the fact that these actions admit at most two toric extensions, as shown in Proposition~\ref{prop:gcdConditionSingleToricExtension} and Lemma~\ref{EquivariantSO(3)}. 
\end{proof}

\begin{prop}\label{prop:HomotopyTypeH0Actions}
For a Hamiltonian action $\Z_n(a,b;r)$ of type H0, the symplectic centralizer group $\Symp^{\Z_n}(S^2 \times S^2,\oml)$ is homotopy equivalent to
\begin{enumerate}
\item $S^1 \times \SO(3)$ if $a=0$,
\item $\T_{r}\times\Z_{2}$ if $n=2a\neq0$,
\item $\T_r$ in all other cases.
\end{enumerate}
\end{prop}
\begin{proof}
As $\lambda>1$, any symplectomorphism acts trivially on homology, so that $\Symp^{\Z_n}(S^2 \times S^2,\oml) = \Symp_h^{\Z_n}(S^2 \times S^2,\oml)$. By Corollary~\ref{Cor:SinglecurveInclassD2s}, the base of the fibration
\[
\Iso^{\Z_{n}}(\oml,J_{r})\simeq\Stab^{\Z_n}(\overline{D}) \to \Symp^{\Z_n}(\SSS,\oml) \to \mathcal{S}^{\Z_n}_{D_{r}} \simeq U_{r}^{\Z_{n}}
\]
is the full space $\simeq U_{r}^{\Z_{n}}=\jzoml$, and is thus
contractible. It follows that the inclusion $\Stab^{\Z_n}(\overline{D}) \to \Symp^{\Z_n}(\SSS,\oml)$ is a homotopy equivalence. The statement follows from Lemma~\ref{lemma:EquivariantIsometryGroups}.
\end{proof}

\begin{prop}\label{prop:HomotopyTypeH1Actions}
For a Hamiltonian action $\Z_n(1,b;r)$ of type H1 or H3, the symplectic centralizer group $\Symp^{\Z_n}(S^2 \times S^2,\oml)$ is homotopy equivalent to the torus $\T_r$.
\end{prop}
\begin{proof}
As before, since $\lambda>1$, $\Symp^{\Z_n}(S^2 \times S^2,\oml) = \Symp_h^{\Z_n}(S^2 \times S^2,\oml)$. Since an action of type H1 or H3 extends to a single toric action $\T_r$, the base of the fibration
\[
\Iso^{\Z_{n}}(\oml,J_{r})\simeq\Stab^{\Z_n}(\overline{D}) \to \Symp^{\Z_n}(\SSS,\oml) \to \mathcal{S}^{\Z_n}_{D_{r}} \simeq U_{r}^{\Z_{n}}\]
is again contractible. Since $r\neq0$, and since $n\neq 2=2a$, Lemma~\ref{lemma:EquivariantIsometryGroups} shows that the group $\Iso^{\Z_{n}}(\oml,J_{r})$ is isomorphic to $\T_r$.
\end{proof}

%%%%%%%%%%%%%%%%%%%%%%%%%%%%%%%%%%%%%%%%%%%%%%%%%%%%%%%%%%%%%%%%%%%%%%%%%%%%%%%%
\section{Homogeneity property of \texorpdfstring{$\Z_n(a,b;0)$}{Z\_n(a,b;r)} actions on \texorpdfstring{$(\SSS,\oml$)}{S2xS2}}
%%%%%%%%%%%%%%%%%%%%%%%%%%%%%%%%%%%%%%%%%%%%%%%%%%%%%%%%%%%%%%%%%%%%%%%%%%%%%%%%

In this section, we investigate the homogeneity property for cyclic actions of the form $\Z_n(a,b;0)$ on the product $(\SSS,\oml)$. 
Given an almost complex structure $J$ in the invariant stratum $U_0^{\Z_n}=\jzoml \cap U_0$, recall that there is a two dimensional family of $J$-holomorphic curves in the class $B$, and that given a fixed point $p_0$, there is exactly one of these curves passing through $p_0$. The subgroup $\Symp_{h,p_0}^{\Z_n}(\SSS,\oml)$ of equivariant symplectomorphisms stabilizing a fixed point $p_0$ acts on the space of $\Z_n$-invariant symplectic spheres in class $B$ passing through $p_0$, and we can study this action as done in Section 3.5.3 of \cite{ChPin-Memoirs}. We can then recover the action of the full centralizer group $\phi \in \Symp_h^{\Z_n}(\SSS,\oml)$ on the stratum $U_0^{\Z_n}$ from its action on the set of fixed points. Now recall that the moment polytope of the $\T_0$ action on the product $(\SSS,\om_\lambda)$ is a rectangle with sides of lengths $\lambda\geq 1$ and $1$.

\begin{figure}[H]
\centering       
\begin{tikzpicture}
\node[left] at (0,2) {$Q=(0,1)$};
\node[left] at (0,0) {$P=(0,0)$};
\node[right] at (3,2) {$R= (1,1)$};
\node[right] at (3,0) {$S=(1,0)$};
\node[above] at (1.5,2) {$B$};
\node[right] at (3.015,1) {$F$};
\node[left] at (0,1) {$F$};
\node[below] at (1.5,0) {$B$};
\draw (0,2) -- (3,2) ;
\draw (0,0) -- (0,2) ;
\draw (0,0) -- (3,0) ;
\draw (3,2) -- (3,0) ;
\end{tikzpicture}
\caption{Moment polytope for $\T_0$ action.}
\label{hirZ1lambda>1}
\end{figure}
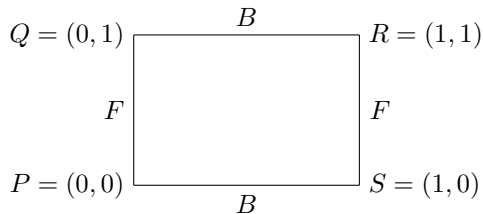

The following result on toric actions plays a crucial role in our arguments. 

\begin{prop}[\cite{ChPin-Memoirs} Propostion 2.16]\label{prop:NormalizersHirzebruch}
Let $N(\T_{r})$ denote the normaliser of the torus $\T_r$ in  the symplectomorphism group $\Symp(\SSS,\oml)$ and let $C(\T_{r})$ denote its centraliser. Then, the Weyl group $W(\T_{r})=N(\T_{r})/C(\T_{r})$ is isomorphic to
\[W(\T_{r})\cong
\begin{cases}
\text{the dihedral group~} D_{4}, & \text{~when~} r=0 \text{~and~} \lambda=1;\\
\text{the dihedral group~} D_{2}\simeq\Z_{2} \times \Z_2, & \text{~when~} r=0 \text{~and~} \lambda>1;\\
\text{the dihedral group~} D_{1}\simeq\Z_{2},& \text{~when~} r\geq 1.
\end{cases}\]
\end{prop}

\begin{cor}\label{Cor:Z_nactionspreservedbyWeylgroup}
Let $W(\T_{0})$ denote the Weyl group of the torus $\T_0$ in $\Symp(\SSS,\oml)$ in the special case $\lambda=1$. Then 
\begin{enumerate}
\item All elements in the normaliser $N(T_{0})$ are equivariant with respect to the action $\Z_2(1,1,0)$.
\item All symplectomorphisms $\phi \in N(\T_{0}) $ in the normaliser satisfying $\phi(P)=S$, $\phi(S)=P$, and $\phi(Q)=R$, $\phi(R)=Q$ are equivariant with respect to $\Z_n(a,b;0)$ actions with $2a \equiv 0\pmod{n}$.
\item All symplectomorphisms $\phi \in N(\T_{0}) $ in the normaliser satisfying $\phi(P)=Q$, $\phi(Q)=P$, $\phi(R)=S$ and $\phi(S)=R$ are equivariant with respect to $\Z_n(a,b;0)$ actions with $2b \equiv 0\pmod{n}$. 
\item All symplectomorphisms $\phi \in N(\T_{0}) $ in the normaliser satisfying $\phi(P)=R$, $\phi(Q)=Q$, $\phi(R)=P$ and $\phi(S)=S$ are equivariant with respect to $\Z_n(1,n-1;0)$. Also, these symplectomorphisms take the class $B$ to the class $F$ and vice-versa.
\item All symplectomorphisms $\phi \in N(\T_{0}) $ in the normaliser satisfying $\phi(P)=P$, $\phi(Q)=S$, $\phi(R)=R$ and $\phi(S)=Q$ are equivariant with respect to $\Z_n(n-1,1;0)$. Also, these symplectomorphisms take the class $B$ to the class $F$ and vice-versa.
\end{enumerate}
\end{cor}
\begin{proof}
Let $\Delta_0\subset\mathfrak{t}_0^*\simeq\R^2$ be the moment polytope of $\T_0$ centered at the origin, and let $D_4$ be its symmetry group. There is an orthogonal representation $\rho: D_4 \rightarrow  \Aut(\mathfrak{t}_0^*)\simeq\AGL(2,\Z)$ generated by the rotation of angle $\pi/2$ and the reflection along the horizontal axis.  
As explained in Section~2.2.2 of~\cite{ChPin-Memoirs}, each symmetry of $\Delta_0$ corresponds to an element in the normalizer $N(\T_0)$.  
The results follow from explicit computations of how the vector $(a,b) \pmod{n}$ transforms under the $D_4$ action. 
\end{proof}

We now find simple conditions on the parameters $a$ and $b$ ensuring the existence of a fixed point $p_0$ such that $\phi(p_0) = p_0$ for all $\phi \in \Symp_{h}^{\Z_n}(\SSS,\oml)$. Recall that any fixed point $p_0$ whose weights are different from the weights at all other fixed points have this property. A case-by-case inspection of the weights given in Table~\ref{table_weights_Z_n} gives sufficient conditions under which some fixed points have unique weights.

\begin{lemma}\label{lemma:ConditionsFixedPointModN}
Fix a $\Z_n(a,b;0)$ action. The fixed points $P$ and $R$ have unique weights if
\[
2a \not\equiv 0\pmod{n},~ 2b \not\equiv 0\pmod{n}, \text{~and~} a \not\equiv -b\pmod{n}.
\]
Similarly, the fixed points $Q$ and $S$ have unique weights if  
\[
2a \not\equiv 0\pmod{n},~ 2b \not\equiv 0\pmod{n}, \text{~and~} a \not\equiv b\pmod{n}.
\]
\end{lemma}
 
\begin{cor}\label{Cor:ChenWilczinskiInequalityU_01}
Fix an action of the form $\Z_n(a,b;0)$. Under the assumptions $2a \not\equiv 0\pmod{n}$ and $2b \not\equiv 0\pmod{n}$,  

either $\phi(R)=R$ or $\phi(Q)=Q$ for all $\phi \in \Symp_h^{\Z_n}(\SSS,\oml)$.  
\end{cor}
\begin{proof}
If $a\equiv-b\pmod{n}$ and $a\equiv b\pmod{n}$, then $2a\equiv 2b\equiv 0\pmod{n}$.
\end{proof}

It is now convenient to subdivide the set of $\Z_n(a,b;0)$ actions into the following four types:
\begin{enumerate}[label=(Z{\arabic*}), start=0]
\item $\Z_n(a,b;0)$ such that $\gcd(a,n) \neq 1$ and $\lambda > 1$.
\item $\Z_n(a,b;0)$ actions such that $2a \not\equiv 0\pmod{n}$ and  $2b \not\equiv 0\pmod{n}$ with $\lambda =1$.
\item $\Z_n(a,b;0)$ actions such either  $2a \equiv 0\pmod{n}$ or  $2b \equiv 0\pmod{n}$, with $n\neq 2$ and $\lambda =1$.
\item $\Z_2(a,b;0)$ with $\lambda=1$.
\end{enumerate}
Note that the case $\Z_n(a,b;0)$ such that $\gcd(a,n) \neq 1$ and $\lambda =1$ is subsumed in Z1,Z2,Z3 actions. We restrict ourselves to these subfamilies for the following reasons. For Z0 actions, the $\gcd(a,n)\neq 1$ condition and the positivity of intersections imply that there is a single toric extension of $\Z_n(a,b;0)$ independently of the value of $\lambda\geq 1$. For the other types of actions, we cannot determine the  toric extensions of the cyclic action unless $\lambda=1$, in which case the only toric extension is $\T_0$. The differences between Z1, Z2, and Z3 actions is the way $\Symp_{h}^{\Z_n}(\SSS,\oml)$ acts on fixed points. In the case of Z1 actions, by Corollary~\ref{Cor:ChenWilczinskiInequalityU_01}, there always exists a fixed point $p_0$ that is left invariant under all symplectomorphisms $\phi \in \Symp^{\Z_n}_{h}(\SSS,\oml)$. On the other hand, for actions of type Z2, we shall explicitly produce symplectomorphisms that swap all fixed point. Finally, actions of type Z3, that is, $\Z_2$ actions on the product $(\SSS,\om_1)$ with factors of equal areas, admit discrete symmetries that are not present for other cyclic actions.

%%%%%%%%%%%%%%%%%%%%%%%%%%%%%%%%%%%%%%%%%%%%%%%%%%%%%%%%%%%%%%%%%%%%%%%%%%%%%%%%
\subsection{Actions of type Z1}\label{Subsection:F_0action}

\begin{lemma}\label{Lemma:CurvethroughQRTypeF_0}
Fix an $\Z_n(a,b;0)$ action of type Z1. Without loss of generality assume that $R$ is fixed by all $\phi \in \Symp_h^{\Z_n}(\SSS,\oml)$. Then every curve in class $B$ that passes through $R$ also passes through $Q$.
\end{lemma}
\begin{proof}
Let $\overline{B}_{QR}$ denote the standard curve in class $B$ that passes through $Q$ and $R$. Assume without loss of generality that there exists a curve $\overline{B}$ that passes through $R$ and $P$. Choose an invariant almost complex structure $J \in \jzoml$ such that $\overline{B}_{QR}$ is invariant under $J$ and an invariant almost complex structure $J' \in \jzoml$ such that $\overline{B}$ is invariant under $J'$. As both $\overline{B}_{QR}$ and $\overline{B}$, pass through the point $R$ and as $B \cdot B =0$, we know by positivity of intersections that $J \neq J'$. As $\lambda =1$, $\jzoml$ is contractible and hence we can choose a smooth path $J_t \in \jzoml$ connecting $J$ and $J'$. Further for each $t$, there exists an invariant curve $\overline{B}_t$ in class $B$ passing through $R$. By equivariant symplectic isotopy theorem, there exists a family of  equivariant symplectomorphisms $\psi_t \in \Symp^{\Z_n}_0(\SSS,\oml)$ $0\leq t\leq 1$ such that $\psi_0 = \id$ and $\psi_t(\overline{B_{QR}}) = \overline{B}_t$ for all $t$. This is a contradiction as for $\Z_n$ action of type Z1, the fixed points $Q$ and $R$ are isolated fixed points but the action of $\psi_t \in \Symp^{\Z_n}_0(\SSS,\oml)$ on the space of fixed points takes $Q$ to $R$.
\end{proof}

\begin{remark}\label{rmk:FixedRQforgenericF_0}
Note that if  $a\equiv b\pmod{n}$ then $a=b=1$ (as $\gcd(a,b)=1$, $0\leq a<n$, and $0 \leq b < n$). Similarly if $a\equiv -b~\pmod{n}$ then $a=1$ and $b=n-1$, or $a=n-1$ and $b=1$. Consequently, for actions of type Z1 which are not of the form $\Z_n(1,1;0)$, $\Z_n(1,n-1;0)$, or $\Z_n(n-1,1;0)$, both the points $Q$ and $R$ have unique weights and hence are fixed by all equivariant symplectomorphisms. 
\end{remark}

Fix a $\Z_n(a,b;0)$ action such that the conditions $2a \not\equiv 0\pmod{n}$ 
and  $2b \not\equiv 0\pmod{n}$. Let $p_0$ denote a fixed point with unique weight (whose existence is guaranteed by Corollary~\ref{Cor:ChenWilczinskiInequalityU_01}). Then, by analogous arguments to the ones given in Section~3.5.3.1 of ~\cite{ChPin-Memoirs} we can prove that there is a sequence of evaluation fibrations
\[\Stab^{\Z_n}(B_{p_0}) \to \Symp^{\Z_n}_{h}(\SSS,\oml) \longtwoheadrightarrow \mathcal{S}^{\Z_n}_{B,p_0} \mathbin{\textcolor{blue}{\xrightarrow{\text{~~~$\simeq$~~~}}}}  U_{0}^{\Z_n} = \jzoml \simeq \{*\}\]
\[\Fix^{\Z_n}(B_{p_0}) \to \Stab^{\Z_n}_{p_0}(B_{p_0}) \longtwoheadrightarrow  \Symp^{\Z_n}(B_{p_0}) \mathbin{\textcolor{blue}{\xrightarrow{\text{~~~$\simeq$~~~}}}} S^1\]
\[\Fix^{\Z_n} (N(B_{p_0})) \to \Fix^{\Z_n}(B_{p_0}) \longtwoheadrightarrow  \Aut^{\Z_n}(N(B_{p_0})) \mathbin{\textcolor{blue}{\xrightarrow{\text{~~~$\simeq$~~~}}}} S^1 \]
\[\Stab^{\Z_n}(\overline{F}) \cap \Fix^{\Z_n}(N(B_{p_0})) \to \Fix^{\Z_n}(N(B_{p_0})) \longtwoheadrightarrow  \overline{\mathcal{S}^{\Z_n}_{F,p_0}} \mathbin{\textcolor{blue}{\xrightarrow{\text{~~~$\simeq$~~~}}}} \mathcal{J}^{\Z_n}(B_{p_0})\]
\[\Fix^{\Z_n}(\overline{F}) \to \Stab^{\Z_n}(\overline{F}) \cap \Fix^{\Z_n}(N(B_{p_0})) \longtwoheadrightarrow  \Symp^{\Z_n}(\overline{F}, N(p_0))  \mathbin{\textcolor{blue}{\xrightarrow{\text{~~~$\simeq$~~~}}}} \left\{*\right\}\]
\[\left\{*\right\} \mathbin{\textcolor{blue}{\xleftarrow{\text{~~~$\simeq$~~~}}}} \Fix^{\Z_n}(N(B_{p_0} \vee \overline{F})) \to \Fix^{\Z_n}(\overline{F}) \longtwoheadrightarrow  \Aut^{\Z_n}(N(B_{p_0} \vee \overline{F}))  \mathbin{\textcolor{blue}{\xrightarrow{\text{~~~$\simeq$~~~}}}} \left\{*\right\}\]
where the subscript $p_0$ denote either diffeomorphisms that fix $p_0$ or space of curves passing through $p_0$ (see Section~\ref{not} for a complete description of the notation). As before, this proves homogeneity.

\begin{thm}\label{thm:Z_nIntersectionWithZeroStratum}
For a $\Z_n(a,b;0)$ action of type Z1, we have homotopy equivalences
\[\Symp_h^{\Z_n}(\SSS,\oml)/\T_0 \simeq U_0^{\Z_n}\simeq *.\]
Consequently, $~\Symp^{\Z_n}_h(\SSS,\oml) \simeq \T_0$.\qed
\end{thm}
Furthermore,
\begin{thm}\label{Thm:FullhomotopytypeF_0}
For a $\Z_n(a,b;0)$ action of type Z1, the homotopy type of the entire symplectomorphism group $\Symp^{\Z_n}(\SSS,\oml)$ is one of the following:
\begin{enumerate}
\item $\Symp^{\Z_n}(\SSS,\oml) \simeq \T_0$ when $\Z_n(a,b;0)$ is not one of the following action $ \Z_n(1,1;0)$, $\Z_n(1,n-1;0)$, $\Z_n(n-1,1;0)$.
    
\item $\Symp^{\Z_n}(\SSS,\oml) \simeq \T_0 \times \Z_2$ for the actions $\Z_n(1,1;0)$, $\Z_n(1,n-1;0)$, and $\Z_n(n-1,1;0)$.
\end{enumerate}
\end{thm}
\begin{proof}
     Let  $ \Aut_{c_1,\oml}(H^2(\SSS))$ denote the space of automorphisms of $H^2(\SSS)$ that preserve the first Chern class and the symplectic form. Consider the left-exact sequence 
\begin{equation*}
   1\to \Symp_h^{\Z_n}(\SSS,\oml)\to  \Symp^{\Z_n}(\SSS,\oml)\to \Aut_{c_1,\oml}(H^2(\SSS))
\end{equation*}

In the case of $(\SSS,\oml)$ with $\lambda = 1$, the group $\Aut_{c_1,\oml}^{\Z_n}\left(H_2(\SSS,\Z)\right)$ is equal to $\Z_2$ and is generated by the symplectomorphism that swaps the two $S^2$ factors. For actions of the form $\Z_n(1,1;0)$, $\Z_n(1,n-1;0)$, and $\Z_n(n-1,1;0)$, by points (1), (4), and (5) in Corollary ~\ref{Cor:Z_nactionspreservedbyWeylgroup}, there exists an equivariant symplectomorphism that swaps the classes $B$ and $F$. Hence for these actions the above short exact sequence is also right-exact and we get that $\Symp^{\Z_n}(\SSS,\oml) \simeq \T_0 \times \Z_2$. For actions of type Z1  that are not of the form $ \Z_n(1,1;0)$, $\Z_n(1,n-1;0)$, or $\Z_n(n-1,1;0)$, we claim that there is no $\Z_n(a,b;0)$-equivariant symplectomorphism that swaps the classes $B$ and $F$. We prove this by contradiction. Suppose $\psi \in\Symp^{\Z_n}(\SSS,\oml)$ satifies $\psi_*B=F$ and $\psi_*F=B$ and let $\overline{B}_{QR}$ denote the standard $B$-curve which is holomorphic with respect to the product complex structure. Then $\psi(\overline{B}_{QR})$ is in the class $F$ and needs to pass through $2$ of the $4$ isolated fixed points of the $\Z_n(a,b;0)$ action. Hence $\psi$ satisfies one of the following properties:
\begin{enumerate}
\item Either $Q$ is fixed by $\psi$ and $\psi(R)$ is either the point $P$ or $S$, or
\item the point $R$ is fixed by $\psi$ and $\psi(Q)$ is either the point $P$ or $S$.
\end{enumerate}
In either cases this contradicts Remark~\ref{rmk:FixedRQforgenericF_0}.
\end{proof}

%%%%%%%%%%%%%%%%%%%%%%%%%%%%%%%%%%%%%%%%%%%%%%%%%%%%%%%%%%%%%%%%%%%%%%%%%%%%%%%%
\subsection{Actions of type Z2}

In this section, we prove the homogeneity property for actions of the form $\Z_n(a,b,0)$ such that either $2a \equiv 0\pmod{n}$ or  $2b \equiv 0\pmod{n}$, $n\neq 2$ and $\lambda=1$.\\ 

Note that because we assume $\gcd(a,b)=1$ and $a,b \in \{0,\cdots,n-1\}$, the only effective actions for which $2a$ and $2b$ are simultaneously congruent to $0\pmod{n}$ are $\Z_{2}(1,0;0)$, $\Z_{2}(0,1;0)$, and $\Z_2(1,1;0)$, which are not of type Z2 since $n=2$. Consequently, 
\begin{enumerate}
\item either $2a$ or $2b$ is congruent to $0\pmod{n}$, but not both, and 
\item Z2 actions have isolated fixed points unless they are of the form $\Z_n(1,0;0)$ or $\Z_n(0,1;0)$, $n\neq2$.
\end{enumerate}
Note also that, by symmetry, it is enough to consider Z2 actions for which $2a\equiv0\pmod{n}$. 

In particular, $\gcd(a,n)\neq 1$. As for actions of type H0, Lemma~\ref{LR} and positivity of intersections imply the following statement.

\begin{lemma}\label{Lemma:gcd(a,n)not1-uniquecurveQR}
Consider an action $\Z_{n}(a,b;0)$ of type Z2 with  $2a \equiv 0\pmod{n}$. The $\T_0$ invariant curve in class $B$ that passes through the fixed points $Q$ and $R$ is the only $\Z_{n}(a,b;0)$ invariant curve in class $B$ passing through the fixed point $Q$.\qed
\end{lemma}

\begin{lemma}\label{Lemma:SwapRandQ}
Consider $\Z_n(a,b;0)$ action of type Z2 with $2a \equiv 0\pmod{n}$,  $2b \not\equiv 0\pmod{n}$, $n\neq 2$ and $\lambda=1$. Given any $\phi \in \Symp_h^{\Z_n}(\SSS,\oml)$, either $\phi$ fixes the two fixed points $Q$ and $R$, or it swaps them. Moreover, there exist equivariant symplectomorphisms that swap $Q$ and $R$.
\end{lemma} 
\begin{proof}
The first statement follows from comparing the  weights at the four fixed points $P$, $Q$, $R$, and $S$. 
For the second statement, observe that the symplectic involution in the Weyl group of $\T_0$ that is obtained by reflecting the moment polytope about its central vertical axis swaps $Q$ and $R$. 
\end{proof}

%%%%%%%%%%%%%%%%%%%%%%%%%%%%%%%%%%%%%%%%%%%%%%%%%%%%%%%%%%%%%%%%%%%%%%%%%%%%%%%%
\subsubsection{Actions of type Z2 with isolated fixed points}\label{Subsection:F_1actionisolatedfixedpoint}

In this subsection we calculate the homotopy type of $\Symp^{\Z_n}(\SSS,\oml)$ for Z2 actions of the form 
$\Z_{2a}(a,b;0)$ with $2a\equiv 0$, $a\not\equiv0$, $2b\not\equiv 0$ and $\lambda=1$.\\

Fix $p_0$ to be the point $Q$. Let $\Symp^{\Z_n}_{h,p_0}(\SSS,\oml)$ denote the space of $\Z_n$-equivariant symplectomorphisms that act trivially on homology and pointwise fix the point $Q$, and let $\mathcal{S}^{\Z_n}_{B,p_0}$ denote the space of $\Z_n$-invariant curves in class $B$ passing through $Q$. Then there is a sequence of fibrations:
\[\Stab^{\Z_n}(B_{p_0}) \to \Symp^{\Z_n}_{h,p_0}(\SSS,\oml) \longtwoheadrightarrow \mathcal{S}^{\Z_n}_{B,p_0} \mathbin{\textcolor{blue}{\xrightarrow{\text{~~~$\simeq$~~~}}}}  U_{0}^{\Z_n}\simeq \jzoml \simeq\{*\}\]
\[\Fix^{\Z_n}(B_{p_0}) \to \Stab^{\Z_n}_{p_0}(B_{p_0}) \longtwoheadrightarrow  \Symp^{\Z_n}(B_{p_0}) \mathbin{\textcolor{blue}{\xrightarrow{\text{~~~$\simeq$~~~}}}} S^1\]
\[\Fix^{\Z_n} (N(B_{p_0})) \to \Fix^{\Z_n}(B_{p_0}) \longtwoheadrightarrow  \Aut^{\Z_n}(N(B_{p_0})) \mathbin{\textcolor{blue}{\xrightarrow{\text{~~~$\simeq$~~~}}}} S^1 \]
\[\Stab^{\Z_n}(\overline{F}) \cap \Fix^{\Z_n}(N(B_{p_0})) \to \Fix^{\Z_n}(N(B_{p_0})) \longtwoheadrightarrow  \overline{\mathcal{S}^{\Z_n}_{F,p_0}} \mathbin{\textcolor{blue}{\xrightarrow{\text{~~~$\simeq$~~~}}}} \mathcal{J}^{\Z_n}(B_{p_0})\]
\[\Fix^{\Z_n}(\overline{F}) \to \Stab^{\Z_n}(\overline{F}) \cap \Fix^{\Z_n}(N(B_{p_0})) \longtwoheadrightarrow  \Symp^{\Z_n}(\overline{F}, N(p_0))  \mathbin{\textcolor{blue}{\xrightarrow{\text{~~~$\simeq$~~~}}}} \left\{*\right\}\]
\[\left\{*\right\} \mathbin{\textcolor{blue}{\xleftarrow{\text{~~~$\simeq$~~~}}}} \Fix^{\Z_n}(N(B_{p_0} \vee \overline{F})) \to \Fix^{\Z_n}(\overline{F}) \longtwoheadrightarrow  \Aut^{\Z_n}(N(B_{p_0} \vee \overline{F}))  \mathbin{\textcolor{blue}{\xrightarrow{\text{~~~$\simeq$~~~}}}} \left\{*\right\}\]
The proofs that these maps are fibrations are the same as in Section~\ref{Subsection:F_0action}, and the notation is described in Section~\ref{not}.
 
\begin{thm}\label{Thm:HomotopytypeF_1isolated}
Consider $\Z_n(a,b;r)$ action of type Z2 with $a\not\equiv0$ and $b\not\equiv 0$. Then $\Symp^{\Z_n}_{h}(\SSS,\oml) \simeq \T_0 \times \Z_2$
\end{thm}
\begin{proof}
From the above fibrations we can conclude that  $\Symp^{\Z_n}_{h,p_0}(\SSS,\oml)/\T_0 \simeq \jzoml \simeq\{*\}$. We now consider the fibration
  \begin{equation*}
    \Symp_{h,p_0}^{\Z_n}(\SSS,\oml)\to  \Symp_h^{\Z_n}(\SSS,\oml)\to \Symp_h^{\Z_n}(\SSS,\oml)\cdot p_0
\end{equation*}
where $\Symp_h^{\Z_n}(\SSS,\oml)\cdot p_0$  denotes the orbit of $p_0$ under the action of  $\Symp_h^{\Z_n}(\SSS,\oml)$. By Lemma~\ref{Lemma:SwapRandQ}, the orbit $\Symp_h^{\Z_n}(\SSS,\oml)\cdot p_0$ consists of two points, namely $Q$ and $R$. Hence $\Symp_h^{\Z_n}(\SSS,\oml) \simeq \Symp^{\Z_n}_{h,p_0}(\SSS,\oml) \times \Z_2 \simeq T_0 \times \Z_2$ as a topological space.
\end{proof}

In order to obtain the homotopy type of the entire symplectomorphism group we consider the following exact sequence
\begin{equation*}
1\to \Symp_h^{\Z_n}(\SSS,\oml)\to  \Symp^{\Z_n}(\SSS,\oml)\to \Aut_{c_1,\oml}(H^2(\SSS))\to 1
\end{equation*}
where $ \Aut_{c_1,\oml}(H^2(\SSS))$ denotes the space of automorphisms of $H^2(\SSS)$ that preserve the first Chern class and the symplectic form. In the case of $(\SSS,\oml)$ with $\lambda = 1$, the group $\Aut_{c_1,\oml}\left(H_2(\SSS,\Z)\right)$ is equal to $\Z_2$ and is generated by the symplectomorphism that swaps the two $S^2$ factors. By a simple comparison of weights at $P$ and $R$, for Z2 actions satisfying $a\not\equiv0$ and $b\not\equiv 0$ one can observe that the flip that generates $\Aut_{c_1,\oml}\left(H_2(\SSS,\Z)\right)$ can be realised by a $\Z_n(a,b,0)$ equivariant symplectomorphism iff $a=b=1$ and $n=2$. Hence $ \Symp^{\Z_n}(\SSS,\oml)=\Symp_h^{\Z_n}(\SSS,\oml)\simeq \T_0 \times \Z_2 $.

\begin{thm}\label{Thm:FullHomotopytypeF_1isolated}
For $\Z_n(a,b;0)$ actions of type Z2 with  $a\not\equiv0$ and $b\not\equiv 0$, the equivariant symplectomorphism group $\Symp^{\Z_n}(\SSS,\oml)$ is homotopy equivalent to $\T_0 \times \Z_2$.\qed
\end{thm}
%%%%%%%%%%%%%%%%%%%%%%%%%%%%%%%%%%%%%%%%%%%%%%%%%%%%%%%%%%%%%%%%%%%%%%%%%%%%%%%%
\subsubsection{Actions of type Z2 with fixed surfaces}
The only actions of type Z2 that have fixed surfaces are actions of the form $\Z_n(0,1;0)$ or $\Z_n(1,0;0)$ with $n \neq 2$ and $\lambda =1$. 
\begin{lemma}\label{Lemma:FixedsurfacesF_1}
Consider a $\Z_n(a,b;0)$ action of the form $\Z_n(0,1;0)$ or $\Z_n(1,0;0)$ with $n \neq 2$ and $\lambda =1$. Then the fixed point set consists of two disjoint fixed spheres.
\end{lemma}
\begin{proof}
By symmetry, we can restrict to the action $\Z_n(0,1,0)$. This is the product $\Z_n$ action on $\SSS$ that pointwise fixes $S^2 \times \{0\}$ and $S^2 \times \{\infty\}$. These are the only fixed surfaces. 
\end{proof}

\begin{lemma}\label{Lemma:SwapFixedsurfacesF_1}
Let $B_1$ and $B_2$ denote the two connected components of the fixed point set of a $\Z_n(a,b;0)$ actions with $a\equiv0$ or $b\equiv 0$ mod($n$). Given $\phi \in \Symp_h^{\Z_n}(\SSS,\oml)$, either $\phi(B_1)=B_1$ or $\phi(B_1)=B_2$. 
\end{lemma}
\begin{proof}
By Lemma~\ref{Lemma:FixedsurfacesF_1} there exists exactly two fixed surfaces. As any equivariant symplectomorphism takes a fixed surface to another fixed surface the proof would be complete if we can indeed produce an equivariant symplectomorphism that swaps the two fixed surfaces. This is realized by the equivariant symplectomorphism described in Corollary~\ref{Cor:Z_nactionspreservedbyWeylgroup} which swaps the two surfaces $B_1$ and $B_2$.
\end{proof}

Consider the fibrations
\[\Fix^{\Z_n}(B_{1}) \longrightarrow \Symp^{\Z_n}_{h,B}(\SSS,\oml)  \longtwoheadrightarrow \Symp(B_{1}) \mathbin{\textcolor{blue}{\xrightarrow{\text{~~~$\simeq$~~~}}}} SO(3)\rule{0em}{2em}\]
\[\Stab^{\Z_n}(F_{p_0}) \longrightarrow \Fix^{\Z_n}(B_{1}) \longtwoheadrightarrow \overline{\mathcal{S}^{\Z_n}_{F,p_0}} \mathbin{\textcolor{blue}{\xrightarrow{\text{~~~$\simeq$~~~}}}} \jsom \simeq \{*\}\rule{0em}{2em}\]
\[\Fix^{\Z_n} (F_{p_0}) \longrightarrow \Stab^{\Z_n}(F_{p_0}) \longtwoheadrightarrow \Symp^{\Z_n}(F_{p_0}) \mathbin{\textcolor{blue}{\xrightarrow{\text{~~~$\simeq$~~~}}}} S^1\rule{0em}{2em}\] 
\[\left\{*\right\}\mathbin{\textcolor{blue}{\xleftarrow{\text{~~~$\simeq$~~~}}}} \Fix^{\Z_n}(N(B_{max} \vee F_{p_0})) \longrightarrow  \Fix^{\Z_n}(F_{p_0})\longtwoheadrightarrow \Aut^{\Z_n}(N(B_{max} \vee \overline{F}))
\mathbin{\textcolor{blue}{\xrightarrow{\text{~~~$\simeq$~~~}}}}  \left\{*\right\}\rule{0em}{2em}\]
 
The notation used is described in Section~\ref{not}. The proofs of the above fibrations are the same as in the proof of the fibration in Section 3.5.5.2 in \cite{ChPin-Memoirs}. We refer the reader to that section for more details.
Finally, from the above fibrations and Lemma~\ref{Lemma:SwapFixedsurfacesF_1}   we obtain the following result
\begin{thm}\label{Thm:FullhomotopyF1fixedsurfaces}
Consider a $\Z_n(a,b;0)$ action such that either $a\equiv0$ or $b\equiv 0$ mod($n$) and $n\neq 2$. Then $\Symp^{\Z_n}_{h}(\SSS,\oml)= \Symp^{\Z_n}(\SSS,\oml) \simeq S^1\times \SO(3)\times \Z_2$
\end{thm}
\begin{proof}
The proof is the same as Theorem~\ref{Thm:FullHomotopytypeF_1isolated}.
\end{proof}
%%%%%%%%%%%%%%%%%%%%%%%%%%%%%%%%%%%%%%%%%%%%%%%%%%%%%%%%%%%%%%%%%%%%%%%%%%%%%%%%
\subsection{Actions of type Z3}
There are three actions of type Z3, namely the actions $\Z_2(1,1;0)$, $\Z_2(1,0;0)$ and $\Z_2(0,1,0)$. The first action has 4 isolated fixed points and the second and third actions have two fixed surfaces each.
\\

Then we have the following fibrations

\[\Stab^{\Z_n}(B_{p_0}) \to \Symp^{\Z_n}_{h,p_0}(\SSS,\oml) \longtwoheadrightarrow \mathcal{S}^{\Z_n}_{B,p_0} \mathbin{\textcolor{blue}{\xrightarrow{\text{~~~$\simeq$~~~}}}}  \jsom \cap U_{0}\]
 \[\Fix^{\Z_n}(B_{p_0}) \to \Stab^{\Z_n}_{p_0}(B_{p_0}) \longtwoheadrightarrow  \Symp^{\Z_n}(B_{p_0}) \mathbin{\textcolor{blue}{\xrightarrow{\text{~~~$\simeq$~~~}}}} S^1\]
\[\Fix^{\Z_n} (N(B_{p_0})) \to \Fix^{\Z_n}(B_{p_0}) \longtwoheadrightarrow  \Aut^{\Z_n}(N(B_{p_0})) \mathbin{\textcolor{blue}{\xrightarrow{\text{~~~$\simeq$~~~}}}} S^1 \]
\[\Stab^{\Z_n}(\overline{F}) \cap \Fix^{\Z_n}(N(B_{p_0})) \to \Fix^{\Z_n}(N(B_{p_0})) \longtwoheadrightarrow  \overline{\mathcal{S}^{\Z_n}_{F,p_0}} \mathbin{\textcolor{blue}{\xrightarrow{\text{~~~$\simeq$~~~}}}} \mathcal{J}^{\Z_n}(B_{p_0})\]
\[\Fix^{\Z_n}(\overline{F}) \to \Stab^{\Z_n}(\overline{F}) \cap \Fix^{\Z_n}(N(B_{p_0})) \longtwoheadrightarrow  \Symp^{\Z_n}(\overline{F}, N(p_0))  \mathbin{\textcolor{blue}{\xrightarrow{\text{~~~$\simeq$~~~}}}} \left\{*\right\}\]
\[\left\{*\right\} \mathbin{\textcolor{blue}{\xleftarrow{\text{~~~$\simeq$~~~}}}} \Fix^{\Z_n}(N(B_{p_0} \vee \overline{F})) \to \Fix^{\Z_n}(\overline{F}) \longtwoheadrightarrow  \Aut^{\Z_n}(N(B_{p_0} \vee \overline{F}))  \mathbin{\textcolor{blue}{\xrightarrow{\text{~~~$\simeq$~~~}}}} \left\{*\right\}\]
For the proofs that the above maps are fibrations see Section~\ref{Subsection:F_1actionisolatedfixedpoint}. Thus we conclude that $\Symp^{\Z_n}_{h,p_0}(\SSS,\oml)/\T_0 \simeq \jzoml$.
\begin{thm}
Consider the action $\Z_2(1,1;0)$ on $(\SSS,\oml)$. Then the space of equivariant symplectomorphisms $\Symp^{\Z_n}(\SSS,\oml)$ is homotopy equivalent as a topological space to $\T_0 \times \Z_8$.
\end{thm}
\begin{proof}
Consider the fibration 
\begin{equation*}
\Symp_{h,p_0}^{\Z_n}(\SSS,\oml)\to  \Symp^{\Z_n}_h(\SSS,\oml)\to \Symp^{\Z_n}_h(\SSS,\oml)\cdot p_0
\end{equation*}
where $\Symp^{\Z_n}_h(\SSS,\oml)\cdot p_0 $ denotes the orbit of $p_0$ under $\Symp^{\Z_n}_h(\SSS,\oml)$.   By point (1) in Corollary~\ref{Cor:Z_nactionspreservedbyWeylgroup}, we know that the two symplectomorphisms $\phi_1$ and $\phi_2$ given by  reflecting the moment polytope along the vertical and horizontal axes respectively, are equivariant with respect to the $\Z_2(1,1;0)$ action. Further, they preserve the homology and they act transitively on the 4 fixed points. Hence the orbit of $p_0$ under compositions of these symplectomorphisms consists of the 4 points $\{P,Q,R,S\}$ and as a consequence  $\Symp^{\Z_n}_h(\SSS,\oml) \simeq  \Symp_{h,p_0}^{\Z_n}(\SSS,\oml) \times \Z_4 \simeq \T_0 \times \Z_4$ as topological spaces. 
\\

To obtain the homotopy type of the entire symplectomorphism group we consider the fibration
\begin{equation*}
1\to \Symp_h^{\Z_n}(\SSS,\oml)\to  \Symp^{\Z_n}(\SSS,\oml)\to \Aut_{c_1,\oml}(H^2(\SSS))\to 1
\end{equation*}

We know that $\Aut_{c_1,\oml}(H^2(\SSS)) \cong \Z_2$ Further the symplectomorphism swapping the 2 $S^2$ factors in $\SSS$ is an equivariant symplectomorphism for the action $\Z_2(1,1;0)$ Hence $\Symp^{\Z_n}(\SSS,\oml) \simeq \Symp_h^{\Z_n}(\SSS,\oml) \times \Z_2 \simeq T_0 \times \Z_8$.
\end{proof}
\begin{thm}
Consider the actions $\Z_2(1,0;0)$ or $\Z_2(0,1;0)$ on $(\SSS,\oml)$. Then the space of equivariant symplectomorphisms $\Symp^{\Z_n}(\SSS,\oml)$ is homotopy equivalent to $S^1 \times \SO(3)\times \Z_2$
\end{thm}
\begin{proof}
The proof is exactly the same as in Theorem~\ref{Thm:FullhomotopyF1fixedsurfaces}.
\end{proof}

%%%%%%%%%%%%%%%%%%%%%%%%%%%%%%%%%%%%%%%%%%%%%%%%%%%%%%%%%%%%%%%%%%%%%%%%%%%%%%%%
\subsection{Actions of type Z0}
Using techniques similar to the ones in the above sections, we can prove the following theorem

\begin{thm}\label{Thm:HomtopytypeZn(a,b;0)}
Fix a $\Z_n(a,b;0)$ action satisfying the conditions $\gcd(a,n) \neq 1$ and $\lambda > 1$, then the space of equivariant symplectomorphisms has the following homotopy type
\begin{enumerate}
\item When either $a=0$ or $b=0$, then $\Symp^{\Z_n}(\SSS,\oml) \simeq S^1 \times \SO(3) \times \Z_2$
\item When $2a \not\equiv 0\pmod{n}$ and  $2b \not\equiv 0\pmod{n}$ then $\Symp^{\Z_n}(\SSS,\oml) \simeq \T_0$.
\item When either $2a \equiv 0\pmod{n}$ or  $2b \equiv 0\pmod{n}$ with $a\neq0$ and $b\neq 0$  then $\Symp^{\Z_n}(\SSS,\oml) \simeq \T_0 \times \Z_2$
\end{enumerate}   
\end{thm}
\begin{proof}
By Corollary~\ref{prop:gcdConditionSingleToricExtension}, we know that when $\gcd(a,b) \neq 1$, then the $\Z_n(a,b;0)$ extends only to the torus $\T_0$. By  Lemma~\ref{Lemma:gcd(a,n)not1-uniquecurveQR} we can prove that there is a unique curve in class $B$ passing through $Q$ and $R$. Hence the lemma follows by considering the respective fibrations as in the previous sections. 
\end{proof}

\begin{remark}
Although the above proof works for all actions of the form $\Z_n(a,b;0)$ with $\gcd(a,n) \neq 1$ and $\lambda \geq 1$, we note that the case $\Z_n(a,b;0)$ with $\gcd(a,n) \neq 1$ and $\lambda =1$ is subsumed in the actions of types Z1, Z2 and Z3.
\end{remark}

%%%%%%%%%%%%%%%%%%%%%%%%%%%%%%%%%%%%%%%%%%%%%%%%%%%%%%%%%%%%%%%%%%%%%%%%%%%%%%%%
\section{Stratification and codimension calculations for actions of type H2}
%%%%%%%%%%%%%%%%%%%%%%%%%%%%%%%%%%%%%%%%%%%%%%%%%%%%%%%%%%%%%%%%%%%%%%%%%%%%%%%%
The goal of this section is to outline the proof of the following proposition that is essential in our understanding of the centralizer group of $\Z_n(1,b;r)$ actions of type H2.
\begin{prop}
Let $\Z_n(1,b;r)$ be an action of type H2. Then the contractible space $\jzoml$ decomposes into two disjoint strata 
\[\jzoml= U_r^{\Z_n}\sqcup U_{r'}^{\Z_n}\]
one of which is open and dense while the other is a codimension 2 submanifold.
\end{prop}

We proceed as in Chapter~5 of~\cite{ChPin-Memoirs} where an analogous result is proven for circle actions that admit exactly two toric extensions. In particular, we construct the universal moduli spaces of $\Z_n$ equivariant curves, use them to model strata in $\jzoml$, and show that the normal bundle at an integrable invariant almost complex structure $J\in\jzoml$ is given by the space of $\Z_n$ invariant deformations of $J$. As the proofs to these claims are identical to the corresponding statements for circle actions, we shall only state the precise theorems below and refer the reader back to Chapter~5 of~\cite{ChPin-Memoirs} for more details.

%%%%%%%%%%%%%%%%%%%%%%%%%%%%%%%%%%%%%%%%%%%%%%%%%%%%%%%%%%%%%%%%%%%%%%%%%%%%%%%%
\subsection{Universal equivariant moduli spaces}
Consider a $\Z_n$ action of type H2 that extends to the toric actions $T_r$ and $T_{r'}$. Throughout this section, let $s=2k$ be equal to either $r$ or $r'$, and identify the action with $\Z_n(1,b;s)\subset\T_s$. By Lemma~\ref{lemma:ConditionsForDrCurvesOnlyThroughQandR} we know that every invariant curve with self intersection $-s$ in class $D_s=B-kF$ passes through the fixed points $Q$ and $R$ of the $s^{th}$ Hizerbruch surface $W_s$ as in figure~\ref{hirz}. Let $\overline{D}$ be the $\T_s$ invariant and $J_s$-holomorphic curve $\overline{D}$ in class $D_s$. 
The group $\Z_n(1,b;s)$ acts effectively on $\overline{D}$, and the same is true for all invariant $D_s$ curves. 

\begin{lemma}\label{WellDefinedModuliZn}
Consider a $\Z_n(1,b;s)$ action of type H2 acting on $(\SSS,\oml)$. Given any invariant $D_s$ curve $C$, the $\Z_n$ action induced on $C$ is effective.
\end{lemma}
\begin{proof}
Follows from Lemma \ref{weightZ_n}. 
\end{proof} 

Any invariant curve $C$ in class $D_s$ is $J$-holomorphic for some $J\in U_s^{\Z_n}$. Any $J$-holomorphic par\-a\-me\-tri\-zation $u: (S^2,j_0) \to (M,J)$ defines an effective $\Z_n$ action on $S^2$ that extends to a unique holomorphic circle action. 
As all circle subgroups of $\PSL(2,\C)$ are conjugate, there is an automorphism $\Lambda\in\PSL(2,\C)$ that brings the restricted $\Z_n$ action on $C$ to the standard action of $\Z_n$ on $S^2$.  
In other words, any invariant curve $C$ in class $D_s$ is the image of some $\Z_n$-equivariant $J$-holomorphic map $u:\big(S^2, j_0)\to \big(\SSS, J, \big)$. 

The pair $(u,J)$ belongs to the equivariant universal moduli space of class $C^l$
\begin{multline*}
\mathcal{M}^{\Z_n}(D_s , \jsoml): = \{(u,J) ~|~ \text{u is  $\Z_n$-equivariant,~ somewhere injective,~$J$-holomorphic,}\\
\text{and represents the class $D_s$}\}
\end{multline*}
where $\Z_n$ acts in the standard way on $S^2$ and where it acts as $\Z_n(1,b;s)$ on $\SSS$. Note that the image of the projection map
\[\pi: \mathcal{M}^{\Z_n}(D_s , \jzomll) \to \jzomll\]
is precisely the stratum $U_s^{\Z_n}$. The next proposition follows from the same line of arguments used in the proofs of Theorem~5.8 and Corollary~5.12 in~\cite{ChPin-Memoirs}.
\begin{prop}\label{prop:Projection map is a sub-immersion}
Fix a $\Z_n(a,b;s)$ action of type H2 on $(\SSS,\oml)$. The moduli space $\mathcal{M}^{\Z_n}(D_s , \jzomll)$ is a Banach manifold, the projection map $\pi: \mathcal{M}^{\Z_n}(D_s , \jzomll) \to \jzomll$ is a smooth sub-immersion whose image, which is equal to the stratum $U_s^{\Z_n}$, is a Banach submanifold of $\jzomll$.\qed
\end{prop}
Let $I^{\Z_n}_{\oml,l} \subset  \jzomll$ be the subspace of invariant, compatible, and complex (integrable) structures. As in the case of Hamiltonian $S^1$ actions, the normal bundle to the stratum $U_s^{\Z_n}$ at an integrable $J$ can be identified with the space of $\Z_n$-invariant deformations of $J$.
\begin{lemma}\label{trans_strata_Z_n}
Fix a $\Z_n$ action of type H2 on $M= (\SSS,\oml)$ . Further let $i: I^{\Z_n}_{\oml,l} \hookrightarrow  \jzomll$ denote the inclusion and  $\pi: \mathcal{M}^{\Z_n}(D_s , \jzomll) \to \jzomll$ denote the projection. Then,
\begin{itemize}
\item $i \pitchfork \pi$,
\item the infinitesimal complement (i.e the fibre of the normal bundle)  of $U_{s,l}^{\Z_n}$ at $J_s \in I^{\Z_n}_{\om,l} \cap U^{\Z_n}_{s}$ can be identified with $H^{0,1}_{J_s}(M,TM)^{\Z_n}$.
\end{itemize}
\end{lemma}
\begin{proof}
The proof is essentially the same as for Lemma 5.14 in \cite{ChPin-Memoirs}.
\end{proof}

% Even isotropy representations
\subsection{Even isotropy representations and codimension calculation}

In the non-equivariant case, the infinitesimal deformations of a complex structure $J\in U_s$ is identified with $H^{0,1}_{J_s}(M,TM) \cong \C^{s-1}$.  Hence, to determine the codimension of $U_s^{\Z_n}$ in $\jzoml$, it suffices to calculate the dimension of the subspace of invariant elements in ${H^{0,1}_{J_s}(M,TM)}$ under the action of~$\Z_n$. The following  Theorem tells us how the isometry group $K(s)\simeq S^{1}\times \SO(3)$ of the Kahler pair $(\oml, J_s)$ acts on the space $H_{J_s}^{0,1}(M,TM)$ of infinitesimal deformations. 
\begin{thm}[Theorem 4.2 in \cite{AGK}]\label{Thm_action on deformations} The Kahler isometry group $K(s)\simeq S^1 \times SO(3)$ of the Hirzebruch surface $W_s$ acts on $H^{0,1}(M,T^{1.0}_{J_s}M)$ via the representation $\Det\otimes\Sym^{s-2}(\C^2)$, where $\Det$ is the standard action of $S^{1}=U(1)$ on $\C$, and where $\Sym(\C^{2})$ is the representation $\mathscr{W}_{k-1}$ of $\SO(3)$ induced by the $(s-2)$-fold symmetric product of the standard representation of $\SU(2)$ on $\C^{2}$.\qed
\end{thm}

To better understand how the subgroup $\Z_n(a,b;s) \subset \T_s\subset K(s)\simeq S^1 \times \SO(3)$ acts on $H^{0,1}_{J_s}(M,TM) \cong \C^{r-1}$, denote by $R(t)$ the element of the maximal torus of $\SO(3)=\PU(2)$ given by
\[R(t)=\begin{pmatrix}1 & 0 & 0\\ 0 & \cos(t) & -\sin(t)\\ 0 & \sin(t) & \cos(t)\end{pmatrix}, ~t\in[0,2\pi)\]
and which lifts to 
\[e(t/2):=\begin{pmatrix}e^{it/2} & 0\\ 0 & e^{-it/2}\end{pmatrix}\in\SU(2).\]
 
A basis of $\Sym^{s-2}$ is given by the homogeneous polynomials $P_{m}=z_{1}^{s-2-m}z_{2}^{m}$ for $m\in\{0,\ldots,s-2\}$. The action of $R(t)$ on $P_{m}$ is
\[R(t)\cdot P_{m}=e(t/2)\cdot P_{m}=e^{i\big(s-2-2m\big)t/2}P_{m}=e^{it(k-1-m)}P_{m}\]
so that the action of $(e^{i\theta},R(t))\subset S^{1}\times\SO(3)$ on $P_{m}$ is
\[(e^{i\theta},R(t))\cdot P_{m}=e^{i\big(\theta+t(k-1-m)\big)}P_{m}\]
Each $P_{m}$ generates an eigenspace for the action of the maximal torus $\T(s)$. In particular, the finite group $\Z_n(a,b;s) \subset S^{1}(a,b;s)$ acts trivially on $P_{m}$ if, and only if, 
\[b+a(k-1-m)=(a,b)\cdot(k-1-m,1)=0\pmod{n}\]
for $m\in\{0,s-2\}$. Equivalently, we must have
\[am+b=(a,b)\cdot(m,1) \equiv 0\pmod{n} \]
for $m\in\{1-k,\ldots,k-1 \}$
\\

As explained in Section~5.3.2 of~\cite{ChPin-Memoirs}, the above calculation is done with respect to the basis of the maximal torus in $K(2n)$ obtained from the construction of the Hirzebruch surface $W_s$ via symplectic reduction. 
This basis differs from the basis that we use for the description of the moment map. The two bases are related through the matrix $\begin{pmatrix}k& -1\\1&0\end{pmatrix}$ which takes the vector $(1,b)$ in the standard basis to the vector $(k-b,1)$ in the basis for the maximal torus in $K(s)$. After this change of basis, we obtain 
\begin{thm}
Given the action $\Z_n(1,b;s)$ on $(\SSS,\oml)$ of type H2, the complex codimension of the stratum $U_s^{\Z_n}$ in $\jzoml$ is given by the number of indices $j \in \{1, \cdots , s-1\}$ such that $j \equiv b\pmod{n}$.\qed
\end{thm}

\begin{cor}
For $\Z_n(1,b;s)$ actions on $(S^2 \times S^2,\oml)$ of type H2, the complex codimension of $U_s^{\Z_n}$ in $\jzoml$ is either $0$ or $1$.
\end{cor}
\begin{proof}
For actions of type H2 we have $n>2\lambda >s$. As $b \in \{1,\cdots, n-1\}$ the number of solutions $j \in \{1, \cdots , s-1\}$ such that $j \equiv b\pmod{n}$ is always either 0 or~1. 
\end{proof}

%%%%%%%%%%%%%%%%%%%%%%%%%%%%%%%%%%%%%%%%%%%%%%%%%%%%%%%%%%%%%%%%%%%%%%%%%%%%%%%%
\section{Homotopy type of centralizers of actions of type \texorpdfstring{H2}{H2} }\label{section:Centralizers-Zn(1,b;r)}
%%%%%%%%%%%%%%%%%%%%%%%%%%%%%%%%%%%%%%%%%%%%%%%%%%%%%%%%%%%%%%%%%%%%%%%%%%%%%%%%

We shall now use the codimension calculations to calculate the homotopy type of $\Symp^{\Z_n}(S^2 \times S^2,\oml)$ for actions of type H2. We follows the arguments of \cite{ChPin-Memoirs} Chapter 4 and we refer the reader to this paper for the proofs.

By Theorem \ref{lemma:Toric-Extensions-(1,b;r)}, a $\Z_n(1,b;r)$ action of type H2 extend to exactly two toric actions $\T_r$ and $\T_{r'}$, where $r'=|2b-r|$ when either $2b-r$ or $r-2b$ is in the range $(0,2\lambda)$, or $r'=r-2b+2n$ when $r-2b+2n$ is in the range $(0,2\lambda)$. As before, let $s=2k$ be either $r$ or $r'$.

\begin{thm}\label{injz}
Consider a $\Z_n(1,b;r)$ action on $(\SSS,\oml)$ of type H2. Then the two inclusions $i:\T_r, \T_{r'} \hookrightarrow \Symp^{\Z_n} (\SSS,\oml) $ induce maps that are injective in homology with coefficients in any field $\mathbb{K}$.
\end{thm}

\begin{proof}
See Lemma 4.2 in \cite{ChPin-Memoirs}.
\end{proof}

Together with the above theorem, the Leray-Hirsch Theorem implies 
\begin{cor}\label{z}
Consider $\Z_n(1,b;r)$ actions on $(\SSS,\oml)$ of type H2. Then 
\[H^*(\Symp_h^{\Z_n}(S^2 \times S^2,\oml), \mathbb{K}) \cong H^*(U_{s}^{\Z_n}, \mathbb{K}) \otimes H^*(\T_s,\mathbb{K})\]
for $s=r$ and $s=r'$.
\end{cor}

\begin{thm}\label{Thm:Z_nrankofhomology}
Consider a $\Z_n(1,b;r)$ actions on $(\SSS,\oml)$ of type H2. Then,
the cohomology groups of $\Symp_h^{\Z_n}(S^2 \times S^2,\oml)$ with coefficients in a field $\mathbb{K}$ are given by  
\[H^p\left(\Symp^{\Z_n}_h(\SSS,\oml), \mathbb{K}\right) = \begin{cases}
\mathbb{K}^4 &p \geq 2\\
\mathbb{K}^3 &p =1 \\
\mathbb{K} &p=0\\
\end{cases}\]
\end{thm}

\begin{proof}
The proof is exactly the same as the proof of Theorem 4.3 in \cite{ChPin-Memoirs}.
\end{proof}

\begin{remark}
As $\lambda >1$ it can be argued as in Theorem 4.1 in \cite{ChPin-Memoirs} that $\Symp_h^{\Z_n}(S^2 \times S^2,\oml) = \Symp^{\Z_n}(S^2 \times S^2,\oml)$ and hence we get the ranks of the cohomology of the entire symplectomorphism group.
\end{remark}

Geometrically, the two tori $\T_r$ and $\T_{r'}$ intersect in $\Symp^{\Z_n}_h(\SSS,\oml)$ along the common circle $S^1(1,b;r)=S(1,b';r')$ where $b'=b$ when either $2b-r$ or $r-2b$ is in the range $(0,2\lambda)$, and $b'=n-b$ when $r-2b+2n$ is in the range $(0,2\lambda)$.
\[
\begin{tikzcd}
S^{1} \arrow{r}{(1,b')} \arrow[swap]{d}{(1,b)} & T^{2}_{r'} \\
T^{2}_{r}  & 
\end{tikzcd}
\]
As explained in Chapter 4 in \cite{ChPin-Memoirs}, we can consider the pushout
\[P:=\pushout(T_r\from S^1\to T_{r'})\]
of the above diagram. This pushout is to be understood in the category of topological groups. The proof of Theorem~4.7 in~\cite{ChPin-Memoirs} applies mutatis mutandis to the present situation to show that the topological group $P$ is homotopically equivalent of the centralizer $\Symp^{\Z_n}_h(\SSS,\oml)$ as topological groups. Moreover,
\begin{thm}\label{thm:HomotopyTypeStabilizerCyclic-SSS-2strata}
Given a $\Z_n(a,b;r) $  on $(S^2 \times S^2,\oml)$ of type H2, then $\Symp^{\Z_n}(S^2 \times S^2,\oml) \simeq \Omega S^3 \times S^1 \times S^1 \times S^1$ as topological spaces, where $\Omega S^3$ denotes the  based loop space of $S^3$.
\end{thm}
\begin{proof}
Similar to the proof of  Theorem 4.2 in \cite{ChPin-Memoirs}.
\end{proof}

%%%%%%%%%%%%%%%%%%%%%%%%%%%%%%%%%%%%%%%%%%%%%%%%%%%%%%%%%%%%%%%%%%%%%%%%%%%%%%%%
\section{Centralizers of Hamiltonian \texorpdfstring{$\Z_n$}{Z\_n} actions on \texorpdfstring{$(\CCC,\oml)$}{CP2\#CP2}}
\label{section:CCC}
%%%%%%%%%%%%%%%%%%%%%%%%%%%%%%%%%%%%%%%%%%%%%%%%%%%%%%%%%%%%%%%%%%%%%%%%%%%%%%%%

In this section we present our results on the homotopy type of the centralizers of Hamiltonian $\Z_n$ actions on the non-trivial bundle $\CCC$. As in the case of $\SSS$, it is shown in \cite{Liat} that every such action extends to Hamiltonian circle actions. The proofs follow from simple modifications of the arguments given in Chapter 5 and Section 6.1 in \cite{ChPin-Memoirs}. We also note that Theorems \ref{Chen}, lemma \ref{lemma:Toric-Extensions-(1,b;r)} all hold even in the case when $r$ is odd.

Recall that  the  odd Hirzebruch surface $W_r$, $r=2k+1$, is defined as a complex submanifold of $ \mathbb{C}P^1 \times \mathbb{C}P^2$ defined by setting
\[
W_r:=  \left\{ \left(\left[x_1,x_2\right],\left[y_1,y_2, y_3\right]\right) \in \mathbb{C}P^1 \times \mathbb{C}P^2 ~|~  x^r_1y_2 - x_2y^r_1 = 0 \right\}
\]
This manifold is diffeomorphic to $\CP^2 \# \overline{\CP^2}$ and is invariant under the toric action defined by 
\[
 \left(u,v\right) \cdot  \left(\left[x_1,x_2\right],\left[y_1,y_2, y_3\right]\right) = \left(\left[ux_1,x_2\right],\left[u^ry_1,y_2,vy_3\right]\right)
~~(u,v) \in \T \]
and whose moment map image is
\[
\begin{tikzpicture}
\node[left] at (0,2) {$Q=(0,1)$};
\node[left] at (0,0) {$P=(0,0)$};
\node[right] at (4,2) {$R= (\lambda - \frac{r+1}{2} ,1)$};
\node[right] at (6,0) {$S=(\lambda + \frac{r-1}{2} ,0)$};
\node[above] at (2,2) {$B-\frac{r+1}{2}F$};
\node[right] at (5,1) {$F$};
\node[left] at (0,1) {$F$};
\node[below] at (3,0) {$B+ \frac{r-1}{2}F$};
\draw (0,2) -- (4,2) ;
\draw (0,0) -- (0,2) ;
\draw (0,0) -- (6,0) ;
\draw (4,2) -- (6,0) ;
\end{tikzpicture}\label{oddHirz}
\]
The homology class $B=[L]$ is the class of a line $L$ in $\CCC$, $E$ is the class of the exceptional divisor, and $F$ refers to the class $B - E$. There is a canonical form which we also call $\oml$ on $\CP^2 \# \overline{\CP^2}$, which has weight $\lambda$ on B and 1 on $F$. Note that with our convention, $\lambda$ must be strictly greater than 1 for the curve in class $E$ to have positive symplectic area.
\\

\begin{lemma}\label{lemma:Odd-Toric-Extensions-(1,b;r)}
Consider a Hamiltonian $\Z_n(1,b;r)$ action on $(\CCC,\oml)$ with $n\geq 2\lambda> r+1$.  
Then, either
\begin{enumerate}
\item $(2b-r+1)$, $(r-2b+1)$, and $(r-2b+2n+1)$ are outside the symplectic range $(0,2\lambda)$ and $\Z_{n}(1,b;r)$ only extends to the toric action $\T_r$, or
\item there is exactly one index $r'$ among $(2b-r)+1$, $(r-2b)+1$, or $(r-2b+2n)+1$ in the range $(0,2\lambda)$. In this case, $\Z_{n}(1,b;r)$ extends to $\T_r$ and $\T_{r'}$, and to no other toric actions.
\end{enumerate}
\end{lemma}
\begin{proof}
The proof is exactly the same as in Lemma~\ref{lemma:Toric-Extensions-(1,b;r)}
\end{proof}

\begin{cor}\label{cor:OddInteresectionwith2strata}
Consider a Hamiltonian $\Z_n(1,b;r)$ action on $(\CCC,\oml)$ with $n\geq 2\lambda> r+1$. Then the $\Z_n(1,b;r)$ actions extends to only one torus if $2\lambda < \min\{|r-2b|+1,\, r-2b+2n +1\}$. If instead $2\lambda > |2b-r|+1$ then the action extends to the two tori $\T_r$ and $\T_{|r-2b|}$. Finally if $2\lambda > r-2b+2n +1$ then the action extends to the two tori $\T_r$ and $\T_{r-2b+2n}$.\qed
\end{cor}

\begin{defn}
An action of the form $\Z_n(1,b;r)$ on the odd Hirzerbruch surfaces is said to be of type OE1  if it only extends to the toric action $\T_{r}$. We say that it is of type OE2 if it extends to exactly two toric actions $\T_{r}$ and~$\T_{r'}$.
\end{defn}

As in the case of $\SSS$, one can prove the homogeneity property for a subset of $\Z_n(a,b;r)$ actions on odd Hirzerbruch surfaces.
\begin{enumerate}[label=(OH{\arabic*}), start=0]
  \item Actions $\Z_n(a,b;r)$ such that $\gcd(a,n) \neq 1$ as considered in Proposition~\ref{prop:gcdConditionSingleToricExtension}
\item Actions $\Z_n(1,b;r)$ of type OE1 satisfying $2b-r\neq n$ and $b\not\in\{0,r\}$
\item Actions $\Z_n(1,b;r)$ of type OE2 satisfying $2b-r\neq n$ with the exclusion of $\Z_{2r}(1,b;r)$, $\Z_{2r+2}(1,r+1;r)$, and $\Z_{2r+2}(1,2r+1;r)$. 
\item Actions of type OE1 with $b\in\{0,r\}$  that is, $\Z_n(1,0;r)$ and $\Z_n(1,r;r)$ with $n\geq 2\lambda>r\geq2$ .
\end{enumerate}

Throughout the rest of the section we assume that the type of the action $\Z_n(a,b;r)$ on $(\CCC,\oml)$ is one of the above. By \cite{Liat} we know that every such action extends to a circle action and hence to a $\T$ action. The proofs of Proposition~\ref{prop:ExistenceOfGoodFixedPointH1H2} and Lemma~\ref{lemma:EvaluationFibrationConfigurations} hold verbatim for actions on the non-trivial bundle as well. As in the case of $\SSS$, we consider the  following sequence of fibrations. 
\[\Stab^{\Z_n}(\overline{D}) \to \Symp^{\Z_n}_{h}(\CCC,\oml) \longtwoheadrightarrow \mathcal{S}^{\Z_n}_{D_{2k+1}} \mathbin{\textcolor{blue}{\xrightarrow{\text{~~~$\simeq$~~~}}}} \jzoml \cap U_{2k+1}\]
\[\Fix^{\Z_n}(\overline{D}) \to \Stab^{\Z_n}(\overline{D}) \longtwoheadrightarrow  \Symp_0^{\Z_n}(\overline{D}) \mathbin{\textcolor{blue}{\xrightarrow{\text{~~~$\simeq$~~~}}}} S^1 ~\text{or}~ SO(3)\]
\[\Fix^{\Z_n} (N(\overline{D})) \to \Fix^{\Z_n}(\overline{D}) \longtwoheadrightarrow  \Aut^{\Z_n}(N(\overline{D})) \mathbin{\textcolor{blue}{\xrightarrow{\text{~~~$\simeq$~~~}}}} S^1\]
\[\Stab^{\Z_n}(\overline{F}) \cap \Fix^{\Z_n}(N(\overline{D})) \to \Fix^{\Z_n}(N(\overline{D})) \longtwoheadrightarrow  \overline{\mathcal{S}^{\Z_n}_{F,p_0}} \mathbin{\textcolor{blue}{\xrightarrow{\text{~~~$\simeq$~~~}}}} \mathcal{J}^{\Z_n}(\overline{D})\simeq \{*\} \]
\[\Fix^{\Z_n}(\overline{F}) \to \Stab^{\Z_n}(\overline{F}) \cap \Fix^{\Z_n}(N(\overline{D})) \longtwoheadrightarrow  \Symp^{\Z_n}(\overline{F}, N(p_0))  \mathbin{\textcolor{blue}{\xrightarrow{\text{~~~$\simeq$~~~}}}} \left\{*\right\}\]
\[\left\{*\right\} \mathbin{\textcolor{blue}{\xleftarrow{\text{~~~$\simeq$~~~}}}} \Fix^{\Z_n}(N(\overline{D} \vee \overline{F})) \to \Fix^{\Z_n}(\overline{F}) \longtwoheadrightarrow  \Aut^{\Z_n}(N(\overline{D} \vee \overline{F}))  \mathbin{\textcolor{blue}{\xrightarrow{\text{~~~$\simeq$~~~}}}} \left\{*\right\}\]
where $\Symp_0^{\Z_n}(\overline{D})$ in the second fibration denotes the identity component of the equivariant symplectomorphism group of $\overline{D}$. As shown in Lemma~\ref{EquivariantSO(3)}, we see that $\Symp^{\Z_n}(\overline{D})$ is homotopic to $S^{1}\times \Z_{2}$ for actions of the type $\Z_{2a}(a,b;r)$. However, comparing the weights at $Q$ and $R$ in fig~\ref{oddHirz}, we see that there is no equivariant symplectomorphism of $\CCC$ that can swap the fixed points $Q$ and $R$. Hence the space $\Stab^{\Z_n}(\overline{D})$ acts transitively only on the connected space $\Symp_0^{\Z_n}(\overline{D})$. This is in sharp contrast to $\SSS$ case where the element $c_2$ in Example~\ref{example:Type c2 is Kahler when n=2a} was an equivariant symplectomorphism for actions of the form $\Z_{2a}(a,b;r)$. 
\\

For the proofs that the above maps are fibrations, we invite the reader to consult Chapter~6 of~\cite{ChPin-Memoirs}.   
\begin{prop}\label{ev_homoCCC}
For a Hamiltonian action $\Z_n(a,b;r)$ of type OH0, the symplectic centralizer group $\Symp^{\Z_n}(\CCC,\oml)$ is homotopy equivalent to
\begin{enumerate}
\item $U(2)$ if $a=0$,
\item $\T_r$ in all other cases.
\end{enumerate}
\end{prop}
\begin{proof}
The proof is exactly the same as in the proof of Proposition~\ref{prop:HomotopyTypeH0Actions} and Theorem 6.1 in \cite{ChPin-Memoirs}.
\end{proof}

\begin{prop}\label{prop:OddHomotopyTypeH1Actions}
For a Hamiltonian action $\Z_n(1,b;r)$ of type OH1 or OH3, the symplectic centralizer group $\Symp^{\Z_n}(\CCC,\oml)$ is homotopy equivalent to the torus $\T_r$.
\end{prop}
\begin{proof}
The proof is exactly the same as in Proposition~\ref{prop:HomotopyTypeH1Actions}.
\end{proof}

\begin{thm}
Consider a $\Z_n(a,b;r)$ action on $\CCC$ of type OH2. Then the complex codimension of the strata $U_r^{\Z_n}$ in $\jzoml$ is given by the number of indices $m \in \{1, \cdots , r-1\}$ such that $m \equiv b \pmod{n}$.
\end{thm}
\begin{proof}
The codimension calculation follows from the following argument. Let $s=2k+1$ be equal to either $r$ or $r'$.

As established in Section 6.2.1 in \cite{ChPin-Memoirs}, the action of the isometry group $K(s)\simeq \U(2)$ on infinitesimal deformations is isomorphic to $\Det^{1-k}\otimes\Sym^{s-2}$. A basis of $\Sym^{s-2}$ is given by the homogeneous polynomials $P_{m}=z_{1}^{s-1-m}z_{2}^{m}$ for $m\in\{0,\ldots,s-2\}$. The action of $D_{\theta,t}$ on $P_{m}$ is
\[D_{\theta,t}\cdot P_{m}=e^{i\big((\theta+t)(1-k)+\theta(s-2-m)+tm\big)}P_{m}\]
so that each $P_{m}$ generates an eigenspace for the action of the maximal torus $T(s)$ generated by $D_{\theta,t}$. In particular, the circle $S^{1}(a,b;s)$ acts trivially on $P_{m}$ if, and only if, 
\[(a-b)(k-m)+b=(a,b)\cdot(n-m,m-k+1)\equiv 0 \pmod{n}\]
Thus the codimension in the balanced basis of the maximal torus of $K(s)$ is given by the number of $m \in \{0,\ldots,s-2\}$ such that $(a-b)(k-m)+b=(a,b)\cdot(k-m,m-k+1)\equiv 0\pmod{n}$. To express this result in the moment map basis, 
we apply the transformation given by the matrix $\begin{pmatrix}k+1& -1\\ k&-1\end{pmatrix}$ which takes the vector $(1,b)$ in the basis for the standard moment polytope to the vector $(k-b+1,\, k-b)$ in the basis of the maximal torus obtained by symplectic reduction. Hence the coefficients $a$ and $b$ in the formula above have to be replaced by $k-b$ and $k-b-1$ to get the correct codimension for the stratum $U_s^{\Z_n}$. After making this substitution, we get that this codimension is given by the number of indices $j \in \{1, \cdots , s-1\}$ such that $j \equiv b\pmod{n}$.
\end{proof}

As before, we have 
\begin{thm}\label{Oddinjz}
Consider the following $\Z_n(1,b;r)$ actions on $(\CCC,\oml)$ of type OH2.
\begin{enumerate}
\item If $r-2b+1$ or $b-2r+1$ is in the range $(0,2\lambda)$, then $r'=|r-2b|$. 
\item If instead $r-2b+2n+1$ is in the range $(0,2\lambda)$, then $r'=r-2b+2n$.
In all cases, the inclusions $i:\T_r, \T_{r'} \hookrightarrow \Symp^{\Z_n} (\CCC,\oml) $ induce maps that are injective in homology with coefficients in any field $\mathbb{K}$.
\end{enumerate}
\end{thm}

\begin{proof}
The proof is similar to the proof of Theorem~\ref{injz}.
\end{proof} 

\begin{thm}
Consider a $\Z_n(1,b;r)$ action on $(\CCC,\oml)$ with $\lambda>1$ of type OH2. Then,
the cohomology groups of $\Symp^{\Z_n}(\CCC,\oml)$ are given by  
\[H^p\left(\Symp^{\Z_n}(\CCC,\oml), \mathbb{K}\right) = \begin{cases}
\mathbb{K}^4 &p \geq 2\\
\mathbb{K}^3 &p =1 \\
\mathbb{K} &p=0\\
\end{cases}\]
\end{thm}
\begin{proof}
Same proof as in Theorem~\ref{Thm:Z_nrankofhomology}.
\end{proof}

Arguing similar to Theorem~\ref{thm:HomotopyTypeStabilizerCyclic-SSS-2strata} we have,

\begin{thm}\label{thm:OddHomotopyTypeStabilizerCyclic-CCC-2strata}
Given a $\Z_n(a,b;r) $  on $(\CCC,\oml)$ of type OH2, then $\Symp^{\Z_n}(\CCC,\oml) \simeq \Omega S^3 \times S^1 \times S^1 \times S^1$ as topological spaces, where $\Omega S^3$ denotes the  based loop space of $S^3$.
\end{thm}

%%%%%%%%%%%%%%%%%%%%%%%%%%%%%%%%%%%%%%%%%%%%%%%%%%%%%%%%%%%%%%%%%%%%%%%%%%%%%%%%
\section{Summary of results}\label{section:Summary}
%%%%%%%%%%%%%%%%%%%%%%%%%%%%%%%%%%%%%%%%%%%%%%%%%%%%%%%%%%%%%%%%%%%%%%%%%%%%%%%%
The following $\Z_n(a,b;r)$ actions with $r\neq0$ on $(\SSS,\oml)$ are tractable:
\begin{enumerate}[label=(H{\arabic*}), start=0]
\item Actions $\Z_n(a,b;r)$ such that $\gcd(a,n) \neq 1$ as considered in Proposition~\ref{prop:gcdConditionSingleToricExtension} and satisfying the extra condition $r \neq 0$.
\item Actions $\Z_n(1,b;r)$ of type E1 satisfying $2b-r\neq n$ and $b\not\in\{0,r\}$.
\item Actions $\Z_n(1,b;r)$ of type E2 satisfying $2b-r\neq n$ and $b\not\in\{0,r\}$, with the exclusion of $\Z_{2r}(1,b;r)$, $\Z_{2r+2}(1,r+1;r)$, and $\Z_{2r+2}(1,2r+1;r)$. 
\item Actions of type E1 with $b\in\{0,r\}$, that is, $\Z_n(1,0;r)$ and $\Z_n(1,r;r)$ with $n\geq 2\lambda>r\geq2$.
\end{enumerate}
The homotopy type of the symplectic centralizers of actions of type H0, H1, H2, and H3 are given in the table bellow.
\noindent
\begin{table}[H]
\begin{tabular}{|p{69mm}|p{18mm}|p{61mm}|}
\hline
$\Z_n(a,b;r)$ action on $(\SSS,\oml)$, $r\neq0$  &\# of strata & $\Symp^{\Z_n}(\SSS,\oml)$\\
\hline
\hline
Type H0 with $a \neq 0$ and $n\neq 2a$& 1 & $\simeq\T$ \\
\hline
Type H0 with $a \neq 0$ and $n=2a$& 1 & $\simeq\T \times \Z_2$ \\
\hline
Type H0 with $a =0$ &1 & $\simeq S^1\times \SO(3)$\\
\hline
Type H1 with $a \neq 0$&1 & $\simeq\T$\\
\hline
Type H1 with $a =0$ &1 & $\simeq S^1\times \SO(3)$\\
\hline
Type H2&2& $\simeq\Omega S^3 \times S^1 \times S^1 \times S^1$\\
\hline
Type H3& 1& $\simeq \T$\\
\hline
\end{tabular}
\caption{Tractable actions $\Z_n(a,b;r)$ on $(\SSS,\oml)$ with $r\neq0$}\label{table:Actions on SSS r not 0}
\end{table}
\begin{remark}
As explained in the introduction, for actions of the form $\Z_n(a,b;r)$ with $\gcd(a,n)=1$, one can reparametrize the action to bring it to an action of the form $\Z_n(1,b';r)$. If this action is of type either H1 or H2, the above table can be used to determine the homotopy type of its centralizer. 
\end{remark}

In the case $r=0$, the following $\Z_n(a,b;0)$ actions on $(\SSS,\oml)$ are tractable and the homotopy type of their centralizers are listed in the table below.

\begin{enumerate}[label=(Z{\arabic*}), start=0]
\item $\Z_n(a,b;0)$ such that $\gcd(a,n) \neq 1$ and $\lambda > 1$.
\item $\Z_n(a,b;0)$ actions such that $2a \not\equiv 0\pmod{n}$ and  $2b \not\equiv 0\pmod{n}$ with $\lambda =1$.
\item $\Z_n(a,b;0)$ actions such either  $2a \equiv 0\pmod{n}$ or  $2b \equiv 0\pmod{n}$, with $n\neq 2$ and $\lambda =1$.
\item $\Z_2(a,b;0)$ with $\lambda=1$.
\end{enumerate}

\noindent
\begin{table}[H]
\begin{tabular}{|p{76mm}|p{11mm}|p{61mm}|}
\hline
$\Z_n(a,b;0)$ action on $(\SSS,\oml)$, $r=0$ & $~\lambda$ & $\Symp^{\Z_n}(S^2 \times S^2,\oml)$\\
\hline
\hline
Type Z1 with $(a,b) \not\in \left\{(1,1),(1,n-1),(n-1,1)\right\}$ & $~\lambda = 1$ &$\simeq \T$  \\
\hline
Type Z1 with $(a,b) \in \left\{(1,1),(1,n-1),(n-1,1)\right\}$ & $~\lambda = 1$ &$\simeq \T \times \Z_2$ \\
\hline
Type Z2 with $a \not\equiv 0$ and $b \not\equiv 0$ & $~\lambda = 1$  &$\simeq \T \times \Z_2$\\
\hline
Type Z2 with either $a \equiv 0$ or $b\equiv 0$ &$~\lambda = 1$ & $\simeq S^1 \times \SO(3) \times \Z_2$\\
\hline
$\Z_2(1,1;0)$ &$~\lambda = 1$ &$\simeq \T \times \Z_8$ \\
\hline
$\Z_2(1,0;0)$ or $\Z_2(0,1;0)$ &$~\lambda=1$ &$\simeq S^1 \times \SO(3) \times \Z_2$\\
\hline
Type Z0 with either $a=0$ or $b=0$ &$~\lambda > 1$ &$\simeq S^1 \times \SO(3) \times \Z_2$ \\
\hline
Type Z0 with $2a\not\equiv0$ or $2b\not\equiv0$ &$~\lambda > 1$ &$\simeq \T$ \\
\hline
Type Z0 with either $2a\equiv0$ or $2b\equiv0$ $a\neq0$, $b\neq0$ &$~\lambda > 1$ &$\simeq \T \times \Z_2$\\
\hline
\end{tabular}
\caption{Tractable actions $\Z_n(a,b;0)$ on $(\SSS,\oml)$}\label{table:Actions on SSS r=0}
\end{table}

The following actions on the non-trivial bundle $(\CCC,\oml)$ are tractable.
\begin{enumerate}[label=(OH{\arabic*}), start=0]
  \item Actions $\Z_n(a,b;r)$ such that $\gcd(a,n) \neq 1$ as considered in Proposition~\ref{prop:gcdConditionSingleToricExtension}
\item Actions $\Z_n(1,b;r)$ of type OE1 satisfying $2b-r\neq n$ and $b\not\in\{0,r\}$
\item Actions $\Z_n(1,b;r)$ of type OE2 satisfying $2b-r\neq n$ with the exclusion of $\Z_{2r}(1,b;r)$, $\Z_{2r+2}(1,r+1;r)$, and $\Z_{2r+2}(1,2r+1;r)$. 
\item Actions of type OE1 with $b\in\{0,r\}$  that is, $\Z_n(1,0;r)$ and $\Z_n(1,r;r)$ with $n\geq 2\lambda>r\geq2$ .
\end{enumerate}

The homotopy types of the corresponding centralizers are given in the table below.

\noindent
\begin{center}
\begin{table}[H]
\begin{tabular}{|p{65mm}|p{17mm}|p{66mm}|}
\hline
Action $\Z_n(a,b;r)$ on $(\CCC,\oml)$  &\# of strata & $\Symp^{\Z_n}(\CCC,\oml)$\\
\hline
\hline
Type OH0 with $a \neq 0$ & 1 & $\simeq \T$ \\
\hline
Type OH0 with $a =0$ &1 & $\simeq U(2)$\\
\hline
Type OH1 with $a \neq 0$&1 & $\simeq \T$\\
\hline
Type OH1 with $a =0$ &1 & $\simeq U(2)$\\
\hline
Type OH2&2& $\simeq \Omega S^3 \times S^1 \times S^1 \times S^1$\\
\hline
Type OH3 & 1 & $\simeq\T$\\
\hline
\end{tabular}
\caption{Tractable actions $\Z_n(a,b;r)$ on $(\CCC,\oml)$}\label{table:Actions on CCC}
\end{table}
\end{center}

%%%%%%%%%%%%%%%%%%%%%%%%%%%%%%%%%%%%%%%%%%%%%%%%%%%%%%%%%%%%%%%%%%%%%%%%%%%%%%%%
\subsection{Discussion on remaining cases}
Because there is no complete classification of $\Z_n(a,b;r)$ actions up to equivariant symplectomorphisms on $(\SSS,\oml)$ and $(\CCC,\oml)$, we are presently unable to determine the homotopy type of the centralizer $\Symp^{\Z_n}(\SSS,\oml)$ in the following (non mutually exclusive) cases: 

\subsubsection{Actions on $(\SSS,\oml)$}\label{section:RemainingSSS}
\begin{enumerate}
\item Actions of the form $\Z_{2r}(1,b;r)$, $\Z_{2r+2}(1,r+1;r)$, or $\Z_{2r+2}(1,2r+1;r)$ which are omitted in type H2.
\item Actions of type E1 and E2 with $2b-r=n$
\item Actions of the form $\Z_n(1,b;r)$ with $n< 2\lambda$.
\item Actions of the form $\Z_n(1,b;r)$ with $2b-r=0$.  In particular, we are unable to deduce the homotopy type of equivariant symplectomorphism when $\Z_n(1,b;r)$ action extends to the tori $\T_r$ and $\T_0$. 

\item Actions of the form $\Z_n(a,b;0)$ with $\gcd(a,n)=1$ and $\lambda >1$.
\end{enumerate}

%%%%%%%%%%%%%%%%%%%%%%%%%%%%%%%%%%%%%%%%%%%%%%%%%%%%%%%%%%%%%%%%%%%%%%%%%%%%%%%%
\subsubsection{Actions on $(\CCC,\oml)$}\label{section:RemainingCCC}
\begin{enumerate}
\item Actions $\Z_{2r}(1,b;r)$, $\Z_{2r+2}(1,r+1;r)$, and $\Z_{2r+2}(1,2r+1;r)$ which are omitted in type OH2.
\item Actions of type $OE_1$ and $OE_2$ with $2b-r=n$.
\item Action of the form $\Z_n(1,b;r)$ with  $n< 2\lambda$.  
\end{enumerate}

%%%%%%%%%%%%%%%%%%%%%%%%%%%%%%%%%%%%%%%%%%%%%%%%%%%%%%%%%%%%%%%%%%%%%%%%%%%%%%%%
\section{Notation}
%%%%%%%%%%%%%%%%%%%%%%%%%%%%%%%%%%%%%%%%%%%%%%%%%%%%%%%%%%%%%%%%%%%%%%%%%%%%%%%%
 Consider a Hamiltonian $\Z_n$ action on $(M,\oml)$ and let $p_0$ be a fixed point for the group action. Given a $\Z_n$ invariant symplectic submanifold $C$, and a $\omega_\lambda$-orthogonal invariant sphere $\overline{F}$ in the homology class $F$ that intersects $C$ at a point $p_0$, we define the following spaces: 

\begin{itemize}\label{not}
\item $N(C)$:= The symplectic normal bundle to a symplectic submanifold $C$. 
\item $\Symp^{\Z_n}_h(M,\oml)$ := The group of $\Z_n$ equivariant symplectomorphisms on $(M,\oml)$ that acts trivially on homology.
\item $\Stab^{\Z_n}(C)$ := The group of all  $\phi \in \Symp^{\Z_n}_h(M,\oml)$ such that $\phi(C) = C$, that is, such that $\phi$ \emph{stabilises} $C$ but does not necessarily act as the identity on $C$.
\item $\Fix^{\Z_n}(C)$ := The group of all  $\phi \in \Symp^{\Z_n}_h(M,\oml)$ such that $\phi|_C = id$, that is, such that \emph{fixes $C$ pointwise}. 
\item $\Fix^{\Z_n}(N(C))$:= The group of all $\phi \in \Symp^{\Z_n}_h(M,\oml)$ such that $\phi|_C = \id$ and  $d\phi|_{N(C)}: N(C) \to N(C)$ is the identity on $N(C)$.
\item $\Aut^{\Z_n}(N(C))$:= The group of $\Z_n$-equivariant fiberwise symplectic automorphisms of the symplectic normal bundle of $C$.
\item $\Aut^{\Z_n}(N(C \vee \overline{F}))$ := The group of $\Z_n$-equivariant fiberwise symplectic  automorphism of the symplectic normal bundle of the crossing divisor $C \vee \overline{F}$ that are identity in a neighbourhood of the wedge point. 
\item $\mathcal{S}^{\Z_n}_{K}$ := The space of unparametrized $\Z_n$-invariant symplectic embedded spheres in the homology class $K$. 
\item $\mathcal{S}^{\Z_n}_{K,p_0}$:= The space of unparametrized $\Z_n$-invariant symplectic embedded spheres in the homology class $K$ passing through $p_0$.
\item $\J_{\oml}^{\Z_n}(C)$ := The space of $\Z_n$-equivariant $\oml$ compatible almost complex structures s.t the curve $C$ is holomorphic.
\item $\Symp^{\Z_n}(C)$:= The space of all $\Z_n$-equivariant symplectomorphisms of the curve C.
\item $\Fix^{\Z_n}(N(C \vee {\overline{F}}))$ := The space of all $\Z_n$-equivariant symplectomorphisms that are the identity in the neighbourhood of $C \vee {\overline{F}}$.
\item $\Symp^{\Z_n}({\overline{F}}, N(p_0))$ := equivariant symplectomorphism of the sphere $\overline{F}$ that are the identity in an open set of $\overline{F}$ around $p_0$. 
\item $\overline{\mathcal{S}^{\Z_n}_{F,p_0}}$:= The space of unparametrized $\Z_n$-invariant symplectic spheres in the homology class $F$ that are equal to a fixed curve ${\overline{F}}$ in a neighbourhood of $p_0$.
\item $\Symp^{\Z_n}_{p_0,h}(M,\oml)$:=  The group  of all  $\phi \in \Symp^{\Z_n}_{h}(M,\oml)$ fixing $p_0$.
\item $\Stab^{\Z_n}_{p_0}(C)$:=  The group of all  $\phi \in \Stab^{\Z_n}(C)$ such that $\phi(p_0) = p_0$. \end{itemize}
All the above spaces are equipped with the $C^\infty$ topology.

%%%%%%%%%%%%%%%%%%%%%%%%%%%%%%%%%%%%%%%%%%%%%%%%%%%%%%%%%%%%%%%%%%%%%%%%%%%%%%%%
% References
%%%%%%%%%%%%%%%%%%%%%%%%%%%%%%%%%%%%%%%%%%%%%%%%%%%%%%%%%%%%%%%%%%%%%%%%%%%%%%%%

\bibliographystyle{amsalpha}
\bibliography{bibliography}

%%%%%%%%%%%%%%%%%%%%%%%%%%%%%%%%%%%%%%%%%%%%%%%%%%%%%%%%%%%%%%%%%%%%%%%%%%%%%%%%

\end{document}